\makeatletter
\documentclass[12pt,a4paper]{article}
\usepackage[fleqn]{amsmath}
\usepackage{amscd,amssymb}
\usepackage{amsthm}
{
\theoremstyle{plain}
\newtheorem{th@@r@m}{Theorem}[section]
\newtheorem{l@mm@}[th@@r@m]{Lemma}
\newtheorem{pr@p@s@t@@n}[th@@r@m]{Proposition}
\newtheorem{c@r@ll@ry}[th@@r@m]{Corollary}
}
{
\theoremstyle{definition}
\newtheorem{d@f@n@t@@n}[th@@r@m]{Definition}
\newtheorem{r@m@rk}[th@@r@m]{Remark}
\newtheorem{@x@mpl@}[th@@r@m]{Example}
\newtheorem{pr@bl@m}[th@@r@m]{Problem}
}
\newif\ifInsection
\renewcommand{\section}{\@startsection {section}{1}{\z@}%
                                   {-28pt \@plus -3pt \@minus -1pt}%
                                   {16pt \@plus1pt}%
                                   {\normalfont\large\bfseries}}
\newdimen\setsub@nh@b@
\def\@seccntformat#1{%
\hbox to \setsub@nh@b@{\S\csname the#1\endcsname.}\ignorespaces
}
{
\theoremstyle{plain}
\newtheorem*{MTh@@r@m}{Main Theorem}
\newtheorem*{TH@@r@m}{Theorem}
}
\newenvironment{thm}{\Inthmtrue\begin{th@@r@m}}{\end{th@@r@m}\Inthmfalse}
\newenvironment{lem}{\Inthmtrue\begin{l@mm@}}{\end{l@mm@}\Inthmfalse}
\newenvironment{prop}{\Inthmtrue\begin{pr@p@s@t@@n}}{\end{pr@p@s@t@@n}\Inthmfalse}
\newenvironment{cor}{\Inthmtrue\begin{c@r@ll@ry}}{\end{c@r@ll@ry}\Inthmfalse}

\newenvironment{rem}{\Inthmtrue\begin{r@m@rk}}{\end{r@m@rk}\Inthmfalse}

\newenvironment{pf}{\Inthmtrue\begin{proof}}{\end{proof}\Inthmfalse}

\newenvironment{THM}{\Inthmtrue\begin{TH@@r@m}}{\end{TH@@r@m}\Inthmfalse}
\newenvironment{claim}[1][Claim]{
   \trivlist
  \thm@preskip\topsep
  \thm@postskip\thm@preskip
  \@topsep \thm@preskip               
  \@topsepadd \thm@postskip           
   \item[\hskip\labelsep\itshape#1\@addpunct{.}]\ignorespaces
}{\endtrivlist}
\newenvironment{pfofclaim}[1][Proof of Claim]{\par
  \normalfont
  \topsep6\p@\@plus6\p@ \trivlist
  \item[\hskip\labelsep\itshape
    #1\@addpunct{.}]\ignorespaces
}{%
  \endtrivlist
}
\newcounter{subthm}[th@@r@m]
\newif\ifInthm
\def\MultiNumbering{
 \def\Inthmlabel{\theth@@r@m.\thesubthm}
 \def\c@equation{\ifInthm\c@subthm\else\c@th@@r@m\fi}
 \def\theequation{\ifInthm\Inthmlabel\else\theth@@r@m\fi}
}
\def\TagsOnRight{\tagsleft@false}
\def\joken{\refstepcounter{subthm}%
  \edef\@currentlabel{\theth@@r@m.\thesubthm}%
     \ignorespaces{\@currentlabel}}
\evensidemargin=\oddsidemargin
\def\SkipAmount{
\lineskiplimit=3pt
\lineskip=3pt plus 1pt
}
\newcommand{\Twoside}{\@twosidetrue  \@mparswitchtrue}
\let\bib\bibitem
\def\lb#1{\label{#1}}
\newif\ifInFtn
\newif\ifFtnMark
  
  \def\NoFootnoteMark{\FtnMarkfalse}
\renewcommand{\@makefntext}[1]{%
  \parindent 1em%
  \ifFtnMark
  	\noindent\hb@xt@1.8em{\hss\@makefnmark}#1
   \else
   	\noindent\hb@xt@1.8em{}#1
    \fi}
\def\lqq{\lq\lq} 
\def\rqq{\rq\rq}

\def\indno#1{\par\noindent\hskip3em\llap{\upshape(#1)}\hskip1em\ignorespaces}
\def\no#1{{\upshape(#1)}\hskip1em\ignorespaces}

\def\RHeads#1{\markright{\hfill\textsl{#1}\hspace{2em}}}
\newif\ifInEnum
\def\parenumerate{%
	\ifnum \@enumdepth >\thr@@\@toodeep\else
    \advance\@enumdepth\@ne
      \edef\@enumctr{enum\romannumeral\the\@enumdepth}%
      \expandafter
      \list
        \csname label\@enumctr\endcsname
			  {\setlength{\leftmargin}{0pt}\setlength{\labelwidth}{1.2em}
         \setlength{\labelsep}{.5em}
		     \addtolength{\itemindent}{\labelwidth}
         \addtolength{\itemindent}{\labelsep}
         \addtolength{\itemindent}{\parindent}
         \setlength{\topsep}{\itemsep}
         \setlength{\parsep}{0pt}
         \usecounter\@enumctr\def\makelabel##1{\hss\llap{\upshape##1}}}%
  \fi}

\renewcommand{\enumerate}{%
  \ifnum \@enumdepth >\thr@@\@toodeep\else
    \advance\@enumdepth\@ne
    \edef\@enumctr{enum\romannumeral\the\@enumdepth}%
      \expandafter
      \list
        \csname label\@enumctr\endcsname
			  {\setlength{\leftmargin}{0pt}\setlength{\labelwidth}{1.2em}
         \setlength{\labelsep}{.5em}
		     \addtolength{\leftmargin}{\labelwidth}
         \addtolength{\leftmargin}{\labelsep}
         \addtolength{\leftmargin}{\parindent}
         \setlength{\topsep}{\itemsep}
         \setlength{\parsep}{0pt}
           \usecounter\@enumctr\def\makelabel##1{\hss\llap{##1}}}%
  \fi}

\newenvironment{enumr}{\InEnumtrue

	\begin{enumerate}}{\end{enumerate}\InEnumfalse}
\newenvironment{parenumr}{\InEnumtrue

	\begin{parenumerate}}{\end{parenumerate}\InEnumfalse}
\newenvironment{enuma}{\InEnumtrue

	\begin{enumerate}}{\end{enumerate}\InEnumfalse}

\def\Thm#1{Theorem #1}

\def\Cor#1{Corollary #1}

\def\Prop#1{Proposition #1}
\def\Props#1{Propositions #1}

\def\Lem#1{Lemma #1}
\def\Lems#1{Lemmas #1}
\def\Df#1{Definition #1}

\def\Ex#1{Example #1}
\def\Rem#1{Remark #1}

\def\case#1{\par\noindent{\it Case} {#1}.\quad\ignorespaces }

\def\ch#1{\operatorname{char}(#1)}
\def\ee#1#2#3{#1 \leq #2 \leq #3}
\def\ne#1#2#3{#1 < #2 \leq #3}
\def\en#1#2#3{#1 \leq #2 < #3}
\def\nn#1#2#3{#1 < #2 < #3}
\def\ds#1{\mathop\bigoplus_{#1}}
\def\dss#1#2{\mathop\bigoplus_{#1}^{#2}}
\def\sm#1{\mathop\sum_{#1}}      

\def\smm#1#2{\mathop\sum_{#1}^{#2}}
 
\def\prdd#1#2{\mathop\prod_{#1}^{#2}}

\def\ep#1{\mathop{\bigwedge}^{\kern-0.3pt#1}\;}

\def\Hm#1#2{H^{#1}_\mathfrak{m}(#2)}

\def\Homcal#1#2#3{{{\mathcal H}\kern-1.2pt\textit{om}}\,_{#1}(#2,#3)}

\def\Ext#1#2#3#4{\operatorname{Ext}^{#1}_{#2}(#3,#4)}  
\def\Extcal#1#2#3#4{{{\mathcal E}\kern-1pt\textit{xt}}\,^{#1}_{#2}(#3,#4)}

\def\Imasup#1#2{\operatorname{Im}^{#1}(#2)}
\def\Coksup#1#2{\operatorname{Coker}^{#1}(#2)}

\def\Proj#1{\operatorname{Proj}(#1)}
\def\Spec#1{\operatorname{Spec}(#1)}

\def\GL{\textit{GL}}
\def\gradplain#1#2{[#1]_{#2}}
\def\gradn#1#2{\lt[#1\rt]_{#2}}
\def\grad#1#2{\bigl[#1\bigr]_{#2}}

\def\directlim#1{
\ooalign
{\hfil\hbox{$\varinjlim$}\hfil\crcr
\leavevmode\hfil\raise-10pt\hbox{\scriptsize$#1$}\hfil}
}
\def\len#1#2{l_{#1}(#2)}
\def\dep#1#2{\operatorname{depth}_{#1}(#2)}
\def\depm#1{\operatorname{depth}_\mfrak(#1)}

\def\rk#1#2{\operatorname{rank}_{#1}(#2)}
\def\rank{\operatorname{rank}}
\def\dim#1{\operatorname{dim}(#1)}
\def\dimm#1#2{\operatorname{dim}_{#1}(#2)}

\def\dimk#1{\operatorname{dim}_k(#1)}

\def\row#1#2#3{{({#1}_{#2},{\ldots}\,,{#1}_{#3})}} 

\def\trow#1#2#3{{{}^t\kern-2pt({#1}_{#2},{\ldots}\,,{#1}_{#3})}}
\def\tpose#1{{{}^t\kern-1pt#1}}
\def\trowiv#1#2#3#4{{{}^t\kern-1.2pt({#1},{#2},{#3},{#4})}}

\def\seq#1#2#3{{{#1}_{#2},{\ldots}\,,{#1}_{#3}}} 

\def\ddd{,{\ldots}\,,}

\def\op{\oplus}

\def\ot{\otimes}

\def\tm{\times}    

\def\opdot{{
\,\ooalign
{\hfil\raise 0.35ex\hbox{$\centerdot$}\hfil\crcr
\leavevmode\hbox{$\oplus$}}\,
}}

\def\indexonsmile#1{
\ooalign
{\hspace{-3pt}\raise -0.3ex\hbox{\scriptsize$\smallsmile$}\crcr
\leavevmode\hfil\raise0.45ex\hbox{\scriptsize$#1$}\hfil}
}
\def\lra{\longrightarrow}
\def\xr@ght@rr@w[#1]#2{\xrightarrow[#1]{#2}}
\def\xxr@ght@rr@w#1{\xrightarrow{#1}}
\def\xra{\@ifnextchar[{\xr@ght@rr@w}{\xxr@ght@rr@w}}
\def\xl@ft@rr@w[#1]#2{\xleftarrow[#1]{#2}}
\def\xxl@ft@rr@w#1{\xleftarrow{#1}}
\def\xla{\@ifnextchar[{\xl@ft@rr@w}{\xxl@ft@rr@w}}

\def\la{\langle}
\def\ra{\rangle}

\def\cc{\circ}
\def\lt{\left}
           
\def\rt{\right}

\def\sset{\subset}

\def\bslash{\backslash}

\let\ggreater\gg
\let\lless\ll

\def\Pbf{\text{\bf P}}

\def\Zbf{\text{\bf Z}}

\def\Sfrak{\mathfrak{S}}

\def\afrak{\mathfrak{a}}

\def\mfrak{\mathfrak{m}}

\def\pfrak{\mathfrak{p}}

\def\Ical{\mathcal{I}}

\def\Kcal{\mathcal{K}}

\def\Ocal{\mathcal{O}}

\def\Wcal{\mathcal{W}}

\def\Atild{{\tilde A}}

\def\Ctild{{\tilde C}}
\def\Dtild{{\tilde D}}

\def\Ftild{{\tilde F}}
\def\Gtild{{\tilde G}}
\def\Htild{{\tilde H}}

\def\Utild{{\tilde U}}
\def\Vtild{{\tilde V}}

\def\Ztild{{\tilde Z}}

\def\btild{{\tilde b}}
\def\ctild{{\tilde c}}

\def\etild{{\tilde e}}

\def\vtild{{\tilde v}}

\def\lgtild{{\tilde \lg}}

\def\tAtild{{^t\kern-0.2em \Atild}}
\def\tA{{^t\kern-0.2em A}}
\def\tB{{^t\kern-0.2em B}}
\def\tF{{^t\kern-0.2em F}}
\def\tFtild{{^t\kern-0.2em \Ftild}}
\def\tPhg{{^t\kern-0.1em\Phg}}
\def\Jlsup#1{{{}^{#1}\kern-0.2em J}}
\def\Kcallsup#1{{{}^{#1}\kern-0.1em\Kcal}}

\def\Hbar{{\bar H}}

\def\Zbar{{\bar Z}}

\def\dbar{{\bar d}}

\def\hbar{{\bar h}}

\def\ubar{{\bar u}}
\def\vbar{{\bar v}}

\def\xbar{{\bar x}}

\def\lgbar{{\bar \lg}}

\def\sgbar{{\bar \sg}}
\def\tgbar{{\bar \tg}}

\def\fhat{{\hat f}}
\def\ghat{{\hat g}}

\def\Dcheck{{\Check D}}

\def\Icheck{{\Check I}}

\def\Ucheck{{\Check U}}

\def\ccheck{{\Check c}}

\def\echeck{{\Check e}}

\def\ag{\alpha}
\def\bg{\beta}
\def\dg{\delta}
\def\gg{\gamma}
\def\kg{\kappa}

\def\lg{\lambda}
\def\mg{\mu}
\def\ng{\nu}
\def\sg{\sigma}
\def\tg{\tau}
\def\eg{\varepsilon}
\def\ig{\iota}
\def\pg{\pi}

\def\thg{\theta}
\def\Thg{\varTheta}
\def\zg{\zeta}
\def\Gg{\varGamma}
\def\Dg{\varDelta}
\def\Phg{\varPhi}

\def\og{\omega}

\def\Xg{\varXi}
\def\xg{\xi}
\def\rg{\rho}
\def\phg{\varphi}
\def\psg{\psi}
\def\hg{\eta}
\def\chg{\chi}

\def\mod{\operatorname{mod}}

\def\Imr{\text{\rm Im}}

\def\MATr{\text{\rm MAT}}

\def\fort{\text{for}}

\def\andt{\text{and}}

\def\andr{\text{\rm and}}

\def\witht{\text{with}}

\def\ift{\text{if}}

\def\ort{\text{or}}

\def\fallt{\text{for all}}

\def\inr{\text{\rm in}}

\def\rankr{\text{\rm rank}}

\def\prr{\text{\rm pr}}

\def\pt{\text{\it pt}}

\def\Bs{\textit{Bseq}}
\def\irrd{{\textit{irrd}}}
\def\CM{Cohen-Macaulay}

\def\Gor{Gorenstein}
\def\qBbm{quasi-Buchsbaum}

\def\Bbm{Buchsbaum}
\def\Wei{Weierstrass}

\def\fg{finitely generated}

\def\Bou{Bourbaki}
\def\Grob{Gr\"obner}

\newenvironment{smallmatrixl}{\null\,\vcenter\bgroup
  \Let@\restore@math@cr\default@tag
  \baselineskip8\ex@ \lineskip1.5\ex@ \lineskiplimit\lineskip
  \ialign\bgroup$\m@th\scriptstyle##$\hfill&&\thickspace
  $\m@th\scriptstyle##$\hfill\crcr
}{%
  \crcr\egroup\egroup\,%
}
\newenvironment{matrixl}{%
  \hskip -\arraycolsep\array{*\c@MaxMatrixCols l}%
}{%
  \endarray \hskip -\arraycolsep
}
\newenvironment{smallbmatrix}{\lt[\begin{smallmatrix}}
{\end{smallmatrix}\rt]}
\newcommand{\cdotsfor}[1]{%
  \ifx[#1\@xp\scdots@for\else\cdots@for\@ne{#1}\fi}
\def\scdots@for#1]{\cdots@for{#1}}
\def\cdots@for#1#2{\multicolumn{#2}c%
  {\m@th\dotsspace@1.5mu\mkern-#1\dotsspace@
   \xleaders\hbox{$\m@th\mkern#1\dotsspace@\cdot\mkern#1\dotsspace@$}%
           \hfill
   \mkern-#1\dotsspace@}%
   }

\def\clml#1#2{
\lt(\begin{matrixl}
{#1} \\ {#2} 
\end{matrixl}\rt)
}                                   
\def\clmlt#1#2{
\text{\footnotesize$\lt( 
\begin{matrixl}  {#1} \vspace{0.2ex}\\ 
{#2} \end{matrixl} \rt)$}
}
\def\clmlsc#1#2{
\lt(\begin{smallmatrixl}
{#1}\\ {#2} 
\end{smallmatrixl}\rt)
}
\def\clm#1#2{
\begin{pmatrix}  {#1} \\ {#2} \end{pmatrix} 
}                                   

\def\clmsc#1#2{
\lt( \smallmatrix  {#1} \\ 
{#2}
\endsmallmatrix \rt)
}
\def\bclmsc#1#2{
\lt[ \smallmatrix  {#1} \\ 
{#2}
\endsmallmatrix \rt]
}
\def\bclm#1#2{
\bmatrix  {#1} \\ {#2} \endbmatrix 
}                                   

\def\kx#1{k[x_{#1}\ddd x_4]}

\def\ktwo{{k[x_3,x_4]}}
\def\kthree{{k[x_2,x_3,x_4]}}

\def\kthreez{{k[z_2,z_3,z_4]}}

\def\Kthreez{{K[z_2,z_3,z_4]}}
\def\seqtwo{{x_3,x_4}}
\def\seqthree{{x_2,x_3,x_4}}
\def\seqfour{{x_1,x_2,x_3,x_4}}

\def\seqfourz{{z_1,z_2,z_3,z_4}}
\def\supbr#1#2{{{#1}^{[#2]}}}
\def\supang#1#2{{{#1}^{\langle#2\rangle}}}

\def\suplpar#1#2{{{}^{(#2)}\kern-1pt{#1}}}
\def\sublpar#1#2{{{}_{(#2)}\kern-1pt{#1}}}
\def\BR#1{B_R(#1)}

\newcommand{\kernsub}[1]{\kern-0.2em{{}_{#1}}}
\newcommand{\kernnsub}[1]{\kern-0.3em{{}_{#1}}}

\def\Pthree{{\Pbf^3}}
\def\overcc#1{{\ooalign{\hfil\raise1.7ex\hbox{\scriptsize$\cc$}\hfil\crcr $#1$}}}
\def\Uovercc{{\overcc{U}}}
\def\Uoverccnatural#1{{\overcc{U}{}^\natural_{#1}}}

\def\Vovercc{{\overcc{V}}}
\def\Wovercc{{\overcc{W}}}

\def\overtildcc#1{{\ooalign{\hfil\raise2.2ex\hbox{\scriptsize$\cc$}\hfil\crcr $#1$}}}

\edef\today{\ifcase\month\or
  January\or February\or March\or April\or May\or June\or
  July\or August\or September\or October\or November\or December\fi
  \space\number\day, \number\year}
\def\@maketitle{%
  \newpage
  \null
  \vskip 2em%
  \begin{center}%
    {\LARGE \@title \par}%
    \vskip 1.5em%
    {\large
      \lineskip .5em%
      \begin{tabular}[t]{c}%
        \@author
      \end{tabular}\par}%
    \vskip .5em%
    {\footnotesize \@date}%
  \end{center}%
  \par
  \vskip 1.5em}

\def\Watashidesu{Mutsumi A{\scshape masaki}\vspace{1ex}\\\tokoro\\\email}

\def\@date{February 8, 2005}
\makeatother
\begin{document}

\SkipAmount
\allowdisplaybreaks
\NoFootnoteMark

\RHeads{Basic Sequence of an Integral Curve}
\title{Inequalities satisfied by the basic sequence\\
of an integral curve in $\Pthree$}
\author{\Watashidesu}
\maketitle
\footnotetext{
2000 Mathematics Subject Classification. 
Primary 14H50; Secondary 13P10, 13D02.
}

\begin{abstract}
We show several new inequalities found recently 
that the basic sequence of 
the saturated homogeneous ideal $I$ of an integral curve in $\Pthree$ 
must satisfy. 
Then we compare our results with Cook's assertions on the 
generic initial ideal of $I$, carrying out numerical computations
with the use of a computer. The outcome is that we have obtained 
new restrictions on the generic initial ideal of $I$.
When $\ch{k}=0$, the basic sequence of a homogeneous ideal $I$ in a 
polynomial ring over an infinite field $k$ is a sequence of the degrees of 
the minimal generators of the generic initial ideal of $I$ arranged in a 
suitable order. 
\end{abstract}

\section*{Introduction}

Let $X$ be a curve in $\Pthree:=\Pbf^3_k$, namely, a locally \CM\ 
equidimensional closed
subscheme of $\Pthree$ of dimension one, and let $I$ denote the
saturated homogeneous ideal of $X$ in a polynomial ring 
$R:=k[\seqfour]$, where $k$ is an infinite field and 
the linear forms $\seqfour$ are chosen
sufficiently generally. Then there is what we call a \Wei\ basis
$\{\ e^i_l\ |\ \ee{1}{i}{3},\ \ee{1}{l}{m_i}\ \}$  of $I$ that is characterized by some
properties connected with a representation of the multiplication of
$e^i_l$ by the linear forms $\seqfour$ (see \cite{A5} and \cite{A7}). 
More precisely,
\begin{equation*}
\lt\{\begin{aligned}
{}&\supbr{I}{1}=I,\quad \supbr{I}{3}=\supang{I}{3},\ 
\supbr{I}{i}=\supang{I}{i}\op\supbr{I}{i+1}\ (i=1,\ 2)\\
{}&\text{as $\kx{i+1}$-modules},\ \andt\ 
x_i\supbr{I}{i+1}\sset(\seq{x}{i+1}{4})\supang{I'}{i}\op\supbr{I}{i+1},
\end{aligned}\rt.
\end{equation*}
where $\supbr{I}{i}$ is a \fg\ graded $\kx{i}$-submodule of $I$ and 
$\supang{I}{i}$ (resp. $\supang{I'}{i}$) is a \fg\ graded free $\kx{i}$-submodule 
(resp. $\kx{i+1}$-submodule) of $I$ with free basis 
$\{\ e^i_l\ |\ \ee{1}{l}{m_i}\ \}$
for each $i\ (\ee{1}{i}{3})$. We have $m_1=1$ and $\deg(e^1_1)=a=m_2$. Put $b:=m_3$.
Let $n_l:=\deg(e^2_l)$ $(\ee{1}{l}{a})$ and
$n_{a+l}:=\deg(e^3_l)$ $(\ee{1}{l}{b})$.
 We may assume that
the sequence $n_1\ddd n_a$ (resp. $n_{a+1}\ddd n_{a+b}$) is
nondecreasing after changing the order of $e^2_1\ddd e^2_a$
(resp. $e^3_1\ddd e^3_b$) if necessary.
The sequence $\BR{I}:=(a;\seq{n}{1}{a};\seq{n}{a+1}{a+b})$ is determined
uniquely by $X\sset\Pthree$ and is called the basic sequence of $I$ or $X$
(see \cite{A3} and \cite{A8} -- \cite{A7}). 
Note that $X$ is projectively \CM\ if and only if
$b=0$ (i.e. the subsequence $(\seq{n}{a+1}{a+b})$ is vacuous).
In general, $a\leq n_1\leq\cdots\leq n_a$ and 
$n_1\leq n_{a+1}\leq\cdots\leq n_{a+b}$.
There are a number of applications of basic sequence to the study of
curves in $\Pthree$ (see \cite{A2} -- \cite{A9}).

\par
In this paper, expanding the ideas in our earlier work \cite[Section 1]{A1},
we show several new inequalities that the components of
$\BR{I}$ must satisfy when $X$ is integral and $\BR{I}$
takes special forms (see Sections \ref{sec4} and \ref{sec5}).
Our main results are \Thm{\ref{proc118}}, \Cor{\ref{proc119}},
\Props{\ref{proc12}, \ref{proc34}}, and \Thm{\ref{proc39}}. 
They can be summarized together with our related earlier results 
in the following manner.

\begin{THM}
Assume that $X$ is contained in an irreducible surface of degree $a\geq 2$
and that $b\geq1$.
Let $n'_1\ddd n'_\og$ $(\og\geq1)$ be a strictly increasing sequence of integers
such that $\{n'_1\ddd n'_\og\}=\{\ n_i\ |\ \ee{1}{i}{a}\ \}$ and 
let $t_l:=\#\{\ i\ |\ n_i=n'_l,\ \ee{1}{i}{a}\ \}$.
Then, the following inequalities hold, where an additional assumption $\ch{k}=0$ is
needed in our proofs of \eqref{cond110} -- \eqref{cond112}.
\begin{parenumr}
\item\lb{cond115}
We have $n_i\leq n_{i+1}\leq n_i+1$ for all $\ee{1}{i}{a-1}$.
\item\lb{cond110}
Suppose that $\og>1$.
Let $m$ be an integer with $\en{1}{m}{\og}$, $t:=\smm{l=1}{m}t_l$, and 
$a'$ be an integer with $\ee{t+1}{a'}{a}$.
If $n_i=n_t+i-t$ for all $\ee{t}{i}{a'}$  and $n_{a+1}<n_{a'}$, then
$t(t-1)/2 > p$, where $p:=\max\{\ i\ |\ n_{a+i}<n_{a'},\ \ee{1}{i}{b}\ \}$.
\item\lb{cond111}
Let $a'$ be an integer with $\ee{2}{a'}{a}$. 
If $n_i=n_1+i-1$ for all $\ee{1}{i}{a'}$, then $n_{a+1}\geq n_{a'}$.
\item\lb{cond112}
Let $a'$ be an integer with $\ee{3}{a'}{a}$. 
If $n_i=n_1+i-2$ for all $\ee{2}{i}{a'}$, then $n_{a+1}\geq n_{a'}$.
\item\lb{cond113}
If $n_i=n_1+i-1$ for all $\ee{1}{i}{a}$,
then $b>2$ and $n_{a+j}\neq n_{a+1}+j-1$ for some $\ee{1}{j}{b}$. 
If further there is an integer  $b_0$ with $\nn{0}{b_0}{b}$ such that
$n_{a+b_0}+1<n_{a+b_0+1}$, then $b_0>2$ and $n_{a+j}\neq n_{a+1}+j-1$
for some $\ee{1}{j}{b_0}$.
\item\lb{cond114}
Suppose that $b=1$ and that $X$ does not contain any line as an
irreducible component.  
Delete the terms $n'_{2}\ddd n'_{\og}$ from the sequence $(a,\seq{n}{1}{a},n_{a+1})$
and then rearrange the remaining terms so that they make a nondecreasing sequence.
Let $(n''_1\ddd n''_{a+3-\og})$ denote the resulting sequence.
Then 
\begin{equation*}
n_{a+1}\leq
a-2+\smm{j=1}{a}n_j -\smm{\gg=2}{\og}n'_\gg
-\smm{i=1}{a-\og}n''_i.
\end{equation*}
\end{parenumr}
\end{THM}

\par
Let us explain the motivation of our study in this paper.
For a sequence $\Bs=(a;n_1\ddd n_a;\seq{n}{a+1}{a+b})$ satisfying
$0<a\leq n_1\leq\cdots\leq n_a$ and 
$n_1\leq n_{a+1}\leq\cdots\leq n_{a+b}$ with $b\geq0$, put
\begin{align*}
&D(\Bs):=\smm{l=1}{a}n_l-\frac{1}{2}a(a-1)-b,\\
&G(\Bs):=1+\smm{l=1}{a}\frac{1}{2}n_l(n_l-3)-\smm{l=1}{b}n_{a+l}+b
-\frac{1}{6}a(a-1)(a-5).
\end{align*}
Denoting the degree and the arithmetic genus of $X$ by
$d(X)$ and $g(X)$ respectively,
we have $d(X)=D(\BR{I})$ and $g(X)=G(\BR{I})$ (see \cite[\Rem{1.9}]{A3}).
With this in mind, we are interested in the following problem.
{\it For each pair of nonnegative integers $d,\ g$, give a good characterization
of a sequence
$\Bs=(a;n_1\ddd n_a;\seq{n}{a+1}{a+b})$ for which there is
an integral curve $X$ in $\Pthree$ such that $\Bs=\BR{I}$,
$d=D(\BR{I})$, and $g=G(\BR{I})$.}
The case $b=0$ was settled a long time ago in \cite{GP1}.
But we do not have any complete general answer to this problem as yet.

\par
When $\ch{k}=0$, the reduced \Grob\ basis of $I$ with respect to 
the reverse lexicographic order becomes a \Wei\ basis of $I$ such that 
$\inr^x(e^1_1)=x_1^a$, $\inr^x(e^2_l)=x_1^{a-l}x_2^{\bg_{a-l}}$ $(\ee{1}{l}{a})$,
$\inr^x(e^3_{l})=x_1^{t_l}x_2^{u_l}x_3^{\bg_{t_lu_l}}$ $(\ee{1}{l}{b})$,
where $1\leq\bg_{a-1}<\bg_{a-2}<\cdots<\bg_0$ and
$t_l<a,\ u_l<\bg_{t_l},\ \bg_{t_lu_l}>0$ (see \cite{A3}, \cite{A8} and \cite{A7}).
Note that $n_l={a-l}+\bg_{a-l}$ $(\ee{1}{l}{a})$,
$n_{a+l}={t_l}+{u_l}+\bg_{t_lu_l}$ $(\ee{1}{l}{b})$. 
One may therefore expect to see some relation between our results and
the main assertions in Cook's paper \cite{C}, though it is pointed out 
that there is an error in the proof of the main theorem of \cite{C} 
(see \cite[Section 4]{DS}). 
There is however no more than
one result apparently common to both, namely, the inequalities
$n_i\leq n_{i+1}\leq n_i+1$ for all $\ee{1}{i}{a-1}$ (see \cite[\Cor{1.2}]{A1}
or \Lem{\ref{proc26}}).
In order to compare our results with Cook's main assertions, 
we have carried out numerical computations with the 
use of a computer which show us sequences $(a;\seq{n}{1}{a};\seq{n}{a+1}{a+b})$
satisfying the conditions 
stated in this paper and the sequences coming from the
generic initial ideals satisfying the assertions in Cook's paper.
The outcome is that neither implies the other in general.
By the same computations we can observe that, in a certain narrow range of
the pairs of degree and genus, our
results in this paper seem effective for the classification of
integral curves in $\Pthree$ (see Section \ref{sec7}).

\par
As in our previous papers, our arguments depend heavily on the knowledge of 
the relation matrices $\lg_2$ and $\lg_3$ appearing in the standard free resolution 
which starts with $e^1_1,e^2_1,\ldots,e^2_a,e^3_1\ddd e^3_b$ (see \cite{A2}, \cite{A8}). 

\par
In Section \ref{sec1}, basic knowledge of \Wei\ bases and standard free resolutios are summarized.
In Section \ref{sec2},  new technical results on $\lg_2$ which play an inportant role in the proof of
our main theorem  are given. We have also included some 
of the points, seen in \cite[Section 1]{A1} and necessary for our purposes,
in a more adaptable form for the convenience of the readers.
In Sections \ref{sec3} and \ref{sec8}  the properties of the matrix $\lg_3$ are studied
for the purpose of analyzing the structure of $\Ext{3}{R}{R/I}{R}$ as $\ktwo$-module.
In Sections \ref{sec4} and \ref{sec5}, our main results are given with full use of the results of
the previous sections. Finally in Section \ref{sec7}, 
comparison of our results and Cook's main assertions are carried out briefly with
numerical cmputations.

\section*{Notation}

\begin{enumr}
\item
Given a set of polynomials, say $Z$, we denote by $\MATr(Z)$ the
set of matrices with entries in $Z$.

\item
The symbol $\op$ will be used in the following two senses\,:
\begin{enuma}
\item
$E' \op
E'' = \{\ (e',e'') \ |\ e' \in E', \ e'' \in E'' \ \}$, 
\item
$E' \sset E, \ E''
\sset E, \ E' \cap E'' = 0, \ E' \op E'' = \{\ e' + e'' \ |\ e' \in  E', \ e'' 
\in E'' \ \} \sset E$. 
\end{enuma}
Usually the context will make it clear 
which it means.
But when direct sums in both meanings appear in a single formula
simultaneously, we will use another symbol $\opdot$ instead of
$\op$ to express the direct sum in the first sense.

\item
Given a graded module $E = \ds{t} [E]_t$, integers $p,q \ (q \geq 0)$ 
and  a sequence of integers $\dbar = (\seq{d}{1}{n})$, we set 
$-\dbar:=(-d_1\ddd -d_n)$, $\dbar + p :=
(d_1 + p \,\ddd d_n + p)$, and
$E(\dbar) := \dss{l=1}{n}E(d_l)$ (in the  sense (a)), 
where $[E(d_l)]_t := [E]_{{d_l}+t}$.  

\item
Let $F=(f_{ij})$ be a matrix whose entries are homogeneous 
polynomials and $\Delta=(d_{ij})$ be a matrix of integer coefficients 
of the same size as $F$. We write $\Dg(F)=\Delta$ to mean that 
$\deg(f_{ij})=d_{ij}$ for all $i,j$ such that $f_{ij}\neq0$.

\item
For a matrix $\Delta=(d_{ij})$ with entries in $\Zbf$ and an
integer $n$, we denote by $\Delta+n$ the matrix $(d_{ij}+n)$.

\item
For an integer $n$, we denote by $1_n$ the $n\times n$ unit
matrix.

\item
Let $F$ be an $m\times n$ matrix in a commutative ring 
$R$ and $P$ be a subring of $R$.
The image of the linear map from $P^n$ to $R^m$ over $P$
defined by multiplication by $F$ will be denoted by $\Imasup{P}{F}$.

\item
Let $F$ be a matrix, and $\seq{s}{1}{p}$ and $\seq{t}{1}{q}$ be strictly
increasing sequences
of integers. We denote by $F\clmlsc{\seq{s}{1}{p}}{\seq{t}{1}{q}}$
the matrix obtained from $F$ by deleting its $s_i$-th 
rows $(\ee{1}{i}{p})$ and $t_j$-th columns $(\ee{1}{j}{q})$.

\end{enumr}

\section{\mathversion{bold}\Wei{} bases and standard free resolutions of homogeneous ideals 
defining curves in $\Pthree$}
\lb{sec1}

We summarize here some properties of \Wei\ bases and free resolutions 
of homogeneous ideals
defining curves in $\Pthree$ which are necessary to prove our main theorems, 
along with elementary results and arguments from 
\cite{A2} -- \cite{A1} and \cite{A8} -- \cite{A7}
for the  convenience of the readers.

\par
Let $y_1,y_2,y_3,y_4$ be indeterminates over an infinite field $k$, 
$R$ the polynomial ring $k[y_1,y_2,y_3,y_4]$, $\mfrak$ the maximal
ideal $(y_1,y_2,y_3,y_4)$,
$\gg_{ij}\ (\ee{1}{i}{4},\ \ee{1}{j}{4})$ elements of
$k$ such that the matrix $\Gg:=(\gg_{ij})$ is invertible, and $x_1,x_2,x_3,x_4$ 
elements of $R$ satisfying
$y_i=\smm{j=1}{4}\gg_{ji}x_j$ $(\ee{1}{i}{4})$.
\par
Let $I$ be a homogeneous ideal in $R$ of height 2 such 
that $\depm{R/I}\geq1$. 
Then, for a sufficiently general choice of  $\Gg$,
there are \fg\ graded $\kx{i}$-submodules $\supbr{I}{i}\sset I$
$(\ee{1}{i}{3})$ and \fg\ graded free $\kx{i}$-submodules 
$\supang{I}{i}\sset I$ $(\ee{1}{i}{3})$
such that 
\begin{equation}
\lt\{\begin{aligned}
{}&\supbr{I}{1}=I,\quad \supbr{I}{3}=\supang{I}{3},\ 
\supbr{I}{i}=\supang{I}{i}\op\supbr{I}{i+1}\ (i=1,\ 2)\\
{}&\text{as $\kx{i+1}$-modules},\ \andt\ 
x_i\supbr{I}{i+1}\sset(\seq{x}{i+1}{4})\supang{I}{i}\op\supbr{I}{i+1}
\end{aligned}\rt.
\lb{eq200}
\end{equation}
by the results of \cite[Section 2]{A5}.
For each $ \ee{1}{i}{3}$, let 
$\{\ e^i_l\ |\ \ee{1}{l}{m_i}\ \}$ be a free basis of
$\supang{I}{i}$ and $a:=\min\{\ l\ |\ \gradplain{I}{l}\neq0\ \}$. 
Since $\rk{R}{I}=1$ and $\dim{R/I}=2$, 
we have $m_1=1$ and $\deg(e^1_1)=a=m_2$.
Put $b:=m_3$, 
$n_l:=\deg(e^2_l)\ (\ee{1}{l}{a})$, 
$n_{a+l}:=\deg(e^3_l)\ (\ee{1}{l}{b})$. We may assume that
the sequence $n_1\ddd n_a$ (resp. $n_{a+1}\ddd n_{a+b}$) is
nondecreasing after changing the order of $e^2_1\ddd e^2_a$
(resp. $e^3_1\ddd e^3_b$) if necessary.
When the matrix $\Gg$ is chosen sufficiently generally, 
the sequence $\BR{I}:=(a;n_1\ddd n_a;n_{a+1}\ddd n_{a+b})$ is determinned
uniquely by $I$ not depending on $\Gg$, and is called the basic sequence of
$I$ (see \cite[\Df{2.13}]{A5}, \cite[\Df(1.5)]{A8}, \cite[\Df{1.4}]{A3}).
It is also called the basic sequence of $\Proj{R/I}$.
Note that $a\leq n_1$. Recall further 
that $b=0$ if and only if $\depm{R/I}=2$ and that
$a\leq n_1\leq n_{a+1}$ when $b>0$ 
(see \cite[\Thm{(1.1)} and \Prop{(1.6)}]{A8}).
By \eqref{eq200}, we have
\begin{equation}
I=R e^1_1\op\lt(\dss{l=1}{a}\kthree e^2_l\rt)
\op\lt(\dss{l=1}{b}\ktwo e^3_l\rt)
\lb{eq40}
\end{equation}
as $\ktwo$-module.
If the system of generators $\{e^1_1,e^2_1,\ldots,e^2_a,e^3_1\ddd e^3_b\}$
mentioned above satisfies the additional conditions
\begin{equation}
\lt\{\begin{aligned}
{}&x_1e^2_{l'}\in 
(\seqthree)\kthree e^1_1\op
\lt(\dss{l=1}{a}\kthree e^2_l\rt)\\
{}&\hskip10em\op\lt(\dss{l=1}{b}\ktwo e^3_l\rt)
\quad\fort\quad\ee{1}{l'}{a},\\
{}&x_1e^3_{l'}\in
(\seqtwo)\kthreez e^1_1\op
(\seqtwo)\lt(\dss{l=1}{a}\kthree e^2_l\rt)\\
{}&\hskip10em\op\lt(\dss{l=1}{b}\ktwo e^3_l\rt)
\quad\fort\quad\ee{1}{l'}{b},\\
{}&x_2e^3_{l'}\in
(\seqtwo)\lt(\dss{l=1}{a}\ktwo e^2_l\rt)
\op\lt(\dss{l=1}{b}\ktwo e^3_l\rt)
\quad\fort\quad\ee{1}{l'}{b},
\end{aligned}\rt.
\lb{eq201}
\end{equation}
then we call it a \Wei\ basis of $I$ with respect to $\seqfour$.
One can always find a \Wei\ basis of $I$ (see \cite[\Thm{2.12}]{A5},
\cite[\Thm{2.5}]{A7}, \cite[\Thm{(1.1)}]{A8}).

\begin{rem}\lb{proc40}
Let $\Gg$ be sufficiently general and let
$\inr^x(I)$ be the ideal in $R$ generated by the monomials
\begin{equation*}
\{\ \inr^x(f)\ |\ \text{$f$ is a homogeneous polynomial of $I$}\ \}
\end{equation*}
in $\seqfour$, 
where $\inr^x(f)$ denotes the initial term of $f$ with respect to the 
reverse lexicographic order associated with the variables $\seqfour$.
It is often called the generic initial ideal of $I$.
In the case where $\ch{k}=0$, we can obtain a 
\Wei{} basis 
$\{e^1_1,e^2_1,\ldots,e^2_a,e^3_1\ddd e^3_b\}$
of $I$ with respect to $\seqfour$ practically 
by taking the reduced and minimal \Grob{} basis of $I$ 
with respect to the order mentioned above in such a way that 
\begin{gather*}
\inr^x(e^1_1)=x_1^a,\quad 
\inr^x(e^2_l)=x_1^{a-l}x_2^{\bg_{a-l}}\ (\ee{1}{l}{a}),\\
\inr^x(e^3_l)=x_1^{t_l}x_2^{u_l}x_3^{\bg_{t_lu_l}}\ (\ee{1}{l}{b}),
\end{gather*}
where $1\leq\bg_{a-1}<\bg_{a-2}<\cdots<\bg_0$,
$t_l<a,\ u_l<\bg_{t_l},\ \bg_{t_lu_l}>0$ for $\ee{1}{l}{b}$,
and $(t_l,u_l)\neq(t_{l'},u_{l'})$ for $l,l'$ with $l\neq l'$, by strong Borel fixedness.
See \cite[\Thm{(1.1)}, \Lem{(4.4)} and \Prop{(4.5)}]{A8} or 
\cite[\Prop{1.5} and \Thm{2.5}]{A7}  for the detail.
Then $\inr^x(I)$ is minimally generated by
$x_1^a,x_1^{a-1}x_2^{\bg_{a-1}}\ddd x_1x_2^{\bg_1},x_2^{\bg_0}$
and  $x_1^{t_l}x_2^{u_l}x_3^{\bg_{t_lu_l}}$ 
$(\ee{1}{l}{b})$.
We have $a=\deg(x_1^a)$, $n^2_l=\deg(x_1^{a-l}x_2^{\bg_{a-l}})$
$(\ee{1}{l}{a})$, and $n^3_l=\deg(x_1^{t_l}x_2^{u_l}x_3^{\bg_{t_lu_l}})$ 
$(\ee{1}{l}{b})$.
\end{rem}

\begin{lem}\lb{proc108}
Let $I$, $R$ and $\seqfour$ be at the beginning of this section, and let 
$\{e^1_1,e^2_1\ddd e^2_a,e^3_1\ddd e^3_b\}$ be a \Wei\ basis of 
$I$ with respect to $\seqfour$.
\begin{parenumr}
\item\lb{c1081}
Let $p$ be an integer with $\ee{1}{p}{a}$ and 
\begin{equation*}
h\in\lt(\ds{\ee{1}{l}{a},\ l\neq p}\ktwo e^2_l\rt)\op\lt(\dss{l=1}{b}\ktwo e^3_l\rt)
\end{equation*}
a homogeneous element such that $\deg(e^2_p)=\deg(h)$. Put 
$e'{}^2_p:=e^2_p+h$. Then, we obtain another \Wei\ basis
$\{e^1_1,e^2_1\ddd e^2_{p-1},e'{}^2_p,e^2_{p+1}\ddd
e^2_a,e^3_1\ddd e^3_b\}$ of $I$ with respect to $\seqfour$.
\item\lb{c1082}
Let $p$ be an integers with $1\leq p\leq a$ such that $\deg(e^1_1)=\deg(e^2_p)$
and $c$ an element of $k$. Put 
$e'{}^1_1:=e^1_1+ce^2_p$. Then,  the system of generators
$\{e'{}^1_1,e^2_1\ddd e^2_a,e^3_1\ddd e^3_b\}$ is another
\Wei\ basis of $I$ with respect to $\seqfour$.
\item\lb{c1084}
Let $p$ be an integers with $1\leq p\leq a$ and $c$ an element of $k^\ast$.
Put $e'{}^2_p:=ce^2_p$. Then, the system of generators 
$\{e^1_1,e^2_1\ddd e^2_{p-1},e'{}^2_p,e^2_{p+1}\ddd
e^2_a,e^3_1\ddd e^3_b\}$
is also a \Wei\ basis of $I$ with respect to $\seqfour$.
\item\lb{c1085}
Let $p,q$ be an integers with $1\leq p<q\leq a$ such that $\deg(e^2_p)=\deg(e^2_q)$. Put 
$e'{}^2_p:=e^2_q$, $e'{}^2_q:=e^2_p$. Then, exchanging $e^2_p$ and $e^2_q$, we obtain another
\Wei\ basis $\{e^1_1,e^2_1\ddd e^2_{p-1},e'{}^2_p,e^2_{p+1}\ddd
e^2_{q-1},e'{}^2_q,e^2_{q+1}\ddd
e^2_a,e^3_1\ddd e^3_b\}$ of $I$ with respect to $\seqfour$.
\item\lb{c1083}
Let $p,q$ be an integers with $1\leq p<q\leq a$ such that $\deg(e^2_l)=\deg(e^2_p)$
for all $\ee{p}{l}{q}$, and let $G\in\GL(q-p+1,k)$. Put 
$(e'{}^2_p\ddd e'{}^2_q):=(e{}^2_p\ddd e{}^2_q)G$. Then, 
we obtain another
\Wei\ basis $\{e^1_1,e^2_1\ddd e^2_{p-1},e'{}^2_p\ddd e'{}^2_q,e^2_{q+1}\ddd
e^2_a,e^3_1\ddd e^3_b\}$ of $I$ with respect to $\seqfour$.
\end{parenumr}
\end{lem}

\begin{pf}
Recall that homogeneous elements $e^1_1,e^2_1,\ldots,e^2_a,e^3_1\ddd e^3_b$ 
of $I$ form a \Wei\ basis of $I$ with respect to $\seqfour$ if and only if 
all of the conditons \eqref{eq40} and \eqref{eq201} hold.
To prove \eqref{c1081}, put $e'{}^1_1:=e^1_1$, $e'{}^2_p:=e^2_p+h$, 
$e'{}^2_l:=e^2_l\ ( l\neq p,\ \ee{1}{l}{a})$,
$e'{}^3_l:=e^3_l\ (\ee{1}{l}{b})$ and let
\begin{align*}
&\supbr{J}{3}:=\smm{l=1}{b}\ktwo e'{}^3_l,\quad 
\supbr{J}{2}:=\smm{l=1}{a}\kthree e'{}^2_l+\supbr{J}{3},\\
&\supbr{J}{1}:=R e'{}^1_1+\supbr{J}{2}.
\end{align*}
We first show that $\supbr{J}{1}$ is an ideal in $R$. Let $l'$ be an 
arbitrary integer with $\ee{1}{l'}{b}$.
Since $\{e^1_1,e^2_1,\ldots,e^2_a,e^3_1\ddd e^3_b\}$ is a \Wei\ basis, we have
\begin{equation*}
\begin{split}
x_2 e'{}^3_{l'}=x_2 e^3_{l'}
&=\smm{l=1}{a}g^2_l e^2_l+\smm{l=1}{b} g^3_l e^3_l\\
&=-g^2_p h+\smm{l=1}{p-1}g^2_l e^2_l+g^2_p e'{}^2_p+\smm{l=p+1}{a}g^2_l e^2_l+
\smm{l=1}{b} g^3_l e^3_l
\end{split}
\end{equation*}
with suitable $g^2_l\in(x_3,x_4)\ktwo$ $(\ee{1}{l}{a})$ and
$g^3_l\in\ktwo$ $(\ee{1}{l}{b})$ by \eqref{eq201}. 
Since $h\in\lt(\ds{\ee{1}{l}{a},\ l\neq p}\ktwo e^2_l\rt)\op\lt(\dss{l=1}{b}\ktwo e^3_l\rt)$,
this implies that
\begin{equation}
x_2e'{}^3_{l'}\in
(\seqtwo)\lt(\smm{l=1}{a}\ktwo e'{}^2_l\rt)+\smm{l=1}{b}\ktwo e'{}^3_l
\sset\supbr{J}{2}.
\lb{eq220}
\end{equation}
Using this repeatedly, we see also that 
\begin{equation}
x_2^te'{}^3_{l'}\in
(\seqtwo)\lt(\smm{l=1}{a}\kthree e'{}^2_l\rt)+\smm{l=1}{b}\ktwo e'{}^3_l
\sset\supbr{J}{2}.
\lb{eq229}
\end{equation}
for all $t\geq1$.
Similarly,
\begin{equation*}
\begin{split}
x_1 e'{}^3_{l'}=x_1 e^3_{l'}
&=g^1_1e^1_1+\smm{l=1}{a}g^2_l e^2_l+\smm{l=1}{b} g^3_l e^3_l\\
&=-g^2_p h+g^1_1e^1_1+\smm{l=1}{p-1}g^2_l e^2_l+g^2_p e'{}^2_p+\smm{l=p+1}{a}g^2_l e^2_l+
\smm{l=1}{b} g^3_l e^3_l
\end{split}
\end{equation*}
with suitable $g^1_1,\ g^2_l\in(x_3,x_4)\kthree$ $(\ee{1}{l}{a})$ and
$g^3_l\in\ktwo$ $(\ee{1}{l}{b})$ by \eqref{eq201}, so that
\begin{equation}
\begin{split}
x_1e'{}^3_{l'}&\in
(\seqtwo)\lt(\kthree e'{}^1_1+\smm{l=1}{a}\kthree e'{}^2_l\rt)+\smm{l=1}{b}\kthree e'{}^3_l\\
&\sset
(\seqtwo)\lt(\kthree e'{}^1_1+\smm{l=1}{a}\kthree e'{}^2_l\rt)+\smm{l=1}{b}\ktwo e'{}^3_l\\
&\sset\supbr{J}{1}
\end{split}
\lb{eq221}
\end{equation}
by \eqref{eq229}.
For each $l'\ (\ee{1}{l'}{a})$, 
\begin{equation*}
\begin{split}
x_1 e'{}^2_{l'}
&=g^1_1e^1_1+\smm{l=1}{a}g^2_l e^2_l+\smm{l=1}{b} g^3_l e^3_l\\
&=-g^2_p h+g^1_1e^1_1+\smm{l=1}{p-1}g^2_l e^2_l+g^2_p e'{}^2_p+\smm{l=p+1}{a}g^2_l e^2_l+
\smm{l=1}{b} g^3_l e^3_l
\end{split}
\end{equation*}
with suitable $g^1_1\in(\seqthree)\kthree$, $g^2_l\in\kthree$ $(\ee{1}{l}{a})$ and
$g^3_l\in\ktwo$ $(\ee{1}{l}{b})$ again by \eqref{eq201}, whether $l'=p$ or not. Hence
\begin{equation}
\begin{split}
x_1e'{}^2_{l'}&\in
(\seqthree)\kthree e'{}^1_1+\smm{l=1}{a}\kthree e'{}^2_l+\smm{l=1}{b}\ktwo e'{}^3_l\\
&\sset\supbr{J}{1}
\end{split}
\lb{eq222}
\end{equation}
by \eqref{eq229}.
Now, we find by \cite[\Lem{2.7}]{A5} that $\supbr{J}{1}$ is an ideal
in $R$, since we have verified \eqref{eq220}, \eqref{eq221} and \eqref{eq222}.
Moreover, it is clear that the subideal $\supbr{J}{1}$ contains the generators
$e^1_1,e^2_1,\ldots,e^2_a,e^3_1\ddd e^3_b$ of $I$. Hence $\supbr{J}{1}=I$.
Since $\deg(e'{}^i_l)=\deg(e^i_l)$ for all $i,l$, we find by counting 
$\dimm{k}{\gradn{I}{t}}$ $(t\in\Zbf)$ with the use of \eqref{eq40} and 
\begin{equation*}
I=\supbr{J}{1}=R e'{}^1_1+\smm{l=1}{a}\kthree e'{}^2_l+\smm{l=1}{b}\ktwo e'{}^3_l
\end{equation*}
that the above expression must be a direct sum as $\ktwo$-module. Then, the conditions 
\eqref{eq220}, \eqref{eq221} and \eqref{eq222}  imply that 
$\{e'{}^1_1,e'{}^2_1\ddd e'{}^2_a,e'{}^3_1\ddd e'{}^3_b\}$ satisfies
\eqref{eq201} with $e^i_l$ replaced by $e'{}^i_l$. Thus 
$\{e'{}^1_1,e'{}^2_1\ddd e'{}^2_a,e'{}^3_1\ddd e'{}^3_b\}$ is a \Wei\ basis
of $I$ with respect to $\seqfour$. This proves \eqref{c1081}.
One can prove \eqref{c1082} in the same manner. The assertions \eqref{c1084}
and \eqref{c1085} are trivial. Since an element of $\GL(q-p+1,k)$ is a product of 
elementary matrices and a diagonal one, we obtain \eqref{c1083} by
\eqref{c1081} and \eqref{c1084}.
\end{pf}

\par
Let us recall the standard free resolution of $I$ starting with 
a \Wei{} basis 
$\{e^1_1,e^2_1\ddd e^2_a,e^3_1\ddd e^3_b\}$. 
For each triple $i,i',l$ with $1\leq i<i'\leq 3,\ \ee{1}{l'}{m_{i'}}$
($m_2=a,\ m_3=b$), we have 
\begin{equation}
x_i e^{i'}_{l'}=g^1_1 e^1_1+\smm{l=1}{a}g^2_l e^2_l
+\smm{l=1}{b} g^3_l e^3_l
\lb{eq202}\end{equation}
by \eqref{eq40}, where 
$g^1_1\in R,\ g^2_l\in \kthree,\ g^3_l\in \ktwo$.
Put
\begin{equation}
\lg_1:=(e^1_1,e^2_1,\ldots,e^2_a,e^3_1\ddd e^3_b).
\lb{eq205}
\end{equation}
Since \eqref{eq201} also holds, there are matirces with
homogeneous components
\begin{equation}
\lt\{\begin{aligned}
{}&U_{01}\in\MATr((x_2,x_3,x_4)\kthree),\ 
\Uovercc\kernsub{1}\in\MATr(\kthree),\\
{}&U_{02},\ U_2\in\MATr((x_3,x_4)\kthree),
\\
{}&U_{21}, \Uovercc\kernsub{3}, \Uovercc\kernsub{5}\in\MATr(\ktwo),\ 
U_4\in\MATr((x_3,x_4)\ktwo)
\end{aligned}\rt.
\lb{eq32}
\end{equation}
such that
\begin{equation}
\lt\{\begin{aligned}
{}&(e^2_1\ddd e^2_a)x_11_a=
\lg_1\begin{bmatrix}
-U_{01}\\
\Uovercc\kernsub{1}\\
-U_{21}
\end{bmatrix},\\
{}&(e^3_1\ddd e^3_b)x_11_b=
\lg_1\begin{bmatrix}
-U_{02}\\
-U_2\\
\Uovercc\kernsub{3}
\end{bmatrix},\quad
(e^3_1\ddd e^3_b)x_21_b=
\lg_1\begin{bmatrix}
0\\
-U_4\\
\Uovercc\kernsub{5}
\end{bmatrix}.
\end{aligned}\rt.
\lb{eq203}
\end{equation}
Let
\begin{equation}
\lt\{\begin{aligned}
{}&U_1:=x_11_a-\Uovercc\kernsub{1},\quad
U_3:=x_11_b-\Uovercc\kernsub{3},\quad
U_5:=x_21_b-\Uovercc\kernsub{5},\\
{}&\lg_2:=\begin{bmatrix}
U_{01}&U_{02}&0\\
U_1&U_2&U_4\\
U_{21}&U_3&U_5
\end{bmatrix},\quad
\lg_3:=\begin{bmatrix}
-U_4\\
-U_5\\
U_3
\end{bmatrix}.
\end{aligned}\rt.
\lb{eq59}
\end{equation}
Then 
\begin{equation}
U_1-x_11_a\in\MATr(\kthree),\quad 
U_3-x_11_b,\ U_5-x_21_b\in\MATr(\ktwo)
\lb{eq204}
\end{equation}
by \eqref{eq32} and
\begin{equation*}
\lg_1\lg_2=0
\end{equation*}
by \eqref{eq203}. Furthermore,
\begin{equation*}
\lg_2\lg_3=0
\end{equation*}
and the sequence
\begin{equation}
\begin{split}
0\lra R(-\tgbar-2)\xra{\ \lg_3\ } 
&R(-\sgbar-1,-\tgbar-1,-\tgbar-1)\\
\xra{\ \lg_2\ }
&R(-a,-\sgbar,-\tgbar)\xra{\ \lg_1\ }I\lra0
\end{split}
\lb{eq35}
\end{equation}
is exact by  \cite[\Ex{2.8}]{A2} or \cite[(3) of \Cor{(3.11)}]{A8},
where $\sgbar:=(\seq{n}{1}{a})$, $\tgbar:=(\seq{n}{a+1}{a+b})$.
The relation $\lg_2\lg_3=0$ embraces a lot of informations on the
structure of $\Ext{3}{R}{R/I}{R}$.

\begin{lem}\lb{proc200}
Let $I$ be a homogeneous ideal in $R$ of height 2 such 
that $\depm{R/I}\geq1$, $\seqfour\in\gradn{R}{1}$ linear forms,
$\{e^1_1,e^2_1\ddd e^2_a,e^3_1\ddd e^3_b\}$ a \Wei\ basis of $I$ 
with respect to $\seqfour$. Let further $\lg_2$, $\lg_3$, $U_i\ (\ee{1}{i}{5})$,
$U_{01},\ U_{02},\ U_{21}$ be matrices
as in \eqref{eq32} -- \eqref{eq59}. Then
\begin{equation*}
\begin{split}
&\Ext{3}{R}{R/I}{R}\cong\Coksup{R}{\tpose{U_3},\tpose{U_5},\tpose{U_4}}\\
&\cong\lt(\dss{i=1}{b}\ktwo(n_{a+i}+2)\rt)/\lt(
\smm{\rg=0}{b-1}\Imasup{\ktwo}{(\tpose{\Uovercc\kernsub{5}})^\rg\cdot\tpose{U_4}}
\rt),
\end{split}
\end{equation*}
where the second isomorphism stands for an isomorphism over $\ktwo$.
Moreover, 
$x_1v\in(\seqtwo)\Imasup{R}{\tpose{U_3},\tpose{U_5}}+\Imasup{\kthree}{\tpose{U_4}}$
for all columns $v$ of $\tpose{U_4}$.
\end{lem}

\begin{pf}
One can find a proof of the major part in \cite[pp. 802--803]{A3}. But, for the 
convenience of the readers, we reproduce a proof here.
Now, let us consider $\Imasup{R}{\tpose{U_3},\tpose{U_5},\tpose{U_4}}$.
Forgetting degrees, we have
\begin{equation*}
R^{a+2b}=\Imr^R\lt(
\begin{bmatrix}
-\tpose{U_5}&-\tpose{U_4}\\
\tpose{U_3}&\tpose{U_2}\\
\tpose{U_{21}}&\tpose{U_1}\\
\end{bmatrix}
\rt)
\op\lt(R^{b}\opdot\kthree^{a+b}\rt).
\end{equation*}
by \cite[\Lem{1.1}]{A5} and 
\begin{equation*}
(\tpose{U_3},\tpose{U_5},\tpose{U_4})
\begin{bmatrix}
-\tpose{U_5}&-\tpose{U_4}\\
\tpose{U_3}&\tpose{U_2}\\
\tpose{U_{21}}&\tpose{U_1}\\
\end{bmatrix}
=0
\end{equation*}
by the relation $\lg_2\lg_3=0$, so that 
\begin{equation*}
\Imasup{R}{\tpose{U_3},\tpose{U_5},\tpose{U_4}}
=\Imasup{R}{\tpose{U_3}}\op\Imasup{\kthree}{\tpose{U_5},\tpose{U_4}}
\end{equation*}
as $\kthree$-module. At the same time, since 
$-\tpose{U_3}\tpose{U_4}+\tpose{U_5}\tpose{U_2}+
x_1\tpose{U_4}-\tpose{U_4}\tpose{\Uovercc\kernsub{1}}=0$, one finds that each
column of $x_1\tpose{U_4}$ lies in 
$(\seqtwo)\Imasup{R}{\tpose{U_3},\tpose{U_5}}+\Imasup{\kthree}{\tpose{U_4}}$ 
by \eqref{eq32}. This proves the last assertion.
Since 
\begin{align*}
x_2^\ng\cdot v
&=\tpose{U_5}\lt(\smm{i=1}{\ng}x_2^{\ng-i}
(\tpose{\Uovercc}\kernsub{5})^{i-1}\rt)v
+(\tpose{\Uovercc\kernsub{5}})^\ng\cdot v
\\
&\in \Imasup{\kthree}{\tpose{U_5}}\op
\sm{\rg\geq0}\Imasup{\ktwo}{(\tpose{\Uovercc\kernsub{5}})^\rg(\tpose{U_4}}),
\notag
\end{align*}
for all $\ng\in\Zbf$ and for each column $v$ of $\tpose{U_4}$,
we see that 
\begin{equation*}
\Imasup{\kthree}{\tpose{U_5},\tpose{U_4}}=
\Imasup{\kthree}{\tpose{U_5}}\op
\sm{\rg\geq0}\Imasup{\ktwo}{(\tpose{\Uovercc\kernsub{5}})^\rg(\tpose{U_4}})
\end{equation*}
as $\ktwo$-module. Moreover 
\begin{equation*}
\sm{\rg\geq0}\Imasup{\ktwo}{(\tpose{\Uovercc\kernsub{5}})^\rg(\tpose{U_4}})
=\smm{\rg=0}{b-1}\Imasup{\ktwo}{(\tpose{\Uovercc\kernsub{5}})^\rg(\tpose{U_4}})
\end{equation*}
by Hamilton-Cayley's theorem. On the other hand,
\begin{equation*}
R^b=\Imasup{R}{\tpose{U_3}}\op\Imasup{\kthree}{\tpose{U_5}}\op
\ktwo^b
\end{equation*}
as $\ktwo$-module. Hence
\begin{equation*}
\Coksup{R}{\tpose{U_3},\tpose{U_5},\tpose{U_4}}
\cong\ktwo^b/\lt(
\smm{\rg=0}{b-1}\Imasup{\ktwo}{(\tpose{\Uovercc\kernsub{5}})^\rg\cdot\tpose{U_4}}
\rt)
\end{equation*}
as $\ktwo$-module.  
Since $\Ext{3}{R}{R/I}{R}\cong\Coksup{R}{\tpose{U_3},\tpose{U_5},\tpose{U_4}}$
by the free resolution \eqref{eq35}, we get our first assertion, taking the degrees
into account.
\end{pf}

\begin{lem}\lb{proc32}
Let $I=(f,g)$ be a homogeneous ideal in $R$ generated by 
homogeneous polynomials $f,\ g$ of degree $p,\ q$ respectively.
Suppose that $p\leq q$ and that $f,\ g$ form an $R$-regular sequence.
Then the basic sequence of $I$ is $(p;q,q+1\ddd q+p-1)$.
\end{lem}

\begin{pf}
Let $\{e^1_1,e^2_1\ddd e^2_a,e^3_1\ddd e^3_b\}$ be a \Wei\ basis of $I$,
$\lg_2$ tha matrix as in \eqref{eq59}, and $(a;\seq{n}{1}{a};\seq{n}{a+1}{a+b})$
the basic sequence of $I$.
It is clear that $a=p$. 
Since $R/I$ is \CM, we have $b=0$. Hence
\begin{equation*}
\lg_2=\begin{bmatrix}
U_{01}\\
U_1
\end{bmatrix}.
\end{equation*}
Since $I$ is minimally generated over $R$ by two elements, the rank of 
the relation matrix $\lg_2\ (\mod\ (\seqfour))$ must be $a-1$. On the other hand
$U_{01}\equiv0\  (\mod\ (\seqfour))$ by \eqref{eq32}, so that
the rank of $U_1\ (\mod\ (\seqfour))$ must be $a-1$.
Since $\seq{n}{1}{a}$ is a nondecreasing sequence, this is possible only 
when $n_l=n_1+l-1$ for all
$\ee{1}{l}{a}$, considering the degrees of the components of
$U_1$ (cf. \eqref{eq202}, \eqref{eq59}). Moreover,
\begin{equation*}
\rankr_k\lt(U_1\clmlt{1}{\ }\ (\mod\ (\seqfour))\rt)=a-1
\end{equation*}
in this case. Hence $I=(e^1_1,e^2_1)$. From this it follows that
$q=\deg(e^2_1)=n_1$. This proves our assertion.
\end{pf}

\section{\mathversion{bold}Properties of $U_{01}$ and $U_1$}
\lb{sec2}

Let $I$ be a homogeneous ideal in $R=k[y_1,y_2,y_3,y_4]$ of height 2 such 
that $\depm{R/I}\geq1$ and let $\BR{I}=(a;n_1\ddd n_a;n_{a+1}\ddd n_{a+b})$
be its basic sequence. It was proved for the first time in \cite[Corollaire 2.2]{GP2} that
$n_1\ddd n_a$ is connected (i.e. $n_i\leq n_{i+1}\leq n_i+1$ for all
$\ee{1}{i}{a-1}$) if $I$ is prime and $b=0$. Later, the same result was proved
without the assumption that $b=0$ in \cite[\Cor{1.2}]{A1}. In recent papers
\cite{C} and \cite{DS}, enhanced assertions similar to the above connectedness
are given for the case $\ch{k}=0$, considering the generic initial ideal of $I$.
But there seems to be some gaps in the proofs given there. Along the same line
of argument as that of \cite{A1}, we gave a proof of the connectedness 
asserted by Cook for a special case in \cite{A10}.
In this section, we first give a brief but a little bit closer 
account of our method to prove the connectedness of $n_1\ddd n_a$.
Developing our argument further, we next show some crucial properties of 
$U_{01}$ and $U_1$ needed to prove our main theorems, when there is an
irreducible homogeneous polynomial of degree $a$ in $I$.

\par
Let
$\zg_{ij}\ (\ee{1}{i}{4},\ \ee{1}{j}{4})$ be indeterminates over
$R$, $K$ the quotient field of the polynomial 
ring $k[\zg]:=k[\zg_{ij}; \ee{1}{i}{4},\ \ee{1}{j}{4}]$, $z_1,z_2,z_3,z_4$ 
elements of $R_K:=K[y_1,y_2,y_3,y_4]$ such that 
$y_i=\smm{j=1}{4}\zg_{ji}z_j$ $(\ee{1}{i}{4})$, and $I_K:=IR_K$. 
Then, there is a
\Wei{} basis $\{\etild^1_1,\etild^2_1\ddd \etild^2_a,\etild^3_1\ddd \etild^3_b\}$ 
of $I_K$ with respect to $z_1,z_2,z_3,z_4$, since $\dep{\mfrak R_K}{R_K/I_K}\geq1$
(see \cite[\Thm{2.12}]{A5}, \cite[\Thm{2.5}]{A7}).
The homogeneous polynomials $\etild^i_l\ (\ee{1}{i}{3},\ \ee{1}{l}{m_i})\in R_K$
satisfy \eqref{eq200} -- \eqref{eq201} with $I,\ x_i,\ k,\ e^i_l$ replaced by
$I_K,\ z_i,\ K,\ \etild^i_l$ respectively.
By exactly the same method as in the previous section, one obtains 
matrices
\begin{equation*}
\lg_1:=(\etild^1_1,\etild^2_1\ddd \etild^2_a,\etild^3_1\ddd \etild^3_b),\quad
\lgtild_2:=\begin{bmatrix}
\Utild_{01}&\Utild_{02}&0\\
\Utild_1&\Utild_2&\Utild_4\\
\Utild_{21}&\Utild_3&\Utild_5
\end{bmatrix},\quad
\lgtild_3:=\begin{bmatrix}
-\Utild_4\\
-\Utild_5\\
\Utild_3
\end{bmatrix}
\end{equation*}
which give a free resolution
\begin{equation*}
\begin{split}
0\lra R_K(-\tgbar-2)\xra{\ \lgtild_3\ } 
&R_K(-\sgbar-1,-\tgbar-1,-\tgbar-1)\\
\xra{\ \lgtild_2\ }
&R_K(-a,-\sgbar,-\tgbar)\xra{\ \lgtild_1\ }I_K\lra0.
\end{split}
\end{equation*}
Note that these matrices satisfy the conditions corresponding to
\eqref{eq32} and \eqref{eq204} with $x_i,\ k$ replaced by $z_i,\ K$.

\begin{rem}\lb{proc123}
The results on \Wei\ bases described in Section 1 are all valid for 
$\{\etild^1_1,\etild^2_1\ddd \etild^2_a,\etild^3_1\ddd \etild^3_b\}$,
since it is a \Wei\ basis in any case.
\end{rem}

\par
Assume that the matrix 
$\Gg=(\gg_{ij})$ mentioned at the beginning of Section 1 is sufficiently 
general. 
For each pair $i,l$, let
$e^i_l$ denote the polynomial obtained from 
$\etild^i_l$ by the substitution $(\zg_{ij})=\Gg$. 
Then $\{e^1_1,e^2_1\ddd e^2_a,e^3_1\ddd e^3_b\}$
is a \Wei\ basis
of $I$ with respect to $\seqfour$
(see the proof of \cite[\Thm{2.12}]{A5}).
In this case, the  matirices $\lg_2$, $\lg_3$, $U_i\ (\ee{1}{i}{5})$, $U_{01},\ U_{02}$ and 
$U_{21}$ described in Section \ref{sec1} are also
obtained from $\lgtild_2$, $\lgtild_3$, $\Utild_i\ (\ee{1}{i}{5})$, $\Utild_{01},\ \Utild_{02}$ and 
$\Utild_{21}$ respectively by the same substitution.

\begin{lem}\lb{proc106}
With the notation above, suppose that there is an irreducible 
homogeneous polynomial of degree $a$ in $I$.  
Denote by $\echeck{}^i_l,\ \Ucheck_{01},\ \Ucheck_1$ and so on
the polynomials and the matrices obtained 
from $\etild^i_l,\ \Utild_{01},\ \Utild_1$ and so on by the substitution $z_3=z_4=0$.
Suppose there is an integer $s'$ ($\nn{0}{s'}{a}$) such that  
\begin{equation*}
\begin{bmatrix}
\Ucheck_{01}\\
\Ucheck_1
\end{bmatrix}=
\begin{bmatrix}
\Dcheck_{11}&\Dcheck_{12}\\
0&\Dcheck_{22}
\end{bmatrix}
\end{equation*}
with an $(s'+1)\tm s'$ matrix $\Dcheck_{11}$, an $(s'+1)\tm s''$ matrix 
$\Dcheck_{12}$
and an $s''\tm s''$ matrix $\Dcheck_{22}$, where $s''=a-s'$. 
Then $n_{s'+1}=a$. 
\end{lem}

\begin{pf}
Since $\Utild_{02},\ \Utild_2\in\MATr((z_3,z_4)\Kthreez)$ by the condition
corresponding to \eqref{eq32}, it follows from the relation $\lgtild_1\lgtild_2=0$
that $(\echeck{}^3_1\ddd\echeck^3_b)\Ucheck_3=0$.
This implies that $\echeck^3_l=0$ for all $\ee{1}{l}{b}$.
Hence
\begin{equation*}
(\echeck^1_1,\echeck^2_1\ddd\echeck^2_a)
\begin{bmatrix}
\Ucheck_{01}\\
\Ucheck_1
\end{bmatrix}=0
\lb{eq38}
\end{equation*}
again by the relation $\lgtild_1\lgtild_2=0$.
Let $\Icheck_{K}:=(\echeck^1_1,\echeck^2_1\ddd\echeck^2_a)K[z_1,z_2]\sset K[z_1,z_2]$.
Then $R_K/I_K+(z_3,z_4)R_K=K[z_1,z_2]/\Icheck_K$ and this ring must be of
finite length over $K$. By Hilbert-Burch theorem,
\begin{equation*}
\echeck^1_1=\det\lt(
\Ucheck_1\rt),\quad
\echeck^2_i=(-1)^i\det\lt(
\begin{bmatrix}
\Ucheck_{01}\\
\Ucheck_1
\end{bmatrix}\clml{i+1}{\ }\rt)
\quad(\ee{1}{i}{a})
\end{equation*}
up to constant multiplication.
Since there is an irreducible homogeneous polynomial of degree $a$ in $I$  
by hypothesis, $\Icheck_K$ also contains an irreducible homogeneous polynomial
of degree $a$, say $\echeck$, by \Lem{\ref{proc25}} below. 
Let $i_0:=0$ if $n_1>a$ and $i_0:=\max\{\ i\ |\ n_i=a,\ \ee{1}{i}{a}\ \}$ otherwise. 
We have
$\echeck=\ccheck_0\echeck^1_1+\smm{i=1}{i_0}\ccheck_i\echeck^2_i$ 
with $\ccheck_0,\ccheck_i\in K$, where we understand
$\smm{i=1}{i_0}\ccheck_i\echeck^2_i=0$
if $i_0=0$. If $n_{s'+1}>a$, then $i_0\leq s'$, so that the above formula implies thta
$\echeck$ must be divisible by $\det(\Dcheck_{22})$, which is a contradiction.
Hence $n_{s'+1}=a$. In consequence, $a=n_1=\cdots=n_{s'+1}$.
\end{pf}

\begin{lem}\lb{proc25}
Let $f(y_1,y_2,y_3,y_4)$ be an irreducible homogeneous polynomial in
$R$. Then 
\begin{equation*}
\fhat:=f(\zg_{11}z_1+\zg_{21}z_2,\zg_{12}z_1+\zg_{22}z_2,
\zg_{13}z_1+\zg_{23}z_2,\zg_{14}z_1+\zg_{24}z_2)
\end{equation*}
is irreducible in $k[\zg][z_1,z_2]$. 
\end{lem}

\begin{pf}
Suppose to the contrary and let $\fhat_1,\ \fhat_2\in k[\zg][z_1,z_2]$ be
polynomials of positive degree such that 
\begin{equation}
\fhat=\fhat_1\fhat_2.
\lb{eq27}
\end{equation}
Since $\fhat$ is homogeneous of degree $\deg(f)$ in 
$\zg_{ij}\ (\ee{1}{i}{4},\ \ee{1}{j}{4})$ and in $z_1,z_2$, we see that 
$\fhat_1$ and $\fhat_2$ are also homogeneous in 
$\zg_{ij}\ (\ee{1}{i}{4},\ \ee{1}{j}{4})$ and in $z_1,z_2$.
Let $\dg_1$ (resp. $\dg_2$) be the degree of $\fhat_1$ (resp. $\fhat_2$)
in $z_1,z_2$. Consider the equality \eqref{eq27} in $k[\zg][z_1,z_2]_{z_1}$.
Since $\zg_{1i}+\zg_{2i}(z_2/z_1)$ $(\ee{1}{i}{4})$ and  
$\zg_{ji}\ (\ee{1}{i}{4},\ \ee{2}{j}{4})$ are
algebraically independent over $k[z_1,z_2]_{z_1}$, and since
\begin{equation*}
\fhat=z_1^{\deg(f)}f(
\zg_{11}+\zg_{21}(z_2/z_1),\zg_{12}+\zg_{22}(z_2/z_1),
\zg_{13}+\zg_{23}(z_2/z_1),\zg_{14}+\zg_{24}(z_2/z_1)),
\end{equation*}
we see that $\fhat_1$ or $\fhat_2$ must be a unit in 
$k[\zg][z_1,z_2]_{z_1}$ by the irreducibility of $f$. 
We may assume without any loss of generality
that $\fhat_1=c_1z_1^{\dg_1}$ $(c_1\in k)$.
Likewise, considering \eqref{eq27} in $k[\zg][z_1,z_2]_{z_2}$,
we find that $\fhat_1$ or $\fhat_2$ must be a unit in 
$k[\zg][z_1,z_2]_{z_2}$. In consequence,  
$\fhat_2=c_2z_2^{\dg_2}$ $(c_2\in k)$.
Thus, $\fhat$ is of degree zero in $\zg_{ij}\ (\ee{1}{i}{4},\ \ee{1}{j}{4})$,
which is absurd.
\end{pf}

With the notation above, let 
$n'_1\ddd n'_\og$ $(\og\geq1)$ be a strictly increasing sequence of integers
such that $\{n'_1\ddd n'_\og\}=\{\ n_i\ |\ \ee{1}{i}{a}\ \}$ and 
let $t_l:=\#\{\ i\ |\ n_i=n'_l,\ \ee{1}{i}{a}\ \}$. When $\og>1$, we have
\begin{equation*}
\begin{bmatrix}
\Utild_{01}\\
\Utild_1
\end{bmatrix}=
\begin{bmatrix}
\ast&\ast&\hdotsfor{2}&\ast\\
\Dtild_1&\ast&\hdotsfor{2}&\ast\\
\Ctild_1&\Dtild_2&\ast&\hdotsfor{1}&\ast\\
0&\Ctild_2&\Dtild_3&\ast&\ast\\
\vdots&\ddots&\ddots&\ddots&\ast\\
0&\dots&0&\Ctild_{\og-1}&\Dtild_\og
\end{bmatrix},
\lb{eq36}
\end{equation*}
where $\Ctild_l$ is a $t_{l+1}\tm t_l$ matrix whose entries are
homogeneous of degree $n'_l+1-n'_{l+1}\leq0$ for each $\en{1}{l}{\og}$
and $\Dtild_l$ is a $t_l\tm t_l$ matrix whose entries are linear forms
in $z_1,z_2,z_3,z_4$ over $K$ for each $\ee{1}{l}{\og}$.

\begin{lem}\lb{proc26}
Let the notation be as above. Suppose that there is an irreducible 
homogeneous polynomial of degree $a$ in $I$. Then, 
$n_i\leq n_{i+1}\leq n_i+1$ for all $\ee{1}{i}{a-1}$. If further $\og>1$, then
$\Ctild_l\neq0$ for all $\en{1}{l}{\og}$.
\end{lem}

\begin{pf}
When $\og=1$, our assertion is
trivial. Assume $\og>1$. 
Suppose $\Ctild_{\og_0}=0$ for some $\en{1}{\og_0}{\og}$.
Then
\begin{equation*}
\begin{bmatrix}
\Utild_{01}\\
\Utild_1
\end{bmatrix}=
\begin{bmatrix}
\Dtild_{11}&\Dtild_{12}\\
0&\Dtild_{22}
\end{bmatrix}
\lb{eq37}
\end{equation*}
with an $(s'+1)\tm s'$ matrix $\Dtild_{11}$, an $(s'+1)\tm s''$ matrix 
$\Dtild_{12}$
and an $s''\tm s''$ matrix $\Dtild_{22}$, where $s'=\smm{l=1}{\og_0}t_l>0$,
$s''=a-s'>0$. Denote by $\Ucheck_{02},\ \Ucheck_1,\ \Dcheck_{ij}$ the matrices obtained 
from $\Utild_{02},\ \Utild_1,\ \Dtild_{ij}$ respectively by the substitution $z_3=z_4=0$.
Then,
\begin{equation*}
\begin{bmatrix}
\Ucheck_{01}\\
\Ucheck_1
\end{bmatrix}=
\begin{bmatrix}
\Dcheck_{11}&\Dcheck_{12}\\
0&\Dcheck_{22}
\end{bmatrix}.
\end{equation*}
But $n_{s'+1}=n'_{\og_0+1}>a$ by our choice of $\og_0$ and $s'$.
This contradicts \Lem{\ref{proc106}}.
Thus we have shown that $\Ctild_l\neq0$ for all $\en{1}{l}{\og}$.
In consequence, $n'_{l+1}=n'_l+1$ for all $\en{1}{l}{\og}$
and $\Ctild_l$ is a matrix in $K$ different from zero.
Hence, $n_i\leq n_{i+1}\leq n_i+1$ for all $\ee{1}{i}{a-1}$.
\end{pf}

\begin{lem}\lb{proc115}
Let the notation and the assumption be the same as in \Lem{\ref{proc26}}.
Suppose $\og>1$ and let $m$ be an ingeger with $\en{1}{m}{\og}$.
Then, there is a $\Gtild_m\in\GL(t_m,K)$ such that
\begin{align*}
&{\begin{bmatrix}
\Gtild_m^{-1}&0\\
0&1_{t_{m+1}}
\end{bmatrix}
\begin{bmatrix}
\Dtild_m\\
\Ctild_m
\end{bmatrix}
\Gtild_m
=\begin{bmatrix}
\Dtild'_m&\Dtild''_m\\
0&\Ctild''_m
\end{bmatrix},}\\
&(\Dtild'_m,\Dtild''_m)-z_11_{t_m}\in\MATr(\Kthreez),\ \andt\
\Ctild''_m\in\MATr(K),
\end{align*}
where the ranks of $\Ctild''_m$ and $\Ctild_m$ are the same and coincide with the number of 
the columns of $\Ctild''_m$.
\end{lem}

\begin{pf}
It is clear that there is a $\Gtild_m\in\GL(t_m,K)$ satisfying $\Ctild_m\Gtild_m=(0,\Ctild''_m)$ 
with a matrix $\Ctild''_m\in\MATr(K)$ such that the ranks of $\Ctild''_m$ and $\Ctild_m$ are 
the same and coincide with the number of the columns of $\Ctild''_m$. Using $\Gtild_m$, 
we have the desired relations, 
since $\Dtild_m-z_11_{t_m}\in\MATr(\Kthreez)$. This proves our assertion.
\end{pf}

\begin{lem}\lb{proc114}
Let the notation and the assumption be the same as in \Lem{\ref{proc106}}
and let $q$ be a positive integer with $q<a$. 
Put $i_0:=0$ if $n_1>a$ and $i_0:=\max\{\ i\ |\ n_i=a,\ \ee{1}{i}{a}\ \}$ otherwise. 
Then, at least one of the
following three cases occurs.
\begin{parenumr}
\item\lb{c1141}
The height of the ideal in $K[z_1,z_2]$ generated by the maximal minors of 
the matrix $\Ucheck_1\clmsc{\phantom{p}}{q+1\ddd a}$ is greater than one.
\item\lb{c1142}
We have $n_1=a$ and there is a $\Gtild\in\GL(a+1,K)$ representing a row operation 
of $\bclmsc{\Ucheck_{01}}{\Ucheck_1}$ 
which adds multiples of $\Ucheck_{01}$ to the rows of $\Ucheck_1$ such that 
the height of the ideal in $K[z_1,z_2]$ generated by the maximal minors of 
$\lt(\Gtild\bclmsc{\Ucheck_{01}}{\Ucheck_1}\rt)\clmlsc{1}{q+1\ddd a}$
is greater than one, where 
$\Gtild=
\begin{smallbmatrix}
\Gtild'&0\\
0&1_{a-i_0}
\end{smallbmatrix}$
with $\Gtild'\in\GL(i_0+1,K)$.
\item\lb{c1143}
We have $n_1=a$ and there is a $\Gtild\in\GL(a,K)$ such that 
the height of the ideal in $K[z_1,z_2]$ generated by the maximal minors of 
$\bclmsc{\Gtild\Ucheck_1}{\Ucheck_{01}}\clmlsc{1}{q+1\ddd a}$
is greater than one, where 
$\Gtild=
\begin{smallbmatrix}
\Gtild'&0\\
0&1_{a-i_0}
\end{smallbmatrix}$
with $\Gtild'\in\GL(i_0,K)$.
\end{parenumr}
\end{lem}

\begin{pf}
Consider first the case $n_1>a$. 
In this case, as the only polynomial in $\Icheck_K$ of degree $a$
up to constant factors, the  $\echeck^1_1=\det(\Ucheck_1)$ 
must be irreducible by hypothesis and \Lem{\ref{proc25}}.
If the height of the ideal generated by the maximal minors of
$\Ucheck_1\clmsc{\phantom{p}}{q+1\ddd a}$ were one, then 
$\det(\Ucheck_1)$ would have to be divisible by a homogeneous 
polynomial of degree between 1 and $a-1$, in contradiction with the irreducibility of 
$\echeck^1_1$. Hence we have \eqref{c1141}.
Next assume $n_1=a$. Then, the ideal $\Icheck_K$ contains a linear combination 
$\echeck=
\ccheck_0\echeck^1_1+\ccheck_1\echeck^2_1+\cdots+\ccheck_{i_0}\echeck^2_{i_0}$ 
$(\ccheck_i\in K,\ \ee{0}{i}{i_0})$ which is irreducible, 
since there is an irreducible homogeneous polynomial of
degree $a$ in $\Icheck_K$ by \Lem{\ref{proc25}}. 
Put 
$\Gtild:=\begin{smallbmatrix}
\Gtild'&0\\
0&1_{a-i_0}
\end{smallbmatrix}$ with
\begin{equation*}
\Gtild':=\begin{bmatrix}
1&0\\
{\begin{matrix}
\ccheck_1/\ccheck_0\\
\vdots\\
\ccheck_{i_0}/\ccheck_0
\end{matrix}}&\text{\Large \ $1_{i_0}$\ }
\end{bmatrix}^{-1}
=\begin{bmatrix}
1&0\\
{\begin{matrix}
-\ccheck_1/\ccheck_0\\
\vdots\\
-\ccheck_{i_0}/\ccheck_0
\end{matrix}}&\text{\Large \ $1_{i_0}$\ }
\end{bmatrix}
\end{equation*}
if $\ccheck_0\neq0$, and
\begin{equation*}
\Gtild':=\begin{bmatrix}
\ccheck_1&\ccheck_{12}&\dots&\ccheck_{1i_0}\\
\vdots&\vdots&&\vdots\\
\ccheck_{i_0}&\ccheck_{i_02}&\dots&\ccheck_{i_0i_0}
\end{bmatrix}^{-1}
\end{equation*}
otherwise, where $\ccheck_{ij}\in K$ $(\ee{1}{i}{i_0},\ \ee{2}{j}{i_0})$ are chosen
so that the inverse matrix exists. 
Then 
\begin{equation*}
(\echeck^1_1,\echeck^2_1\ddd\echeck^2_a)\Gtild^{-1}
\lt(\Gtild\bclm{\Ucheck_{01}}{\Ucheck_1}\rt)=0,\quad
(\echeck^1_1,\echeck^2_1\ddd\echeck^2_a)\Gtild^{-1}=
(\echeck/\ccheck_0,\ldots)
\end{equation*}
or
\begin{equation*}
(\echeck^2_1\ddd\echeck^2_a,\echeck^1_1)
\begin{bmatrix}
\Gtild^{-1}&0\\
0&1
\end{bmatrix}
\bclm{\Gtild\Ucheck_1}{\Ucheck_{01}}=0,\quad
(\echeck^2_1\ddd\echeck^2_a,\echeck^1_1)
\begin{bmatrix}
\Gtild^{-1}&0\\
0&1
\end{bmatrix}=
(\echeck,\ldots)
\end{equation*}
accordingly. By Hilbert-Burch theorem, we must have
$\echeck/\ccheck_0=\det\lt(\Gtild\bclmsc{\Ucheck_{01}}{\Ucheck_1}\clmsc{1}{\ }\rt)$
if $\ccheck_0\neq0$, and 
$\echeck=
\det\lt(\bclmsc{\Gtild\Ucheck_1}{\Ucheck_{01}}\clmlsc{1}{\ }\rt)$
otherwise. Hence the case \eqref{c1142} or \eqref{c1143} occurs for the same reason
as in the case $n_1>a$.
\end{pf}

\begin{lem}\lb{proc121}
Let $p,q$ be positive integers with $p>q$, $\seq{\ng}{1}{p}$ a sequence of 
integers, and $U$ a $p\tm q$ matrix with components in $k[x_1,x_2]$ such that
its columns are homogeneous elements of $\dss{i=1}{p}k[x_1,x_2](-\ng_i)$ whose 
degrees form a nondecreasing sequence $\seq{\mg}{1}{q}$. 
Suppose that $\mg_i=\ng_i+1$ for all $\ee{1}{i}{q}$, the $(i,j)$ component of $U$
lies in $k[x_2]$ for all $i>q$, and that 
the height of the ideal in $k[x_1,x_2]$ generated 
by the maximal minors of $U$ is greater than one.
Suppose further that
\begin{align}
&U-\bclm{x_11_q}{0}\in\MATr(k[x_2])\quad \ort
\lb{eq242}\\
&U-
{\begin{bmatrix}
0&0\\
0&x_11_{q-1}\\
0&0
\end{bmatrix}}
\in\MATr(k[x_2]).
\lb{eq243}
\end{align}
Then there are permutations $\ig'\in\Sfrak_p,\ \ig\in\Sfrak_q$ and matrices
$G'\in\GL(p,k[x_2])$, $G\in\GL(q,k[x_2])$ with homogeneous components 
satisfying the following conditions.
\begin{parenumr}
\item\lb{c1211}
The $j$th column of $G'UG$ is a homogeneous element of 
$\dss{i=1}{p}k[x_1,x_2](-\ng_{\ig'(i)})$ of degree $\mg_{\ig(j)}$ for all $\ee{1}{j}{q}$.
\item\lb{c1212}
The $(j+p-q,j)$ component of $G'UG$ lies in
$k^\ast x_2^{\mg_{\ig(j)}-\ng_{\ig'(j+p-q)}}$ for all $\ee{1}{j}{q}$ and the
$(i,j)$ component is zero for all $i,j$ with  $i>j+p-q$.
\item\lb{c1213}
For every $\ee{1}{l}{q}$,
\begin{equation*}
\smm{j=l}{q}(\mg_{\ig(j)}-\ng_{\ig'(j+p-q)})\leq
\smm{j=l}{q}\mg_{j}-\kg(\ng,l,p,q)+(\mg_q-\mg_l),
\end{equation*}
where
\begin{equation*}
\kg(\ng,l,p,q):=\min\lt\{\ \smm{j=l}{q}\ng_{r_j}\ \lt|\ 
\begin{matrixl}
\text{$r_l\ddd r_q$ are distinct integers with}\\
\text{$\ne{j}{r_j}{p}$ for all $\ee{l}{j}{q}$}
\end{matrixl}
\rt.\rt\}.
\end{equation*}
\end{parenumr}
\end{lem}

\begin{pf}
First of all, $U\clmlsc{1\ddd q}{\ }\neq0$, since the ideal in $k[x_1,x_2]$ generated 
by the maximal minors of $U$ is of height greater than one by hypothesis. 
We prove our assertion by induction on $q$. 
When  $q=1$, there is an $i_1\ (\ne{1}{i_1}{p})$ such that 
the $(i_1,1)$ component is $c_{i_11}x_2^{\mg_1-\ng_{i_1}}$ with $c_{i_11}\in k^\ast$. 
Just moving the $i_1$th row to the last, we obtain our assertion in this case.
Suppose $q>1$ and that our assertion
is true for smaller values of $q$.
Let 
$\ubar_{i_1}$ $(\ne{q}{i_1}{p})$ be a row of $U$
different from zero.
Since the components of $\ubar_{i_1}$ are homogeneous elements
of $k[x_2]$, they are of the form $cx_2^l$ $(c\in k,\ l\geq0)$.
We can write $\ubar_{i_1}=(0\ddd 0,c_{i_1 j_1}x_2^{\mg_{j_1}-\ng_{i_1}},\ast\ddd\ast)$, 
where $c_{i_1 j_1}x_2^{\mg_{j_1}-\ng_{i_1}}$ is the 
$(i_1,j_1)$ component of $U$ with $c_{i_1j_1}\in k^\ast$.
There is therefore an invertible 
homogeous matrix $G_1\in\GL(q,k[x_2])$ representing a column operation on $U$,
which makes no change in the first $j_1$ columns, such that 
$\ubar_{i_1} G_1=(0\ddd 0,c_{i_1 j_1}x_2^{\mg{j_1}-\ng{i_1}},0\ddd0)$
and the colums of $U G_1$ are still
homogeneous elements of $\dss{i=1}{p}k[x_1,x_2](-\ng_i)$ of degrees $\seq{\mg}{1}{q}$.
Let $G_2\in\GL(q,k)$ be the matrix representing the permutation that moves 
the $j_1$th column to the last and the $j$th column to the $(j-1)$th $(\ne{j_1}{j}{q})$.
Put $G'_1:=
\begin{smallbmatrix}
G_2^{-1}G_1^{-1}&0\\
0&1_{p-q}
\end{smallbmatrix}$ if \eqref{eq242} holds or $j_1>1$ and \eqref{eq243} holds.
In the case $j_1=1$ and \eqref{eq243} holds, put
$G'_1:=
\begin{smallbmatrix}
G_2^{-1}&0\\
0&1_{p-q}
\end{smallbmatrix}$.
Let further $G'_2\in\GL(p,k)$ be the the matrix representing the permutation that moves 
the $i_1$th row to the last and the $i$th row to the $(i-1)$th $(\ne{i_1}{i}{p})$.
Put $U':=G'_2G'_1UG_1G_2$,
$(\mg'_1\ddd\mg'_{q-1},\mg'_q):=(\mg_1\ddd\mg_{j_1-1},\mg_{j_1+1}\ddd\mg_q,\mg_{j_1})$
and $(\ng'_1\ddd\ng'_{p-1},\ng'_p):=
(\ng_1\ddd\ng_{j_1-1},\ng_{j_1+1}\ddd\ng_q,\ng_{j_1},\ng_{q+1}\ddd
\ng_{i_1-1},\ng_{i_1+1}\ddd\ng_p,\ng_{i_1})$.
Then,
the $j$th column of $U'$ is a homogeneous element of 
$\dss{i=1}{p}k[x_1,x_2](-\ng'_i)$ of degree $\mg'_{j}$ for all $\ee{1}{j}{q}$,
the $(p,q)$ component lies in $k^\ast x_2^{\mg_{j_1}-\ng_{i_1}}=k^\ast x_2^{\mg'_q-\ng'_p}$, 
the $(p,j)$ components are zero for all $j<q$, 
the sequence $\mg'_1\ddd\mg'_{q-1}$ is nondecreasing, 
and $\mg'_i=\ng'_i+1$ for all $\ee{1}{i}{q}$.
Moreover,
$U'-\bclmsc{x_11_q}{0}$ or 
$U'-\begin{smallbmatrix}
0&0\\
0&x_11_{q-1}\\
0&0
\end{smallbmatrix}$ or
$U'-\begin{smallbmatrix}
x_11_{q-1}&0\\
0&0
\end{smallbmatrix}$ lies in $\MATr(k[x_2])$.
The height of the ideal generated by the maximal minors of
$U'\clmsc{p}{q}$ must also be greater than one. Put $U'':=U'\clmsc{p}{q}$.
By the induction hypothesis applied to 
$U''$,
there are permutations $\ig'_0\in\Sfrak_{p-1},\ \ig_0\in\Sfrak_{q-1}$ and matrices
$G'_0\in\GL(p-1,k[x_2])$, $G_0\in\GL(q-1,k[x_2])$ 
such that the $j$th column of $G'_0U''G_0$ is a homogeneous element of 
$\dss{i=1}{p-1}k[x_1,x_2](-\ng'_{\ig'_0(i)})$ of degree $\mg'_{\ig_0(j)}$ for all 
$\ee{1}{j}{q-1}$, 
the $(j+p-q,j)$ component of $G'_0U''G_0$ lies in
$k^\ast x_2^{\mg'_{\ig_0(j)}-\ng'_{\ig'_0(j+p-q)}}$ for all $\ee{1}{j}{q-1}$, the
$(i,j)$ component is zero for all $i,j$ with  $i>j+p-q$,
and 
\begin{equation}
\smm{j=l}{q-1}(\mg'_{\ig_0(j)}-\ng'_{\ig'_0(j+p-q)})\leq
\smm{j=l}{q-1}\mg'_{j}-\kg(\ng',l,p-1,q-1)+(\mg'_{q-1}-\mg'_l)
\lb{eq241}
\end{equation}
for every $\ee{1}{l}{q-1}$.
There are therefore permutations $\ig',\ \ig$ and matrices
$G',\ G$ such that \eqref{c1211} and \eqref{c1212} hold.
Here, we have 
$(\mg_{\ig(1)}\ddd \mg_{\ig(q)})=(\mg'_{\ig_0(1)}\ddd\mg'_{\ig_0(q-1)},\mg'_q)$,
$(\ng_{\ig'(1)}\ddd \ng_{\ig'(p)})=(\ng'_{\ig'_0(1)}\ddd\ng'_{\ig'_0(p-1)},\ng'_p)$,
$G':=\begin{smallbmatrix}
G'_0&0\\
0&1
\end{smallbmatrix}G'_2G'_1$, and
$G:=G_1G_2\begin{smallbmatrix}
G_0&0\\
0&1
\end{smallbmatrix}$.
We want to verify \eqref{c1213}. 
The case $l=q$ is easy since $\mg_{\ig(q)}=\mg_{j_1}$ and $\ng_{\ig'(p)}=\ng_{i_1}$.
Let $l$ be an integer with $\ee{1}{l}{q-1}$ and put
\begin{equation*}
\dg:=\smm{j=l}{q-1}\mg'_{j}-\kg(\ng',l,p-1,q-1)+(\mg'_{q-1}-\mg'_l)+(\mg'_q-\ng'_p).
\end{equation*}
Since $\smm{j=l}{q}(\mg_{\ig(j)}-\ng_{\ig'(j+p-q)})\leq\dg$ 
by \eqref{eq241}, it is enough to show the following inequality.

\begin{claim}
With the notaion above, 
$\dg\leq\smm{j=l}{q}\mg_{j}-\kg(\ng,l,p,q)+(\mg_q-\mg_l)$.
\end{claim}

\begin{pfofclaim}
Let $r_l\ddd r_{q-1}$ be distinct integers
satisfying $\ne{j}{r_j}{p-1}$ for all $\ee{l}{j}{q-1}$ such that 
$\kg(\ng',l,p-1,q-1)=\smm{j=l}{q-1}\ng'_{r_j}$. Then,
\begin{equation*}
\dg=\smm{j=l}{q-1}\mg'_{j}-
\smm{j=l}{q-1}\ng'_{r_j}+(\mg'_{q-1}-\mg'_l)+(\mg_{j_1}-\ng_{i_1}).
\end{equation*}
Given integers $j,j'\ (j< j')$ and a sequence
$s_1,s_2,\ldots$, denote by $T(s,j,j')$ the subsequence $s_{j+1}\ddd s_{j'}$ with
order being forgotten. Notice that 
$T(\ng',j,p-1)\bslash\{\ng'_q=\ng_{j_1}\}=T(\ng,j+1,p)\bslash\{\ng_{i_1}\}$ 
for $\ee{j_1-1}{j}{q-1}$ and that 
$T(\ng',j,p-1)=T(\ng,j,p)\bslash\{\ng_{i_1}\}$ 
for $\ee{1}{j}{j_1-1}$. Moreover
\begin{equation*}
\kg(\ng,l,p,q)=\min\lt\{\ \smm{j=l}{q}\ng_{r_j}\ \lt|\
\begin{matrixl}
\text{$\ng_{r_j}$ is a term of $T(\ng,j,p)$ for all $\ee{l}{j}{q}$}\\
\text{and $r_l\ddd r_q$ are distinct}
\end{matrixl}
\rt.\rt\}.
\end{equation*}
If $1\leq l\leq j_1-1$ and $q\notin\{\seq{r}{l}{q-1}\}$, then
$\smm{j=l}{j_1-1}\ng'_{r_j}+\ng_{i_1}+\smm{j=j_1}{q-1}\ng'_{r_j}\geq\kg(\ng,l,p,q)$
by the above observation, so that
\begin{align*}
\dg&=\smm{j=l}{j_1-1}\mg'_{j}+\mg_{j_1}+\smm{j=j_1}{q-1}\mg'_{j}-
\smm{j=l}{j_1-1}\ng'_{r_j}-\ng_{i_1}-\smm{j=j_1}{q-1}\ng'_{r_j}+(\mg'_{q-1}-\mg'_l)\\
&=\smm{j=l}{j_1-1}\mg_{j}+\mg_{j_1}+\smm{j=j_1}{q-1}\mg_{j+1}-
\smm{j=l}{j_1-1}\ng'_{r_j}-\ng_{i_1}-\smm{j=j_1}{q-1}\ng'_{r_j}+(\mg'_{q-1}-\mg'_l)\\
&\leq\smm{j=l}{q}\mg_{j}-\kg(\ng,l,p,q)+(\mg_q-\mg_l).
\end{align*}
If $1\leq l\leq j_1-1$ and $q=r_{j_2}$ for some $\ee{l}{j_2}{q-1}$, let $r'_l\ddd r'_{q-1}$ be the
sequence obtained from $\seq{r}{l}{q-1}$ by replacing $r_{j_2}$ with $p$ 
(i.e. $r'_{j_2}=p$ and $r'_j=r_j$ for $j\neq j_2$).
Since $\ng'_{r_{j_2}}=\ng_{j_1}$, $\ng'_{r'_{j_2}}=\ng'_p=\ng_{i_1}$, and
$\ng'_{r'_{j_1-1}}$ is a term of $T(\ng',j_1-1,p)\bslash\{\ng'_q=\ng_{j_1}\}=T(\ng,j_1,p)$,
we find similarly that
\begin{align*}
\dg&=\smm{j=l}{j_1-2}\mg_{j}+\mg_{j_1-1}+\mg_{j_1}+
\smm{j=j_1}{q-1}\mg_{j+1}-
\smm{j=l}{j_1-2}\ng'_{r'_j}-\ng_{j_1}-\ng'_{r'_{j_1-1}}-
\smm{j=j_1}{q-1}\ng'_{r'_j}\\
&\hspace{5em}+(\mg'_{q-1}-\mg'_l)\\
&\leq\smm{j=l}{q}\mg_{j}-\kg(\ng,l,p,q)+(\mg_q-\mg_l).
\end{align*}
If $\ee{j_1}{l}{q-1}$ and $q\notin\{\seq{r}{l}{q-1}\}$, then
\begin{align*}
\dg&=\mg_{j_1}+\smm{j=l}{q-1}\mg_{j+1}-\ng_{i_1}-
\smm{j=l}{q-1}\ng'_{r_j}+(\mg'_{q-1}-\mg'_l)\\
&\leq\smm{j=l}{q}\mg_{j}-\kg(\ng,l,p,q)+(\mg_q-\mg_l).
\end{align*}
If $\ee{j_1}{l}{q-1}$ and $q=r_{j_2}$ for some $\ee{l}{j_2}{q-1}$, let $r'_l\ddd r'_{q-1}$ be the
sequence obtained from $\seq{r}{l}{q-1}$ by replacing $r_{j_2}$ with $p$.
Since $\ng'_{r_{j_2}}=\ng_{j_1}$, $\ng'_{r'_{j_2}}=\ng'_p=\ng_{i_1}$, $\mg'_{q-1}=\mg_q$,
$\mg'_l=\mg_{l+1}$, and
$\mg_{j_1}-\ng_{j_1}=\mg_l-\ng_l=\mg_{l+1}-\ng_{l+1}=1$,
\begin{align*}
\dg&=\mg_{j_1}+\smm{j=l}{q-1}\mg_{j+1}-\ng_{j_1}-
\smm{j=l}{q-1}\ng'_{r'_j}+(\mg'_{q-1}-\mg'_l)\\
&=\mg_l+\smm{j=l}{q-1}\mg_{j+1}-\ng_l-
\smm{j=l}{q-1}\ng'_{r'_j}+(\mg_q-\mg_{l+1})\\
&=\mg_l+\smm{j=l}{q-1}\mg_{j+1}-\ng_{l+1}-
\smm{j=l}{q-1}\ng'_{r'_j}+(\mg_l-\ng_l)+(\ng_{l+1}-\mg_{l+1})+(\mg_q-\mg_l)\\
&=\mg_l+\smm{j=l}{q-1}\mg_{j+1}-\ng_{l+1}-
\smm{j=l}{q-1}\ng'_{r'_j}+(\mg_q-\mg_l)\\
&\leq\smm{j=l}{q}\mg_j-\kg(\ng,l,p,q)+(\mg_q-\mg_l).
\end{align*}
\end{pfofclaim}
Thus, we have \eqref{c1213}, too.
\end{pf}

\begin{lem}\lb{proc122}
Let $(\zg_{ij})$, $\Gg$, $U_{01}$, $U_1$ $\Utild_{01}$ and $\Utild_1$ be the matrices mentioned at 
the begininig of
this section such that $U_{01}$ and $U_1$ are obtained from $\Utild_{01}$ and $\Utild_1$ respectively
by the substitution 
$(\zg_{ij})=\Gg$, where $\Gg$ is chosen sufficiently generally.
Let $q$ be a positive integer with $q<a$. Then, there are a 
$a\tm (a+1)$ matrix $G'$ with components in $k[x_2]$ and a $G\in\GL(q,k[x_2])$ such that the matrix
$(u_{ij}):=G'\lt(\bclmsc{U_{01}}{U_1}\clmsc{\phantom{1}}{q+1\ddd a}\rt)G$ satisfies the following
conditions.
\begin{parenumr}
\item\lb{c1221}
Each $u_{j+a-q\,j}$ is a homogeneous element of $k[x_2]\bslash\{0\}+(\seqtwo)R$ 
for all $\ee{1}{j}{q}$ and $u_{ij}\in(\seqtwo)R$ for all $i,j$ with $i>j+a-q$.
\item\lb{c1222}
The degree of $\prdd{j=l}{q}u_{j+a-q\,j}$ is less than or equal to 
$q-l+1$ for all $\ee{1}{l}{q}$.
\end{parenumr}
\end{lem}

\begin{pf}
Let the notation be the same as in \Lem{\ref{proc114}}. 
One of the three cases stated there occurs.
Put
\begin{equation*}
\Ucheck:=
\lt\{
\begin{aligned}
{}&\Ucheck_1\clm{\phantom{1}}{q+1\ddd a}\ \text{for the case \eqref{c1141}},\quad\\
{}&\Gtild\bclm{\Ucheck_{01}}{\Ucheck_1}\clml{1}{q+1\ddd a}\ 
\text{for the case \eqref{c1142}}.
\end{aligned}\rt.
\end{equation*}
Then, $\Ucheck-\bclmsc{z_11_q}{0}\in\MATr(K[z_2])$.
In the case \eqref{c1143} of \Lem{\ref{proc114}}, we should be a little bit
more careful.
When $q\geq i_0$, put 
\begin{equation*}
\Ucheck':=\lt(\bclm{\Gtild\Ucheck_1}{\Ucheck_{01}}\Gtild^{-1}\rt)\clml{1}{q+1\ddd a}
=\lt(\bclm{\Gtild\Ucheck_1}{\Ucheck_{01}}\clml{1}{q+1\ddd a}\rt)
\begin{bmatrix}
\Gtild'{}^{-1}&0\\
0&1_{q-i_0}
\end{bmatrix},
\end{equation*}
and let $\Ucheck$ be the matrix obtained from $\Ucheck'$ by moving the
$i$th row to the $(i+1)$th for $\ee{1}{i}{a-1}$ and the last row to the first.
Then $\Ucheck-\begin{smallbmatrix}
0&0\\
0&z_11_{q-1}\\
0&0
\end{smallbmatrix}\in\MATr(K[z_2])$.
Suppose $q<i_0$. Since $\Gtild'$ is an invertible matrix,
$\rk{K}{\Gtild'\clmlsc{1}{q+1\ddd i_0}}\geq q-1$. There is therefore a
$\Gtild''\in\GL(i_0-1,K)$ and a $\Gtild'''\in\GL(q,K)$ such that 
$\Gtild''\lt(\Gtild'\clmlsc{1}{q+1\ddd i_0}\rt)\Gtild'''$ coincides with $\bclmsc{1_q}{0}$ or
$\begin{smallbmatrix}
0&0\\
0&1_{q-1}\\
0&0\end{smallbmatrix}$.
Put
\begin{equation*}
\Ucheck':=
\begin{bmatrix}
\Gtild''&0\\
0&1_{a+1-i_0}
\end{bmatrix}
\lt(\bclm{\Gtild\Ucheck_1}{\Ucheck_{01}}\clml{1}{q+1\ddd a}\rt)\Gtild'''.
\end{equation*}
Then, since $\Gtild=
\begin{smallbmatrix}
\Gtild'&0\\
0&1_{a-i_0}
\end{smallbmatrix}$, we see
\begin{multline*}
\Ucheck'-
\bclm{z_1\Gtild''\Gtild'\clmlsc{1}{q+1\ddd i_0}\Gtild'''}{0}\\
=\Ucheck'-
\begin{bmatrix}
\Gtild''&0\\
0&1_{a+1-i_0}
\end{bmatrix}
\lt(\bclm{\Gtild z_11_a}{0}\clml{1}{q+1\ddd a}\rt)\Gtild'''
\in\MATr(K[z_2]).
\end{multline*}
Hence $\Ucheck'-\bclmsc{z_11_q}{0}$ or
$\Ucheck'-\begin{smallbmatrix}
0&0\\
0&z_11_{q-1}\\
0&0\end{smallbmatrix}$ lies in $\MATr(K[z_2])$.
Now, in the case where \eqref{c1143} of \Lem{\ref{proc114}} holds with $q<i_0$, 
let $\Ucheck$ be the matrix obtained from $\Ucheck'$ by moving 
the $i$th row to the $(i+1)$th for $\ee{i_0}{i}{a-1}$ and the 
last row to the $i_0$th. In all the above cases, there are a 
$a\tm (a+1)$ matrix $\Gtild'_0$ with components in $K$ and a $\Gtild_0\in\GL(q,K)$ 
such that 
$\Ucheck\equiv
\Gtild'_0\lt(\bclmsc{\Utild_{01}}{\Utild_1}\clmsc{\phantom{1}}{q+1\ddd a}\rt)\Gtild_0\
(\mod\ (z_3,z_4)R_K)$.
Put $\Utild:=\Gtild'_0\lt(\bclmsc{\Utild_{01}}{\Utild_1}\clmsc{\phantom{1}}{q+1\ddd a}\rt)\Gtild_0$
and let $U_\irrd$ (resp. $U$) denote the matrix obtained from $\Utild$ (resp. $\Ucheck$)
by the substitution $(\zg_{ij})=\Gg$. Since the height of the ideal in $K[z_1,z_2]$ generated by the
maximal minors of $\Ucheck$ is greater than one by \Lem{\ref{proc114}}, 
so is the height of the ideal in $k[x_1,x_2]$ generated by the maximal minors of $U$.
Besides, $U_\irrd\equiv U\ (\mod\ (x_3,x_4)R)$. 
Since $U$ is obtained in the manner above, we find that it satisfies all the hypotheses of
\Lem{\ref{proc121}} with $p=a$,
$(\seq{\ng}{1}{p})=(\seq{n}{1}{a})$, and $(\seq{\mg}{1}{q})=(n_1+1\ddd n_q+1)$.
\par
Let $\ig'\in\Sfrak_a,\ \ig\in\Sfrak_q$, $G'\in\GL(a,k[x_2])$, and $G\in\GL(q,k[x_2])$ be 
the permutations and the matrices stated in
\Lem{\ref{proc121}} and let $u_{ij}$ denote the $(i,j)$ component of $G'U_\irrd G$. 
Since the components of $G'U G$ satisfy \eqref{c1212} of 
\Lem{\ref{proc121}} and since $G'U_\irrd G\equiv G'U G\ (\mod\ (x_3,x_4)R)$,
we see that \eqref{c1221} holds and that the degree of $u_{j+a-q\,j}$ is
$\mg_{\ig(j)}-\ng_{\ig'(j+a-q)}$. On the other hand, 
$\smm{j=l}{q}\mg_{j}-\kg(\ng,l,a,q)+(\mg_q-\mg_l)=\smm{j=l}{q}(n_j+1)-\smm{i=l}{q}n_{i+1}
+(n_q-n_l)\leq q-l+1$ for all $\ee{1}{l}{q}$ by a direct computation.
Since the degree of $\prdd{j=l}{q}u_{j+a-q\,j}$
must be $\smm{j=l}{q}(\mg_{\ig(j)}-\ng_{\ig'(j+a-q)})$,
our assertion \eqref{c1222} follows from \eqref{c1213} of \Lem{\ref{proc121}}.
\end{pf}

\section{\mathversion{bold}Matrices that represent
operations of $x_1$ and $x_2$ on the $\ktwo$-module $\Ext{3}{R}{R/I}{R}$}
\lb{sec3}

As seen in Section \ref{sec1}, the module $\Ext{3}{R}{R/I}{R}$,
which is the dual of $\Hm{1}{R/I}$ up to shift in grading,
is isomorphic to $\Coksup{R}{\tpose{U_3},\tpose{U_5},\tpose{U_4}}$.
The operation of $x_1$ (resp. $x_2$) on an element of $\Ext{3}{R}{R/I}{R}$ is 
therefore represented by the matrix 
$\tpose{\Uovercc\kernsub{3}}=x_11_b-\tpose{U_3}$ 
(resp. $\tpose{\Uovercc\kernsub{5}}=x_21_b-\tpose{U_5}$).
In this section, we describe how these matrices are influenced by 
the structure of $\Ext{3}{R}{R/I}{R}$ as an $\ktwo$-module from a computational aspect.
Our arguments below will be applied in Section \ref{sec4} with 
$(V_3,V_5)=(\tpose{U_3},\tpose{U_5})$.

\par
Assume that $\seqfour$ are elements of $\gradn{R}{1}$ such that $R=k[\seqfour]$.
Let $\seq{d}{1}{Q}$ be a nondecreasing sequence of integers, and $V_3$ and $V_5$ 
be matrices giving homogneous homomorphisms 
\begin{equation*}
V_3,\ V_5:\dss{i=1}{Q}R(d_i+1)\lra\dss{i=1}{Q}R(d_i+2)
\end{equation*}
of degree zero 
such that the components of $\Vovercc\kernsub{3}:=x_11_Q-V_3$ and
$\Vovercc\kernsub{5}:=x_21_Q-V_5$ lie in $k[x_3,x_4]$. 
Suppose there are positive integers $p,\ q$ with $\en{1}{p}{Q},\ \ee{1}{q}{Q-p}$
such that $d_p<d_{p+1}=d_{p+2}=\cdots=d_{p+q}$ and if $q<Q-p$ then
$d_{p+q}<d_{p+q+1}$.
Take an element $s\in k$ and put
$\xbar_1:=x_1+sx_2$, $\Wovercc:=\Vovercc\kernsub{3}+s\Vovercc\kernsub{5}$,
$W:=\xbar_11_Q-\Wovercc=V_3+sV_5$. Let
\begin{equation*}
W=:\begin{bmatrix}
C'&C''&C'''\\
D'&D''&D'''
\end{bmatrix},\quad 
A':=\begin{bmatrix}
C''&C'''\\
D''&D'''
\end{bmatrix},\quad
A'':=\begin{bmatrix}
C''\\
D''
\end{bmatrix},
\end{equation*}
where the number of the rows of $(C',C'',C''')$
 (resp. $(D',D'',D''')$) is $p$
(resp. $Q-p$), and the numbers of the columuns of $D',D'',D'''$ are $p,q,Q-p-q$
respectively.
Notice that for each row of $A''$ the degrees of its components are the same
and that $C''$ 
is a matrix whose components are zero or of degree zero in 
$\xbar_1,x_3,x_4$.
Besides, $C'''=0$ since the degrees its components must be negative.
Choose $G_1\in \GL(q,k)$ so that
the columns of $C''G_1$  
different from zero are linearly independent over $k$.
Put
\begin{equation*}
G_2:=
\begin{bmatrix}
1_{p}&0&0\vspace{3pt}\\
0&G_1&0\vspace{3pt}\\
0&0&1_{Q-p-q}
\end{bmatrix}\quad\andt\quad
A:=A'
\begin{bmatrix}
G_1&0\\
0&1_{Q-p-q}
\end{bmatrix}.
\end{equation*}
Observe that $\xbar_1$ appears in $A$ only in the form $c\xbar_1$ with
some $c\in k$. We may therefore write 
\begin{equation}
A=\xbar_1A_1+A_0,
\lb{eq1}
\end{equation}
where $A_1$ (resp. $A_0$) is a matrix with entries in $k$ (resp. $k[x_3,x_4]$)
and $\Dg(A_0)=\Dg(A)=\Dg(A_1)+1$. Let $S$ denote the graded module 
$\dss{i=1}{Q}R(d_i+2)$ and 
$\deg(v)$ the degree of an element $v\in S$.
We will regard the columns of $A$ as 
homogeneous elements of $S$.
Note that the degrees of the first $q$ columns of 
$A$ are the same and equal to $-1-d_{p+1}$, while the degrees of the remaining
columns are smaller than that.
Let $\seq{b}{1}{m}$ be all the columns of $A$ of degree $-1-d_{p+1}$
which do not vanish modulo $(\xbar_1,x_3,x_4)$, and denote the remaining
columns of $A$ by $\seq{a}{1}{n}$.
Actually, $\{\seq{b}{1}{m}\}$ consists of all the
columns $b$ of $A''G_1$ such that at least one of the first $p$ components of $b$ is 
an element of $k$ different from zero. 
Let $b_i=:\tpose{(b_{1i}\ddd b_{Qi})}$ and $b''_i:=\tpose{(b_{1i}\ddd b_{pi})}$.
The components of 
$b''_i$ lie in $k$ and 
the vectors $b''_1\ddd b''_m$ are linearly independent over $k$ by the choice of $G_1$.
We construct a matrix $H$ whose columns are  
$\xbar_1^{\ng_1}x_3^{\ng_3}x_4^{\ng_4}a_i$ 
$(\ng_1+\ng_3+\ng_4+\deg(a_i)=-1-d_{p+1},\ \ee{1}{i}{n},\ 
\ng_1,\ \ng_3,\ \ng_4\geq0)$ arranged in a suitable order.
Observe that $\deg(a_i)\leq -1-d_{p+1}$ for all $\ee{1}{i}{n}$.
The first $p$ components of $a_i$ $(\ee{1}{i}{n})$ are therefore zero
and so are the first $p$ rows of $H$.
Since $W-\xbar_11_Q\in\MATr(\ktwo)$, we see $A'-\bclmsc{0}{\xbar_11_{Q-p}}\in\MATr(\ktwo)$.
Hence
\begin{equation*}
G_2^{-1}A-
\begin{bmatrix}
O\\
\xbar_11_{Q-p}
\end{bmatrix}
\in \MATr(\ktwo).
\end{equation*}
This implies that
\begin{multline*}
\dss{i=1}{Q}k[\xbar_1,x_3,x_4](d_i+2)=
\Imasup{k[\xbar_1,x_3,x_4]}{G_2^{-1}A}\op\\
\lt(\lt(\dss{i=1}{p}k[\xbar_1,x_3,x_4](d_i+2)\rt)\opdot
\lt(\dss{i=p+1}{Q}\ktwo(d_i+2)\rt)\rt)
\end{multline*}
by \cite[\Lem{1.1}]{A5}, so that multiplying both sides by 
$G_2$ on the left, we get
\begin{multline*}
\dss{i=1}{Q}k[\xbar_1,x_3,x_4](d_i+2)=\Imasup{k[\xbar_1,x_3,x_4]}{A}\op\\
\lt(\lt(\dss{i=1}{p}k[\xbar_1,x_3,x_4](d_i+2)\rt)\opdot
\lt(\dss{i=p+1}{Q}\ktwo(d_i+2)\rt)\rt).
\end{multline*}
Moreover, since an element of 
$\gradn{k[\xbar_1,x_3,x_4](d_i+2)}{-1-d_{p+1}}$ is zero or lies in $k$
for all $\ee{1}{i}{p}$, we see
\begin{equation}
\begin{split}
&\gradn{\dss{i=1}{Q}k[\xbar_1,x_3,x_4](d_i+2)}{-1-d_{p+1}}\\
&\hspace{6em}=\gradn{\Imasup{k[\xbar_1,x_3,x_4]}{A}}{-1-d_{p+1}}\op
\gradn{\dss{i=1}{Q}\ktwo(d_i+2)}{-1-d_{p+1}}\\
&\hspace{6em}=\la \seq{b}{1}{m},\ H\ra \op
\gradn{\dss{i=1}{Q}\ktwo(d_i+2)}{-1-d_{p+1}},
\end{split}
\lb{eq3}
\end{equation}
where $\la \seq{b}{1}{m},\ H\ra$ denotes the vector space over $k$
spanned by $\seq{b}{1}{m}$ and the columns of $H$.

\par
We will denote a matrix $Z$ with components
in $R$ by $Z(x_1,x_2,x_3,x_4)$ when we want to pay attention to
the variables $x_1,x_2,x_3,x_4$.
Let $\xg,\ \hg$ be parameters over $R$.
For a matrix $Z=Z(x_1,x_2,x_3,x_4)$ with components in $R$,
let  
$\Ztild=\Ztild(\xg,\hg,x_1,x_2,x_3,x_4)
:=Z(x_1,x_2,x_3+\xg\xbar_1,x_4+\hg\xbar_1)$ and 
$\Zbar=\Zbar(\xg,\hg,x_1,x_2)
:=Z(x_1,x_2,\xg\xbar_1,\hg\xbar_1)$.
Observe that 
\begin{equation}
\Zbar(x_3/\xbar_1,x_4/\xbar_1,x_1,x_2)=Z(x_1,x_2,x_3,x_4).
\lb{eq20}
\end{equation}
Now we consider $\seq{b}{1}{m}$ and $H$. Notice that their components are
polynomials in $\xbar_1,x_3,x_4$. With the notation above, the components of 
$\seq{\btild}{1}{m}$ and $\Htild$, therefore, lie in 
$k[\xbar_1,x_3,x_4,\xg,\hg]=k[\xbar_1,x_3,x_4]\ot_k k[\xg,\hg]$. 
Let $T$ be the local ring $k[\xg,\hg]_{(\xg,\hg)}$ and
let $\la \seq{\btild}{1}{m},\ \Htild\ra_T$ denote 
the submodule of 
\begin{equation*}
\gradn{\dss{i=1}{Q}k[\xbar_1,x_3,x_4](d_i+2)
\ot_k T}{-1-d_{p+1}}
\end{equation*}
spanned over $T$ by $\seq{\btild}{1}{m}$ and the columns of $\Htild$,
where $\gradn{S\ot_k T}{\rg}=\gradn{S}{\rg}\ot_k T$ for $\rg\in \Zbf$.
Since $\btild_i(0,0,x_1,x_2,x_3,x_4)=b_i\ (\ee{1}{i}{m})$ and 
$\Htild(0,0,x_1,x_2,x_3,x_4)=H$, we find by
\eqref{eq3} that 
\begin{multline}
\gradn{\dss{i=1}{Q}k[\xbar_1,x_3,x_4](d_i+2)
\ot_k T}{-1-d_{p+1}}\\
=\la \seq{\btild}{1}{m},\ \Htild\ra_T \op
\gradn{\dss{i=1}{Q}\ktwo(d_i+2)\ot_k T}{-1-d_{p+1}}.
\lb{eq4}
\end{multline}
Put
\begin{equation*}
N_p:=0^p\opdot 
\lt(\dss{i=p+1}{Q}k[x_3,x_4](d_i+2)\rt)\sset\lt(\dss{i=1}{Q}k[x_3,x_4](d_i+2)\rt).
\end{equation*}
This is a module over $\ktwo$ consisting of all the elements of 
$\dss{i=1}{Q}k[x_3,x_4](d_i+2)\sset S$
such that the first $p$ components are zero. 
Let $v=v(x_3,x_4)\ (=v(x_1,x_2,x_3,x_4))$ be an element of $\grad{N_p}{-1-d_{p+1}}$.
We can write 
\begin{multline}
\vtild(\xg,\hg,x_1,x_2,x_3,x_4)\\
=
\smm{i=1}{m}\btild_i(\xg,\hg,x_1,x_2,x_3,x_4)\ghat_i+
\Htild(\xg,\hg,x_1,x_2,x_3,x_4)\cdot\tpose{\row{\fhat}{1}{l}}+w
\lb{eq5}
\end{multline}
with $\ghat_i,\ \fhat_i\in T$ and 
$w\in \gradn{\dss{i=1}{Q}\ktwo(d_i+2)\ot_k T}{-1-d_{p+1}}$ by \eqref{eq4},
where $l$ denotes the number of the columns of $H$.

\begin{lem}\lb{proc2}
Let $v$ and $w$ be as above. Suppose that 
$w\equiv 0\ (\mod\ (x_3,x_4))$. Then 
\begin{equation*}
\xbar_1^\ng v\in \Imasup{k[\xbar_1,x_3,x_4]}{\seq{a}{1}{n}}
\op N_p
\end{equation*}
for all $\ng\geq0$.
\end{lem}

\begin{pf}
Since $v\in\gradn{N_p}{-1-d_{p+1}}$ by hypotheses, 
the first $p$ components of $\vtild$
are zero. Besides, since an element of 
$\gradn{k[\xbar_1,x_3,x_4](d_i+2)\ot_k T}{-1-d_{p+1}}$ is zero or lies in
$T$ for $\ee{1}{i}{p}$, it follows from the assumption $w\equiv 0\ (\mod\ (x_3,x_4))$
that the first $p$ components of $w$ are also zero.
On the other hand, the vectors $\btild''_i$
$(\ee{1}{i}{m})$ are linearly independent over $T$ and the first $p$ components of 
the columns of
$\Htild$ are zero, since these properties
are inherited from $b''_1\ddd b''_m$ and $H$.
Hence, $\ghat_i=0$ for all $\ee{1}{i}{m}$. In other words
\begin{equation*}
\vtild(\xg,\hg,x_1,x_2,x_3,x_4)=\Htild(\xg,\hg,x_1,x_2,x_3,x_4)
\cdot\tpose{\row{\fhat}{1}{l}}+w,
\end{equation*}
so that
\begin{equation*}
\vbar(\xg,\hg,x_1,x_2)=\Hbar(\xg,\hg,x_1,x_2)
\cdot\tpose{\row{\fhat}{1}{l}}.
\end{equation*} 
Let $\ng\geq0$ be an integer. Since the denominators of $\fhat_i\ (\ee{1}{i}{l})$ lies 
in $k^\ast+(\xg,\hg)$, there is a polynomial $\psg_0\in(\xg,\hg)^{\ng+1}$
such that $\psg_i:=(1+\psg_0)\fhat_i\in k[\xg,\hg]$ for all $\ee{1}{i}{l}$. Hence
\begin{equation*}
(1+\psg_0)\xbar_1^\ng\vbar(\xg,\hg,x_1,x_2)=\xbar_1^\ng\Hbar(\xg,\hg,x_1,x_2)
\cdot\tpose{\row{\psg}{1}{l}}.
\end{equation*} 
Now substitute $x_3/\xbar_1$ and $x_4/\xbar_1$ for $\xg$
and $\hg$ respectively in this equality.
We find by \eqref{eq20} that
\begin{multline*}
(1+\psg_0(x_3/\xbar_1,x_4/\xbar_1))\xbar_1^\ng v(x_3,x_4)\\=
\xbar_1^\ng H(x_1,x_2,x_3,x_4)\cdot\tpose{(\psg_1(x_3/\xbar_1,x_4/\xbar_1),
\ldots,\psg_l(x_3/\xbar_1,x_4/\xbar_1))}.
\end{multline*}
Write
\begin{equation*}
\psg_j(x_3/\xbar_1,x_4/\xbar_1)=
\sm{\mg\geq0}\psg_{j\mg}(x_3,x_4)/\xbar_1^\mg
\end{equation*}
for $\ee{0}{j}{l}$,
where $\psg_{j\mg}(x_3,x_4)$ is a homogeneous polynomials in 
$x_3,x_4$ of degree $\mg$ for each $j,\mg$. Notice that $\psg_{0\mg}=0$ for
$\mg\leq\ng$. Moreover, we can write 
$a_i=\xbar_1a_{i1}+a_{i0}$ by \eqref{eq1}, where $a_{i1}$ (resp. $a_{i0}$)
is a column of $A_1$ (resp. $A_0$).
Compare terms with no factor $\xbar_1$ in the denominators in the above equality.
Then, since $v\in\dss{i=1}{Q}\ktwo(d_i+2)$ and each column  of $H$ is of the form
$\xbar_1^{\ng_1}x_3^{\ng_3}x_4^{\ng_4}a_i$, we find that
$\xbar_1^\ng v$ is the sum of a finite number of vectors of the forms
\begin{align*}
&\xbar_1^{\ng_1+\ng-\mg}x_3^{\ng_3}x_4^{\ng_4}
\psg_{j\mg}(x_3,x_4) a_i \\
&\hskip3em(\deg(a_i)+\ng_1+\ng_3+\ng_4=-1-d_{p+1},\ 
\ng_1,\ \ng_3,\ \ng_4\geq0,\ \ng_1+\ng\geq\mg,\\
&\hskip6em \ee{1}{i}{n},\ \ee{1}{j}{l})\quad \andr\\
&x_3^{\ng_3}x_4^{\ng_4}
\psg_{j\ng_1+\ng+1}(x_3,x_4) a_{i1} \\
&\hskip3em(\deg(a_i)+\ng_1+\ng_3+\ng_4=-1-d_{p+1},\ 
\ng_1,\ \ng_3,\ \ng_4\geq0,\\ 
&\hskip6em \ee{1}{i}{n},\ \ee{1}{j}{l}).
\end{align*}
Since the first $p$ components of $a_{i1}$
are zero for all $\ee{1}{i}{n}$, our assertion holds. 
\end{pf}

\begin{lem}\lb{proc3}
Assumption being the same as in \Lem{\ref{proc2}}, we have
\begin{equation*}
\Wovercc{}^\ng v\in N_p
\end{equation*}
for all $\ng\geq0$.
\end{lem}

\begin{pf}
For $\ng=0$, our assertion is trivial. Suppose $\ng>0$.
Using the equality $W=\xbar_1 1_Q-\Wovercc$, we see 
\begin{equation*}
\xbar_1^\ng v
=W\lt(\smm{i=1}{\ng}\xbar_1^{\ng-i}\Wovercc{}^{i-1}\rt)v
+\Wovercc{}^\ng v.
\end{equation*}
On the other hand, $\xbar_1^\ng v=\phg_1+\phg_2$ with 
$\phg_1\in\Imasup{k[\xbar_1,x_3,x_4]}{\seq{a}{1}{n}}
\sset \Imasup{k[\xbar_1,x_3,x_4]}{W}$, $\phg_2\in N_p$ by \Lem{\ref{proc2}}.
The vector $\phg_2-\Wovercc{}^\ng v$ is therefore contained in
$\Imasup{k[\xbar_1,x_3,x_4]}{W}$. Since the components of 
$\phg_2-\Wovercc{}^\ng v$ are all contained in $\ktwo$ and $W$ is a
square matrix such that $W=\xbar_1 1_Q-\Wovercc$ with
$\Wovercc\in \MATr(\ktwo)$,
it follows that $\phg_2-\Wovercc{}^\ng v=0$. Hence $\Wovercc{}^\ng v\in N_p$.
\end{pf}

\begin{lem}\lb{proc44}
Let $V_3$, $V_5$, $d_i$ $(\ee{1}{i}{Q})$, 
$p$, $q$, $N_p$, $\xg,\ \hg$,
$T$ be as above. 
Let $\seq{v}{1}{\dg}$ be homogeneous elements of 
$\dss{i=1}{Q}\ktwo(d_i+2)$ and 
$v$ an element of $\gradn{N_p}{-1-d_{p+1}}\cap\Imasup{\ktwo}{\seq{v}{1}{\dg}}$.
Suppose that one of the following two conditions holds:
\begin{parenumr}
\item\lb{cond100}
$d_{p+1}>d_p+1$,
\item\lb{cond101}
$d_{p+1}=d_p+1$ and the number of the  minimal generators of the
$\ktwo$-module 
\begin{equation*}
\dss{i=1}{Q}R(d_i+2)/\Imasup{R}{V_3,V_5,\seq{v}{1}{\dg}}
\end{equation*}
of degree $-1-d_{p+1}$ is $q':=\max\{\ i\ |\ d_{p-i+1}=\cdots=d_{p-1}=d_p\ \}$ 
and remains unchanged for any small
homogeneous transformation of variables 
$x_1,x_2,x_3,x_4$.
\end{parenumr}
\par
Then, for every $s\in k$, the element 
$w$ of $\gradn{\dss{i=1}{Q}\ktwo(d_i+2)\ot_k T}{-1-d_{p+1}}$ defined by the equality
\eqref{eq5} must be congruent to zero modulo $(\seqtwo)$. 
\end{lem}

\begin{pf}
In the case \eqref{cond100}, we get our assertion, taking
the degrees of the components of $v$ and $w$ into account.
Let us go on the case \eqref{cond101}.
Let $s\in k$, and let $\xbar_1$, $W$, $A$, $b_i\ (\ee{1}{i}{m})$, and $H$ be as above.
We want to consider $\Coksup{R\ot_k k[\xg,\hg]}{\Vtild_3,\Vtild_5,\seq{\vtild}{1}{\dg}}$,
where $\Vtild_i=\Vtild(\xg,\hg,x_1,x_2,x_3,x_4)\ (i=3,5)$ and 
$\vtild_j=\vtild_j(\xg,\hg,x_1,x_2,x_3,x_4)\ (\ee{1}{j}{\dg})$.
It follows from the equality
\begin{equation*}
\dss{i=1}{Q}R(d_i+2)=\Imasup{R}{V_3,V_5}+\lt(\dss{i=1}{Q}\ktwo(d_i+2)\rt)
\end{equation*}
that
\begin{equation*}
\dss{i=1}{Q}R(d_i+2)\ot_k T=\Imasup{R\ot_k T}{\Vtild_3,\Vtild_5}+
\lt(\dss{i=1}{Q}\ktwo(d_i+2)\ot_k T\rt).
\end{equation*}
Hence 
\begin{multline}
\lt(\dss{i=1}{Q}R(d_i+2)\ot_k k[\xg,\hg]_\psg\rt)/
\Imasup{R\ot_k k[\xg,\hg]_\psg}{\Vtild_3,\Vtild_5,\seq{\vtild}{1}{\dg}}\\
\cong \lt(\dss{i=1}{Q}\ktwo(d_i+2)\ot_k k[\xg,\hg]_\psg\rt)/(\ktwo\ot_k k[\xg,\hg]_\psg)E''
\lb{eq6}
\end{multline}
over $\ktwo\ot_k k[\xg,\hg]_\psg$ for some $\psg\in k[\xg,\hg]\bslash(\xg,\hg)$, where
\begin{equation*}
E'':=\Imasup{R\ot_k k[\xg,\hg]}{\Vtild_3,\Vtild_5,\seq{\vtild}{1}{\dg}}
\cap \lt(\dss{i=1}{Q}\ktwo(d_i+2)\ot_k k[\xg,\hg]\rt).
\end{equation*}
Let $w$ be the element of
$\dss{i=1}{Q}\ktwo(d_i+2)\ot_k T$ of 
degree $-1-d_{p+1}$ such that the equality \eqref{eq5} holds.
The vectors $\btild_i\ (\ee{1}{i}{m})$ and the columns of $\Htild$  
are contained in
$\Imasup{R\ot_k T}{\Vtild_3}+\Imasup{R\ot_k T}{\Vtild_5}$, and 
$\vtild\in\Imasup{R\ot_k T}{\seq{\vtild}{1}{\dg}}$, so that 
$w\in (\ktwo\ot_k T)E''$. Let 
$\phi_i:=\tpose{(0\ddd0,\overset{\indexonsmile{i}}1,0\ddd0)}$
$(\ee{1}{i}{Q})$ be the canonical bases of 
$\dss{i=1}{Q}\ktwo(d_i+2)\ot_k k[\xg,\hg]$. Then $\deg(\phi_i)=-1-d_{p+1}$
if and only if $\ee{p-q'+1}{i}{p}$. Since the parameters 
$\xg,\ \hg$ correspond to a small homogeneous transformation of
the variables $\seqfour$, our hypothesis implies by \eqref{eq6} that the $q'$ vectors
$\seq{\phi}{p-q'+1}{p}$ must be linearly independent over 
$k(\pfrak)$ in 
$\lt(\lt(\dss{i=1}{Q}\ktwo(d_i+2)\ot_k k[\xg,\hg]\rt)/E''\rt)\ot k(\pfrak)$
for all points $\pfrak\in\Spec{k[\xg,\hg]}$ in a neighborhood of
the origin $\xg=\hg=0$.
We find therefore that the first $p$ components
of any element of $E''$ of degree $-1-d_{p+1}$ must be zero.
Hence $w\equiv 0 \ (\mod\ (x_3,x_4))$.
\end{pf}

\section{\mathversion{bold}Structure of $\Imasup{R}{\tpose{U_4}}$}
\lb{sec8}

To investigate the structure of $\Coksup{R}{\tpose{U_3},\tpose{U_5},\tpose{U_4}}$
as a $\ktwo$-module, one needs to know the
behavior of $\tpose{\Uovercc\kernsub{3}},\ \tpose{\Uovercc\kernsub{5}}$ and their
multiples more intimately. This section is devoted to a description of
some rudimentary properties concerned with it.

\par
Let $V_3,\ V_5,\ \Vovercc\kernsub{3},\ \Vovercc\kernsub{5}$ and 
$d_i$ $(\ee{1}{i}{Q})$ be the same as in the previous section.
Here, we think of $V_3,\ V_5,\ \Vovercc\kernsub{3},\ \Vovercc\kernsub{5}$
as homogeneous linear mappings of degree one from 
$\dss{i=1}{Q}R(d_i+2)$ to itself.

\begin{lem}\lb{proc9}
Let $V_3,\ V_5,\ \Vovercc\kernsub{3},\ \Vovercc\kernsub{5}$, 
$d_i$ $(\ee{1}{i}{Q})$ be as above and let
$\Xg$ denote the module generated over $\ktwo$ by the
columns of 
$[\Vovercc\kernsub{3},\Vovercc\kernsub{5}]
:=\Vovercc\kernsub{3}\Vovercc\kernsub{5}-\Vovercc\kernsub{5}\Vovercc\kernsub{3}$.
Let further $\mg\geq0$, $\ng\geq0$ be integers
and $(\seq{p}{1}{\mg+\ng})$ be a $(\mg+\ng)$-tuple of integers such that 
$\mg=\#\{\ i\ |\ p_i=3\ \}$, $\ng=\#\{\ i\ |\ p_i=5\ \}$. 
Then for every
$v\in\dss{i=1}{Q}\ktwo(d_i+2)$, 
we have 
\begin{equation*}
\lt(\prdd{i=1}{\mg+\ng}\Vovercc\kernsub{p_i}\rt)v\in
\Vovercc\kernsub{3}^\mg\Vovercc\kernsub{5}^\ng v
+\sm{
	\begin{smallmatrix}
	m\geq0,\ n\geq0,\\
  0\leq m+n<\mg+\ng
	\end{smallmatrix}
 }
\Vovercc\kernsub{3}^{m}\Vovercc\kernsub{5}^{n}\Xg,
\end{equation*}
where $\prdd{i=1}{\mg+\ng}\Vovercc\kernsub{p_i}$ stands for
$1_Q$ in the case $\mg=\ng=0$.
\end{lem}

\begin{pf}
Put $l=\mg+\ng$. Let $\mg''$, $\ng''$ be integers with $\mg''\geq0$, $\ng''>0$ and
let $l'':=\mg''+\ng''\leq l$.
Assume that $p_i=3$ for $i\leq\mg''$, $p_i=5$ for $\ee{\mg''+1}{i}{l''}$ and that
$p_{l''+1}=3$ if $l''<l$.
If $l''=l$, then $\mg=\mg'',\ \ng=\ng''$ and 
$\lt(\prdd{i=1}{l}\Vovercc\kernsub{p_i}\rt)v=
\Vovercc\kernsub{3}^\mg\Vovercc\kernsub{5}^\ng v$.
Suppose $l''<l$ and let 
$v':=\lt(\prdd{i=l''+2}{l}\Vovercc\kernsub{p_i}\rt)v\in\dss{i=1}{Q}\ktwo(d_i+2)$. 
Then, using the equality
\begin{equation*}
\Vovercc\kernsub{3}^{\mg''}\Vovercc\kernsub{5}^{\ng''}\Vovercc\kernsub{3}
\cdot v'=
\Vovercc\kernsub{3}^{\mg''}\Vovercc\kernsub{5}^{\ng''-1}
\Vovercc\kernsub{3}\Vovercc\kernsub{5}\cdot v'
-\Vovercc\kernsub{3}^{\mg''}\Vovercc\kernsub{5}^{\ng''-1}
[\Vovercc\kernsub{3},\Vovercc\kernsub{5}]v'
\end{equation*}
repeatedly, we find 
\begin{equation*}
\lt(\prdd{i=1}{l}\Vovercc\kernsub{p_i}\rt)v=
\Vovercc\kernsub{5}^{\mg''}\Vovercc\kernsub{3}^{\ng''}\Vovercc\kernsub{5}
\cdot v'
\in
\Vovercc\kernsub{5}^{\mg''+1}\Vovercc\kernsub{3}^{\ng''}v'
+\sm{
	\begin{smallmatrix}
	m\geq0,\ n\geq0,\\
  0\leq m+n<\mg+\ng
	\end{smallmatrix}
}
\Vovercc\kernsub{5}^{m}\Vovercc\kernsub{3}^{n}\Xg.
\end{equation*}
Since $\mg''+1+\ng''>l''$, our assertion follows by ascending induction on
$l''$.
\end{pf}

\par
For each pair of integers $\mg,\ng$ with $\mg\geq0$, $\ng\geq0$,
put
\begin{equation*}
L_{\mg,\ng}:=\sm{
	\begin{smallmatrix}
	\row{p}{1}{\mg+\ng}\ \text{such that}\\
	\#\{i|p_i=3\}=\mg\ \andt\ \#\{i|p_i=5\}=\ng
	\end{smallmatrix}
}
\prdd{j=1}{\mg+\ng}\Vovercc\kernsub{p_j},
\end{equation*}
where $L_{0,0}=1_Q$.

\begin{cor}\lb{proc100}
Let $V_3,\ V_5,\ \Vovercc\kernsub{3},\ \Vovercc\kernsub{5}$, 
$d_i$ $(\ee{1}{i}{Q})$,
and
$\Xg$ be as in \Lem{\ref{proc9}}.
Then,
\begin{equation*}
\binom{\mg+\ng}{\mg}\Vovercc\kernsub{3}^\mg
\Vovercc\kernsub{5}^\ng v 
\in
L_{\mg,\ng}v
+
\sm{
	\begin{smallmatrix}
	m\geq0,\ n\geq0,\\
  0\leq m+n<\mg+\ng
	\end{smallmatrix}
}
\Vovercc\kernsub{3}^{m}\Vovercc\kernsub{5}^{n}\Xg
\end{equation*}
for every $v\in\dss{i=1}{Q}\ktwo(d_i+2)$ and every pair $\mg, \ng$ of integers
with $\mg\geq0,\ng\geq0$.
\end{cor}

\begin{pf}
With the notation of \Lem{\ref{proc9}}, 
\begin{equation*}
\Vovercc\kernsub{3}^\mg\Vovercc\kernsub{5}^\ng v
\in
\lt(\prdd{i=1}{\mg+\ng}\Vovercc\kernsub{p_i}\rt)v+\sm{
	\begin{smallmatrix}
	m\geq0,\ n\geq0,\\
  0\leq m+n<\mg+\ng
	\end{smallmatrix}
 }
\Vovercc\kernsub{3}^{m}\Vovercc\kernsub{5}^{n}\Xg.
\end{equation*}
Our assertion follows by summing up all the above relations such that
$\#\{\ i\ |\ p_i=3\ \}=\mg\ \andt\ \#\{\ i\ |\ p_i=5\ \}=\ng$.
\end{pf}

\begin{lem}\lb{proc117}
Let $V_3,\ V_5$, $\Vovercc\kernsub{3}$ and $\Vovercc\kernsub{5}$ be
as in \Lem{\ref{proc9}}, and $v$ be an element of $\dss{i=1}{Q}\ktwo(d_i+2)$.
Then,
\begin{equation*}
x_2^\mg v\in\Imasup{\kthree}{V_5}+\Vovercc\kernsub{5}^{\mg}\cdot v,
\quad
x_1^\mg v\in
\Imasup{R}{V_3}+\Vovercc\kernsub{3}^{\mg}\cdot v
\end{equation*}
for all $\mg\geq1$. If further $v\in(x_3,x_4)\dss{i=1}{Q}\ktwo(d_i+2)$, then
\begin{equation*}
x_2^\mg v\in(x_3,x_4)\Imasup{\kthree}{V_5}+\Vovercc\kernsub{5}^{\mg}\cdot v,
\quad
x_1^\mg v\in
(x_3,x_4)\Imasup{R}{V_3}+\Vovercc\kernsub{3}^{\mg}\cdot v
\end{equation*}
for all $\mg\geq1$.
\end{lem}

\begin{pf}
First we have
\begin{equation*}
x_2v=(x_21_Q)v=(V_5+\Vovercc\kernsub{5})v\in
\Imasup{\kthree}{V_5}+\Vovercc\kernsub{5}\cdot v.
\end{equation*}
Inductively,
\begin{equation*}
\begin{split}
x_2^\mg v&=x_2(x_2^{\mg-1}v)\in
x_2(\Imasup{\kthree}{V_5}+\Vovercc\kernsub{5}^{\mg-1}\cdot v)\\
&\sset \Imasup{\kthree}{V_5}+
V_5(\Vovercc\kernsub{5}^{\mg-1}\cdot v)+
\Vovercc\kernsub{5}(\Vovercc\kernsub{5}^{\mg-1}\cdot v)\\
&=\Imasup{\kthree}{V_5}+\Vovercc\kernsub{5}^{\mg}\cdot v
\end{split}
\end{equation*}
for $\mg\geq1$. Likewise,
\begin{equation*}
x_1^\mg v\in
\Imasup{k[x_1,x_3,x_4]}{V_3}+\Vovercc\kernsub{3}^{\mg}\cdot v
\sset\Imasup{R}{V_3}+\Vovercc\kernsub{3}^{\mg}\cdot v
\end{equation*}
for $\mg\geq1$. If $v\in(x_3,x_4)\dss{i=1}{Q}\ktwo(d_i+2)$, then
$v=x_3v'+x_4v''$ with $v',v''\in\dss{i=1}{Q}\ktwo(d_i+2)$. 
Hence
\begin{equation*}
x_2^\mg v=x_3x_2^\mg v'+x_4x_2^\mg v''\in
(x_3,x_4)\Imasup{\kthree}{V_5}+\Vovercc\kernsub{5}^{\mg}\cdot v.
\end{equation*}
The same holds for $x_1^\mg v$.
\end{pf}

\begin{lem}\lb{proc103}
Let $V_3,\ V_5$, $\Vovercc\kernsub{3}$, $\Vovercc\kernsub{5}$ and $\Xg$ be
as in \Lem{\ref{proc9}}.
Let $\seq{v}{1}{\dg}$ be elements of $(x_3,x_4)\dss{i=1}{Q}\ktwo(d_i+2)$ and $P$
the module generated over $\ktwo$ by them. Then
\begin{equation*}
\begin{split}
 (x_3,x_4)\Imasup{R}{V_3,V_5}+&\Imasup{R}{\seq{v}{1}{\dg}}\\
=(x_3,x_4)\Imasup{R}{V_3}&\op(x_3,x_4)\Imasup{\kthree}{V_5}\\
&\op
\sm{\mg\geq0,\ \ng\geq0}
\lt(x_3\Vovercc\kernsub{3}^\mg\Vovercc\kernsub{5}^\ng \Xg+
x_4\Vovercc\kernsub{3}^\mg\Vovercc\kernsub{5}^\ng \Xg+
\Vovercc\kernsub{3}^\mg\Vovercc\kernsub{5}^\ng P\rt).
\end{split}
\end{equation*}
\end{lem}

\begin{pf}
Let $E$ (resp. $M$) denote the right (resp. left) hand side of the above equality.
We want to show first that $E$ is an $R$-module.
With the use of \Lem{\ref{proc117}}, we see
\begin{equation}
\begin{split}
\Imasup{\kthree}{[\Vovercc\kernsub{3},\Vovercc\kernsub{5}]}&=
\sm{\mg\geq0}x_2^\mg
\Imasup{\ktwo}{[\Vovercc\kernsub{3},\Vovercc\kernsub{5}]}\\
&\sset\Imasup{\kthree}{V_5}+
\sm{\mg\geq0}\Vovercc\kernsub{5}^\mg
\Imasup{\ktwo}{[\Vovercc\kernsub{3},\Vovercc\kernsub{5}]}\\
&=\Imasup{\kthree}{V_5}+
\sm{\mg\geq0}\Vovercc\kernsub{5}^\mg \Xg.
\end{split}
\lb{eq50}
\end{equation}
Let $(\seq{p}{1}{\mg+\ng})$ be a $(\mg+\ng)$-tuple of integers such that  
$\mg=\#\{\ i\ |\ p_i=3\ \}$, $\ng=\#\{\ i\ |\ p_i=5\ \}$, and $l$ be an integer with $\ee{1}{l}{\dg}$.
Since $v_l=x_3v'_l+x_4v''_l$ with $v'_l,v''_l\in\dss{i=1}{Q}\ktwo(d_i+2)$, 
\begin{equation}
\begin{split}
\lt(\prdd{i=1}{\mg+\ng}\Vovercc\kernsub{p_i}\rt)v_l
&=
x_3\lt(\prdd{i=1}{\mg+\ng}\Vovercc\kernsub{p_i}\rt)v'_l
+x_4\lt(\prdd{i=1}{\mg+\ng}\Vovercc\kernsub{p_i}\rt)v''_l\\
&\in
x_3\Vovercc\kernsub{3}^\mg\Vovercc\kernsub{5}^\ng v'_l
+x_3\sm{m\geq0,\ n\geq0 }
\Vovercc\kernsub{3}^{m}\Vovercc\kernsub{5}^{n}\Xg\\
&\hskip5em
+x_4\Vovercc\kernsub{3}^\mg\Vovercc\kernsub{5}^\ng v''_l
+x_4\sm{m\geq0,\ n\geq0 }
\Vovercc\kernsub{3}^{m}\Vovercc\kernsub{5}^{n}\Xg\\
&=
\Vovercc\kernsub{3}^\mg\Vovercc\kernsub{5}^\ng v_l
+\sm{m\geq0,\ n\geq0}
\lt(x_3\Vovercc\kernsub{3}^{m}\Vovercc\kernsub{5}^{n}\Xg
+x_4\Vovercc\kernsub{3}^{m}\Vovercc\kernsub{5}^{n}\Xg\rt)
\end{split}
\lb{eq210}
\end{equation}
by \Lem{\ref{proc9}}.
For $w\in\dss{i=1}{Q}\kthree(d_i+2)$,
\begin{equation}
\begin{split}
x_1V_5\cdot w&=(x_11_Q)V_5\cdot w=
(V_3+\Vovercc\kernsub{3})V_5\cdot w\\
&=V_3(V_5\cdot w)+V_5\Vovercc\kernsub{3}\cdot w
+[\Vovercc\kernsub{3},V_5]w\\
&=V_3(V_5\cdot w)+V_5\Vovercc\kernsub{3}\cdot w
-[\Vovercc\kernsub{3},\Vovercc\kernsub{5}]w\\
&\in\Imasup{R}{V_3}+\Imasup{\kthree}{V_5}+
\Imasup{\kthree}{[\Vovercc\kernsub{3},\Vovercc\kernsub{5}]}.
\end{split}
\lb{eq51}
\end{equation}
Hence
\begin{equation}
x_1(x_3,x_4)\Imasup{\kthree}{V_5}\sset E
\lb{eq52}
\end{equation}
by \eqref{eq50} and \eqref{eq51}.
Using \Lem{\ref{proc117}} with $\mg=1$, together with
\Lem{\ref{proc9}} and \eqref{eq210}, we see
\begin{equation}
x_i\lt(x_3\Vovercc\kernsub{3}^\mg\Vovercc\kernsub{5}^\ng \Xg+
x_4\Vovercc\kernsub{3}^\mg\Vovercc\kernsub{5}^\ng \Xg+
\Vovercc\kernsub{3}^\mg\Vovercc\kernsub{5}^\ng P\rt)
\sset E
\lb{eq53}
\end{equation}
for $i=1,\ 2$. Thus $E$ is a module over $R$ by \eqref{eq52} and \eqref{eq53}.
Moreover the columns of the matrices $x_jV_i$ ($i\in\{3,5\},\ j\in\{3,4\}$) and the
elements $\seq{v}{1}{\dg}$ are contained in $E$. Hence $M\sset E$. 
To complete the proof of our assertion, it remains to prove that 
\begin{equation*}
x_3\Vovercc\kernsub{3}^\mg\Vovercc\kernsub{5}^\ng \Xg+
x_4\Vovercc\kernsub{3}^\mg\Vovercc\kernsub{5}^\ng \Xg+
\Vovercc\kernsub{3}^\mg\Vovercc\kernsub{5}^\ng P
\sset M.
\end{equation*}
Let $v=x_3 v'+x_4 v''$ be an arbitrary element of $(x_3,x_4)\dss{i=1}{Q}R(d_i+2)$. 
Then, by the same computation as in the proof of \Lem{\ref{proc117}} with $\mg=1$, we see 
\begin{equation*}
\Vovercc\kernsub{3}\cdot v \in (x_3,x_4)\Imasup{R}{V_3}+x_1 v \quad\andt\quad
\Vovercc\kernsub{5}\cdot v \in (x_3,x_4)\Imasup{R}{V_5}+x_2 v.
\end{equation*}
This means that 
$\Vovercc\kernsub{3}^\mg\Vovercc\kernsub{5}^\ng v\in M$
whenever $v\in M$. Since the columns of 
$[\Vovercc\kernsub{3},\Vovercc\kernsub{5}]=[V_3,V_5]$ are elments of 
$\Imasup{R}{V_3}+\Imasup{R}{V_5}$, the sets $x_j\Xg$ $(j=3,4)$ are
contained in $M$. Besides, it is clear that $P\sset M$.
Hence, $x_3\Vovercc\kernsub{3}^\mg\Vovercc\kernsub{5}^\ng \Xg+
x_4\Vovercc\kernsub{3}^\mg\Vovercc\kernsub{5}^\ng \Xg+
\Vovercc\kernsub{3}^\mg\Vovercc\kernsub{5}^\ng P
\sset M$.
\end{pf}

\begin{cor}\lb{proc111}
With the notation of \Lem{\ref{proc103}},
\begin{equation*}
\begin{split}
&\lt((x_3,x_4)\Imasup{R}{V_3,V_5}+\Imasup{R}{\seq{v}{1}{\dg}}\rt)\cap
\dss{i=1}{Q}\ktwo(d_i+2)\\
&\hskip 4em  =\sm{\mg\geq0,\ \ng\geq0}
\lt(x_3\Vovercc\kernsub{3}^\mg\Vovercc\kernsub{5}^\ng \Xg+
x_4\Vovercc\kernsub{3}^\mg\Vovercc\kernsub{5}^\ng \Xg+
\Vovercc\kernsub{3}^\mg\Vovercc\kernsub{5}^\ng P\rt).
\end{split}
\end{equation*}
\end{cor}

\begin{pf}
Let $v$ be an element of the left hand side of the above equality.
Then, by \Lem{\ref{proc103}}, $v=\phg_1+\phg_2+\phg_3$ with
$\phg_1\in(x_3,x_4)\Imasup{R}{V_3}$,
$\phg_2\in(x_3,x_4)\Imasup{\kthree}{V_5}$ and 
$\phg_3\in\sm{\mg\geq0,\ \ng\geq0}
\lt(x_3\Vovercc\kernsub{3}^\mg\Vovercc\kernsub{5}^\ng \Xg+
x_4\Vovercc\kernsub{3}^\mg\Vovercc\kernsub{5}^\ng \Xg+
\Vovercc\kernsub{3}^\mg\Vovercc\kernsub{5}^\ng P\rt)$.
Since $v-\phg_3=\phg_1+\phg_2$ lies in 
$\lt(\Imasup{R}{V_3}\op\Imasup{\kthree}{V_5}\rt)\cap
\dss{i=1}{Q}\ktwo(d_i+2)=0$, one sees that
$v$ and $\phg_3$ must coincide. Hence the left hand side is contained
in the right hand side. The inclusion of the other direction is obvious.
\end{pf}

\begin{prop}\lb{proc105}
Assume that $\ch{k}=0$.
Let $\Vovercc\kernsub{3}$, $\Vovercc\kernsub{5}$, $\Xg$, $\seq{v}{1}{\dg}$
and $P$ be as in \Lem{\ref{proc103}}. 
Suppose that $\seq{v}{1}{\dg}$ are homogeneous and that there is a homogeneous
submodule $P_0$ of $(\seqtwo)\dss{i=1}{Q}\ktwo(d_i+2)$ over $\ktwo$ satisfying
\begin{equation}
\lt\{
\begin{split}
&\Xg\sset \sm{\ng\geq0}\Vovercc\kernsub{5}^{\ng}P_0+P,\\
&x_iP_0\sset P_0+(x_3,x_4)\Imasup{R}{V_3,V_5}+\Imasup{R}{\seq{v}{1}{\dg}}
\end{split}\rt.
\lb{eq212}
\end{equation}
for $i=1,2$. Then, 
\begin{equation*}
\sm{\mg\geq0,\ \ng\geq0}\Vovercc\kernsub{3}^\mg\Vovercc\kernsub{5}^\ng (P_0+P)
=P_0+\sm{m\geq0,\ n\geq0}L_{m,n}P.
\lb{eq216}
\end{equation*}
\end{prop}

\begin{pf}
Since $v_j=x_3v'_j+x_4v''_j$ with $v'_j,v''_j\in\dss{i=1}{Q}\ktwo(d_i+2)$, 
it follows from \Cor{\ref{proc100}} that
\begin{equation}
\binom{\mg+\ng}{\mg}\Vovercc\kernsub{3}^\mg\Vovercc\kernsub{5}^\ng v_j
\in
L_{\mg,\ng}v_j+
\sm{
	\begin{smallmatrix}
	m\geq0,\ n\geq0,\\
  0\leq m+n<\mg+\ng
	\end{smallmatrix}
}
\lt(x_3\Vovercc\kernsub{3}^{m}\Vovercc\kernsub{5}^{n}\Xg+
x_4\Vovercc\kernsub{3}^{m}\Vovercc\kernsub{5}^{n}\Xg\rt)
\lb{eq232}
\end{equation}
for all $\ee{1}{j}{\dg}$ and all pairs of integers $\mg,\ng$ with $\mg\geq0$, $\ng\geq0$,
in the same manner as in \eqref{eq210}. 
Let $v$ be a homogeneous element of $P_0$.
We have
\begin{equation*}
x_1^\mg x_2^\ng v\in
P_0+(x_3,x_4)\Imasup{R}{V_3,V_5}+\Imasup{R}{\seq{v}{1}{\dg}}
\end{equation*}
by \eqref{eq212}. 
There are homogeneous elements 
$v'\in P_0$, $\phg_1\in(x_3,x_4)\Imasup{R}{V_3}$,
$\phg_2\in(x_3,x_4)\Imasup{\kthree}{V_5}$ and 
$\phg_3\in\sm{m\geq0,\ n\geq0}
\lt(x_3\Vovercc\kernsub{3}^m\Vovercc\kernsub{5}^n \Xg+
x_4\Vovercc\kernsub{3}^m\Vovercc\kernsub{5}^n \Xg+
\Vovercc\kernsub{3}^m\Vovercc\kernsub{5}^n P\rt)$ 
\linebreak
such that
\begin{equation*}
x_1^\mg x_2^\ng v=v'+\phg_1+\phg_2+\phg_3
\end{equation*}
by \Lem{\ref{proc103}}.
On the other hand,
\begin{equation*}
\begin{split}
x_1^\mg x_2^\ng v=x_2^\ng x_1^\mg  v&\in 
\Vovercc\kernsub{5}^\ng \Vovercc\kernsub{3}^\mg v +
\Imasup{R}{V_3}+\Imasup{\kthree}{V_5}\\
&\sset\Vovercc\kernsub{3}^\mg \Vovercc\kernsub{5}^\ng v+
\sm{m\geq0,\ n\geq0}
\lt(x_3\Vovercc\kernsub{3}^m\Vovercc\kernsub{5}^n \Xg+
x_4\Vovercc\kernsub{3}^m\Vovercc\kernsub{5}^n \Xg\rt)\\
&\hspace{5em}+\Imasup{R}{V_3}+\Imasup{\kthree}{V_5}
\end{split}
\end{equation*}
by \Lem{\ref{proc117}} and \Lem{\ref{proc9}}, where we have applied the same
computation as in \eqref{eq210} to $v$. Let $\phg_0\in\Imasup{R}{V_3}+\Imasup{\kthree}{V_5}$
and $\phg'_3\in\sm{m\geq0,\ n\geq0}
\lt(x_3\Vovercc\kernsub{3}^m\Vovercc\kernsub{5}^n \Xg+
x_4\Vovercc\kernsub{3}^m\Vovercc\kernsub{5}^n \Xg\rt)$
be the elements such that 
$x_1^\mg x_2^\ng v=\Vovercc\kernsub{3}^\mg \Vovercc\kernsub{5}^\ng v+\phg'_3+\phg_0$.
Then
\begin{equation*}
\Vovercc\kernsub{3}^\mg \Vovercc\kernsub{5}^\ng v-v'-\phg_3+\phg'_3
\in\lt(\Imasup{R}{V_3}\op\Imasup{\kthree}{V_5}\rt)\cap\dss{i=1}{Q}\ktwo(d_i+2).
\lb{eq214}
\end{equation*}
Since the right hand side is the zero module,
$\Vovercc\kernsub{3}^\mg \Vovercc\kernsub{5}^\ng v=v'+\phg_3-\phg'_3$.
On the other hand, 
\begin{equation*}
\begin{split}
\phg_3-\phg'_3&\in
\sm{m\geq0,\ n\geq0}
\lt(x_3\Vovercc\kernsub{3}^m\Vovercc\kernsub{5}^n \Xg+
x_4\Vovercc\kernsub{3}^m\Vovercc\kernsub{5}^n \Xg+
\Vovercc\kernsub{3}^m\Vovercc\kernsub{5}^n P\rt)\\
&\sset
\sm{m\geq0,\ n\geq0}
\lt(x_3\Vovercc\kernsub{3}^m\Vovercc\kernsub{5}^n \Xg+
x_4\Vovercc\kernsub{3}^m\Vovercc\kernsub{5}^n \Xg\rt)+
\sm{m\geq0,\ n\geq0}L_{m,n}P
\end{split}
\end{equation*}
by \eqref{eq232}. Hence
\begin{equation}
 \Vovercc\kernsub{3}^\mg \Vovercc\kernsub{5}^\ng v\in
 P_0+\sm{m\geq0,\ n\geq0}L_{m,n}P+
 \sm{m\geq0,\ n\geq0}
\lt(x_3\Vovercc\kernsub{3}^m\Vovercc\kernsub{5}^n \Xg+
x_4\Vovercc\kernsub{3}^m\Vovercc\kernsub{5}^n \Xg\rt).
\lb{eq233}
\end{equation}
Let $M_0:= P_0+\sm{m\geq0,\ n\geq0}L_{m,n}P$ and 
$M_1:=\sm{\mg\geq0,\ \ng\geq0}\Vovercc\kernsub{3}^\mg\Vovercc\kernsub{5}^\ng(P_0+P)$. 
By Hamilton-Cayley's theorem, 
they are \fg\ homogeneous $\ktwo$-modules and $M_0\sset M_1$ by \Lem{\ref{proc9}} 
and \eqref{eq212}.
Since $M_1\sset M_0+(\seqtwo)M_1$ by  \eqref{eq212}, \eqref{eq232} and \eqref{eq233},
we find by Nakayama's lemma that $M_0=M_1$.
\end{pf}

\begin{cor}\lb{proc110}
In the same situation as in \Prop{\ref{proc105}}, let
\begin{equation*}
v\in\lt((x_3,x_4)\Imasup{R}{V_3,V_5}+\Imasup{R}{\seq{v}{1}{\dg}}\rt)\cap
\dss{i=1}{Q}\ktwo(d_i+2).
\end{equation*}
Then
\begin{equation*}
\Vovercc\kernsub{5}^\rg v\in
P_0+\sm{m\geq0,\ n\geq0}
L_{m,n}P
\end{equation*}
for all $\rg\geq0$.
\end{cor}

\begin{pf}
Since $\Vovercc\kernsub{5}^\rg v\in(\seqtwo)\Imasup{\kthree}{V_5}+x_2^\rg v$ by \Lem{\ref{proc117}},
$\Vovercc\kernsub{5}^\rg v$ is also an element of
$\lt((x_3,x_4)\Imasup{R}{V_3,V_5}+\Imasup{R}{\seq{v}{1}{\dg}}\rt)\cap
\dss{i=1}{Q}\ktwo(d_i+2)$. 
Our assertion follows from this immediately by \Cor{\ref{proc111}}, \eqref{eq212} and 
\Prop{\ref{proc105}}.
\end{pf}

\section{Inequalities coming from the irreducibility of a polynomial 
of minimal degree}
\lb{sec4}

Let $X$ be a curve in $\Pthree$, that is, a locally \CM\ equidimensional 
closed subscheme of dimension one of a three dimensional projective space $\Pthree$
over an infinite field $k$.
Let further $I$ be the saturated homogeneous ideal of $X$
in a polynomial ring $R$ over $k$ in four indeterminates, 
$(a;n_1,n_2\ddd n_a;n_{a+1}\ddd n_{a+b})$ the basic sequence of $X$
(i.e. the basic sequence of $I$), and $\seqfour$ sufficiently general linear forms of $R$.
We have $R=k[\seqfour]$. Throughout our argument in this section, we will
make use of the notation as in Section \ref{sec1} unless otherwise
specified.

\begin{lem}\lb{proc107}
Assume that $a\geq 3$ and that $X$ is contained in an irreducible 
surface of degree $a$.
Let $\og$ and $\seq{t}{1}{\og}$ be the integers defined in Sectinon \ref{sec2}
and suppose $\og>2$.
Let $m$ be an integer with $\en{1}{m}{\og}$, $t:=\smm{l=1}{m}t_l$, and 
$a'$ an integer with $\ee{t+2}{a'}{a}$.
Suppose $n_i=n_t+i-t$ for all $\ee{t}{i}{a'}$.
Then there is a \Wei{} basis 
$\{e^1_1,e^2_1\ddd e^2_a,e^3_1\ddd e^3_b\}$ of $I$ with respect to $\seqfour$ such 
that the submatrices $U_{01}$ and $U_1$ of the standard relation matrix 
$\lg_2$
among $e^1_1,e^2_1\ddd e^2_a,e^3_1\ddd e^3_b$ 
defined by \eqref{eq203} and \eqref{eq59} satisfies 
\begin{equation}
\begin{bmatrix}
U_{01}\\
U_1
\end{bmatrix}
=
\begin{bmatrix}
\ast&\ast&\hdotsfor{3}&\ast&\dots&\ast\\
D'&\ast&\hdotsfor{3}&\ast&\dots&\ast\\
0&c_{t}&\ast&\hdotsfor{2}&\ast&\dots&\ast\\
0&0&c_{t+1}&\ast&\dots&\ast&\dots&\ast\\
\vdots&\ddots&\ddots&\ddots&\ast&\ast&\dots&\ast\\
0&\dots&0&0&c_{a'-1}&\ast&\dots&\ast\\
0&\dots&0&0&\ast&\ast&\dots&\ast\\
\vdots&&\vdots&\vdots&\vdots&\vdots&&\vdots\\
0&\dots&0&0&\ast&\ast&\dots&\ast
\end{bmatrix}, 
\lb{eq219}
\end{equation}
where $c_i$ is the $(i+1,i)$ component of $U_1$ for each $\ee{t}{i}{a'-1}$, 
$\seq{c}{t}{a'-1}\in k^\ast$, and $D'$ is a $t\tm(t-1)$ matrix.
\end{lem}

\begin{pf}
There is an irreducible homogeneous polynomial of degree $a$
in $I$ since $X$ is contained in an irreducible surface of degree $a$
by hypothesis. 
Let $\zg_{ij}\ (\ee{1}{i}{4},\ \ee{1}{j}{4})$, $\seqfourz$, 
$K$, $R_K$ and $I_K$ be as in the beginning of Section \ref{sec2}.
It is enough to show the corresponding assertion over $K$ for $I_K$.
Let further $\{\etild^1_1,\etild^2_1\ddd \etild^2_a,\etild^3_1\ddd \etild^3_b\}$
be a \Wei\ basis of $I_K$ with respect to $\seqfourz$ and let $\lgtild_2$, $\Utild_{01}$,
$\Utild_1$, $\og$, $\seq{t}{1}{\og}$, $\seq{\Ctild}{1}{\og-1}$, $\seq{\Dtild}{1}{\og}$, 
and so on also be as in Section \ref{sec2}.
By \Lem{\ref{proc26}}, $\Ctild_l\neq0$ for all $\en{1}{l}{\og}$. 
Besides, $\og\geq3$ and $t_l=1$ for $\nn{m}{l}{m+a'-t}$ by hypothesis.
Thererfore $\rk{k}{\Ctild_l}=1$ for $\ee{m}{l}{m+a'-t-1}$, 
$\Ctild_l$ $(\nn{m}{l}{m+a'-t-1})$
are elements of $K^\ast$, and $\Ctild_m$ and $\Ctild_{m+a'-t-1}$ are 
a nonzero row and column
vectors respectively. Put $\ctild_{l+t-m}:=\Ctild_l$ for $\nn{m}{l}{m+a'-t-1}$.
Let $\ctild_i\ (\ee{a'-1}{i}{a'''-1})$ be the components of $\Ctild_{m+a'-t-1}$,
namely,
$\tpose{(\seq{\ctild}{a'-1}{a'''-1})}:=\Ctild_{m+a'-t-1}$. If $\ctild_{a'-1}=0$,
then there is an $i_0$ with $\ne{a'-1}{i_0-1}{a'''-1}$ such that $\ctild_{i_0-1}\neq0$.
Exchanging $\etild^2_{a'}$ and $\etild^2_{i_0}$, we obtain another
\Wei\ basis $\{\etild^1_1,\etild^2_1\ddd \etild^2_{a'-1},\etild'{}^2_{a'},\etild^2_{a'+1}\ddd
\etild^2_{i_0-1},\etild'{}^2_{i_0},\etild^2_{i_0+1}\ddd
\etild^2_a,\etild^3_1\ddd \etild^3_b\}$ $(\etild'{}^2_{a'}:=\etild^2_{i_0},\ \etild'{}^2_{i_0}:=\etild^2_{a'})$
by \eqref{c1085} of \Lem{\ref{proc108}}, 
with respect to which the first component of $\Ctild_{m+a'-t-1}$
is different from zero.
We may therefore assume from the first that $\ctild_{a'-1}\neq0$.
Let $\Gtild_m\in\GL(t_m,K)$ 
be as in \Lem{\ref{proc115}} and let $(0\ddd0,\ctild_{t}):=\Ctild_m\Gtild_m$ 
(i.e. $\Ctild''_m=\ctild_{t}$). Put $t':=\smm{l=1}{m-1}t_l$ and 
\begin{equation*}
(\etild'{}^2_{t'+1}\ddd \etild'{}^2_{t}):=
(\etild^2_{t'+1}\ddd \etild^2_{t})\Gtild_m.
\end{equation*}
By \eqref{c1083} of \Lem{\ref{proc108}} the 
$\{\etild^1_1,\etild^2_{1}\ddd \etild^2_{t'},\etild'{}^2_{t'+1}\ddd 
\etild'{}^2_{t},\etild^2_{t+1}\ddd 
\etild^2_a,\etild^3_1\ddd \etild^3_b\}$
is a \Wei\ basis of $I_K$. 
Since $(\etild^1_1,\etild^2_1,\ldots,\etild^2_a,\etild^3_1\ddd \etild^3_b)\lgtild_2=0$,
we have
\begin{align*}
&(\etild^1_1,\etild^2_{1}\ddd \etild^2_{t'},\etild'{}^2_{t'+1}
\ddd \etild'{}^2_{t},\etild^2_{t+1}\ddd 
\etild^2_a,\etild^3_1\ddd \etild^3_b)
\begin{bmatrix}
1&0\\
0&\Gtild^{-1}
\end{bmatrix}
\lgtild_2\Gtild=0 \quad\witht\\
&\Gtild:=
\begin{bmatrix}
1_{t'}&0&0\\
0&\Gtild_m&0\\
0&0&1_{a+b-t}
\end{bmatrix}.
\end{align*}
Let $\Utild'_i\ (\ee{1}{i}{5})$, $\Utild'_{01},\ \Utild'_{02}$,
$\Utild'_{21}$ and $\lgtild'_2$ 
be the matrices defined by \eqref{eq203} and \eqref{eq59}
with respect to this new \Wei\ basis.
Then, the above equality implies that 
\begin{equation*}
\lgtild'_2=\begin{bmatrix}
1&0\\
0&\Gtild^{-1}
\end{bmatrix}
\lgtild_2\Gtild.
\end{equation*} 
Hence, the matrix $\begin{smallbmatrix}
\Utild'_{01}\\
\Utild'_1
\end{smallbmatrix}$
coincides with the right hand side of \eqref{eq219} with tilde attached over $c_j$ and $D'$ 
by \Lem{\ref{proc115}}. Thus we can choose a desired \Wei\ basis by substituting
$(\zg_{ij})$ with a sufficiently general $\Gg\in\GL(4,k)$.
\end{pf}

\begin{lem}\lb{proc109}
Let the assumption and the notation be as in \Lem{\ref{proc107}}.
Then we can choose
$\{e^1_1,e^2_1\ddd e^2_a,e^3_1\ddd e^3_b\}$ 
so that $\lg_2$ satisfies 
\begin{equation}
U_{21}\clmlt{\phantom{1}}{1\ddd t-1,a',a'+1\ddd a}=0
\lb{eq218}
\end{equation}
in addition to \eqref{eq219}.
\end{lem}

\begin{pf}
By \Lem{\ref{proc107}}, there is a a \Wei\ basis $\{e^1_1,e^2_1\ddd e^2_a,e^3_1\ddd e^3_b\}$ 
of $I$ with respect to $\seqfour$ satisfying \eqref{eq219}.
The condition \eqref{eq218} holds if and only if 
\begin{equation}
x_1e^2_{l'}\in (\seqthree)\kthree e^1_1+\kthree e^2_1+\cdots+\kthree e^2_a
\lb{eq224}
\end{equation}
for all $\nn{t-1}{l'}{a'}$ by the standard method of computing $\lg_2$.
We show by an inductive argument that the above \Wei\ basis can be modified 
so that it also satisfies \eqref{eq224}.
Let $a''$ be an integer with $t-1<a''<a'$. Suppose that \eqref{eq219} is 
satisfied and that \eqref{eq224} holds for all $\nn{t-1}{l'}{a''}$.
We have 
\begin{equation}
x_1e^2_{a''}=g^1_1 e^1_1+g^2_1 e^2_1+\cdots+g^2_ae^2_a+
g^3_1 e^3_1+\cdots+g^3_b e^3_b
\lb{eq225}
\end{equation}
with $g^1_1\in (\seqthree)\kthree,\ g^2_l\in \kthree,\ g^3_l\in \ktwo$.
Let $a'''$ denote $\max\{\ i\ |\ n_i=n_{a'},\ \ee{1}{i}{a}\ \}$. First of all,
$g^2_{a''+1}=-c_{a''}$. Besides, $g^2_l=0$ for $\ne{a''+1}{l}{a}$ if $a''<a'-1$, and
$g^2_l=0$ for $\ne{a'''}{l}{a}$ if $a''=a'-1$, since its degree must be negative.
Put $h:=-(g^3_1 e^3_1+\cdots+g^3_b e^3_b)/c_{a''}$ and 
$e'{}^2_{a''+1}:=e^2_{a''+1}+h$. Then the polynomials 
$e^1_1,e^2_1\ddd e^2_{a''}, e'{}^2_{a''+1},e^2_{a''+2}\ddd e^2_a,e^3_1\ddd e^3_b$ 
form a \Wei\ basis of $I$ by \eqref{c1081} of \Lem{\ref{proc108}}.
Since 
\begin{equation*}
x_1e^2_{a''}=g^1_1 e^1_1+g^2_1 e^2_1+\cdots+g^2_{a''+1}e'{}^2_{a''+1}+
\cdots+g^2_ae^2_a
\end{equation*}
by \eqref{eq225}, this new \Wei\ basis satisfies the condition corresponding to
\eqref{eq224} for $l'=a''$. Since $\deg(e^2_l)=n_1+l-t$ for $\ee{t}{l}{a'}$,
we can reduce \eqref{eq224} to
\begin{equation*}
x_1e^2_{l'}\in (\seqthree)\kthree e^1_1+\kthree e^2_1+\cdots+\kthree e^2_{l'+1}
\end{equation*}
for $\nn{t-1}{l'}{a''}$.
Moreover
\begin{equation*}
\begin{split}
x_1e^2_{l'}\in (\seqthree)&\kthree e^1_1+\kthree e^2_1+\cdots+\kthree e^2_{t}\\
&+\ktwo e^3_1+\cdots+\ktwo e^3_b
\end{split}
\end{equation*}
for $\ee{1}{l'}{t-1}$ by \eqref{eq219}.
No change therefore  occurs in the first $a''-1$ columns of
$\lg_2$, even though we replace $e{}^2_{a''+1}$ with $e'{}^2_{a''+1}$.
Hence the new \Wei\ basis satisfies the condition corresponding to
\eqref{eq224} for all $\ne{t-1}{l'}{a''}$ with the constants $c_j\ (\ee{t}{j}{a''})$ 
unchanged in \eqref{eq219}, since $g^2_{a''+1}=-c_{a''}$.
If $a''=a'-1$, then we are done.
Otherwise, we want to verify also that the constants $c_{a''+1}\ddd c_{a'-1}$ 
do not change at all in the exprssion \eqref{eq219}.
We have
\begin{equation*}
\begin{split}
&x_1e^2_{a''+1}=g'{}^1_1 e^1_1+g'{}^2_1 e^2_1+\cdots+g'{}^2_ae^2_a+
g'{}^3_1 e^3_1+\cdots+g'{}^3_b e^3_b,\\
&x_1e^2_{l'}=g''{}^1_1 e^1_1+g''{}^2_1 e^2_1+\cdots+g''{}^2_ae^2_a+
g''{}^3_1 e^3_1+\cdots+g''{}^3_b e^3_b\quad(\nn{a''+1}{l'}{a'})
\end{split}
\end{equation*}
with $g'{}^1_1,g''{}^1_1\in (\seqthree)\kthree,\ g'{}^2_l,g''{}^2_l\in \kthree,\ 
g'{}^3_l,g''{}^3_l\in \ktwo$.
Then,
\begin{equation*}
\begin{split}
&x_1e'{}^2_{a''+1}=g'{}^1_1 e^1_1+g'{}^2_1 e^2_1+\cdots+
g'{}^2_{a''+1}e'{}^2_{a''+1}+\cdots+g'{}^2_ae^2_a\\
&\hskip8em+g'{}^3_1 e^3_1+\cdots+g'{}^3_b e^3_b+x_1h-g'{}^2_{a''+1}h,\\
&x_1e^2_{l'}=g''{}^1_1 e^1_1+g''{}^2_1 e^2_1+\cdots+
g''{}^2_{a''+1}e'{}^2_{a''+1}+\cdots+g''{}^2_ae^2_a\\
&\hskip8em+g''{}^3_1 e^3_1+\cdots+g''{}^3_b e^3_b-g''{}^2_{a''+1}h.
\end{split}
\end{equation*}
Applying \eqref{eq201} to 
$\{e^1_1,e^2_1\ddd e^2_{a''}, e'{}^2_{a''+1},e^2_{a''+2}\ddd e^2_a,e^3_1\ddd e^3_b\}$, 
we have
\begin{equation*}
x_1h=f'{}^1_1e^1_1+f'{}^2_1 e^2_1+\cdots+f'{}^2_{a''+1}e'{}^2_{a''+1}+
\cdots+f'{}^2_ae^2_a+f'{}^3_1 e^3_1+\cdots+f'{}^3_b e^3_b
\end{equation*}
with $f'{}^1_1,f'{}^2_l\in (\seqtwo)\kthree,\ f'{}^3_l\in \ktwo$.
Likewise, applying \eqref{eq201} repeatedly, we see
\begin{equation*}
x_2^th=f'{}^2_1 e^2_1+\cdots+f'{}^2_{a''+1}e'{}^2_{a''+1}+
\cdots+f'{}^2_ae^2_a+f'{}^3_1 e^3_1+\cdots+f'{}^3_b e^3_b
\end{equation*}
with $f'{}^2_l\in(\seqtwo)\kthree,\ f'{}^3_l\in\ktwo$ for $t>0$.
In the expression of $x_1e'{}^2_{a''+1}$ as a standard
linear combination of 
$e^1_1,e^2_1\ddd e^2_{a''}, e'{}^2_{a''+1},e^2_{a''+2}\ddd e^2_a,e^3_1\ddd e^3_b$, therefore,
the coefficients of $e^1_1,e^2_1\ddd e^2_{a''}, e'{}^2_{a''+1},e^2_{a''+2}\ddd e^2_a$
are congruent to $g'{}^1_1,g'{}^2_1\ddd g'{}^2_a$ modulo $(\seqtwo)\kthree$.
The same holds for the coefficients in the expressions of $x_1e^2_{l'}$ 
$(\nn{a''+1}{l'}{a'})$.
Hence the constants $c_{a''+1}\ddd c_{a'-1}$ do not change at all. 
The zeros appearing in the right hand side of \eqref{eq219} for reasons of degree
also do not change at all.
Thus, our new \Wei\ basis satisfies  \eqref{eq219} and the condition 
corresponding to \eqref{eq224}
for all $\ne{t-1}{l'}{a''}$, with no change in $D'$ and $\seq{c}{t}{a'-1}$. 
We reach our assertion inductively starting with $a''=t$.
\end{pf}

\par
In exactly the same manner as in the arguments of \Lems{\ref{proc107} and \ref{proc109}},
one obtains the lemma below.

\begin{lem}\lb{proc7}
Assume that $a\geq2$ and that $X$ is contained in an irreducible surface of degree $a$.
Let $\og$ and $\seq{t}{1}{\og}$ be the integers defined in Sectinon \ref{sec2}
and suppose $\og>1$.
Put $t:=\smm{l=1}{m}t_l$ for an integer $m$ with $\en{1}{m}{\og}$.
Suppose $t+1\leq a$ and $n_{t+1}=n_{t}+1$.
Then there is a \Wei{} basis 
$\{e^1_1,e^2_1\ddd e^2_a,e^3_1\ddd e^3_b\}$ of $I$ with respect to $\seqfour$ such 
that the submatrices $U_{01}$ and $U_1$ of the standard relation matrix 
$\lg_2$ among $e^1_1,e^2_1\ddd e^2_a,e^3_1\ddd e^3_b$ 
defined by \eqref{eq203} and \eqref{eq59} satisfies the following conditions:
\begin{equation}
\begin{bmatrix}
U_{01}\\
U_1
\end{bmatrix}
=
\begin{bmatrix}
\ast&\ast&\hdotsfor{2}&\ast\\
D'&\ast&\hdotsfor{2}&\ast\\
0&{\begin{smallbmatrix}1_{t-t''}\\0\end{smallbmatrix}}&\ast&\hdots&\ast&\\
0&0&\ast&\hdots&\ast\\
\vdots&\vdots&\vdots&&\vdots\\
0&0&\ast&\hdots&\ast\\
\end{bmatrix},\quad
U_{21}\clmlt{\phantom{1}}{1\ddd t'',t+1,t+2\ddd a}=0,
\lb{eq64}
\end{equation}
where $t''$ is the number of the columns of $D'$ and $t''<t$.
\end{lem}

\begin{pf}
Let $\zg_{ij}$, $z_i$, 
$K$, $R_K$, $I_K$, $\{\etild^1_1,\etild^2_1\ddd \etild^2_a,\etild^3_1\ddd \etild^3_b\}$,
$\lgtild_2$, $\Utild_i\ (\ee{1}{i}{5})$, $\Utild_{01},\ \Utild_{02}$,
$\Utild_{21}$, $\og$, $\seq{t}{1}{\og}$, $\seq{\Ctild}{1}{\og-1}$, $\seq{\Dtild}{1}{\og}$ 
be as in the proof of \Lem{\ref{proc107}}.
There is an irreducible homogeneous polynomial of degree $a$
in $I$ since $X$ is contained in an irreducible surface of degree $a$.
Let $\Gtild_m,\ \Dtild'_m,\ \Ctild''_m$ be as in \Lem{\ref{proc115}} and let
$t''_m$ denote the number of the columns of $\Ctild''_m$ which is equal to its rank.
There is a $\Gtild'_{m+1}\in\GL(t_{m+1},k)$ such that
$\Gtild'_{m+1}\Ctild''_m=
\begin{smallbmatrix}
1_{t''_m}\\
0
\end{smallbmatrix}$.
Besides, $\Gtild'_{m+1}\Dtild_{m+1}\Gtild'_{m+1}\kern-3ex{}^{-1}-z_11_{t_{m+1}}
\in\MATr(\Kthreez)$. Put $t':=\smm{l=1}{m-1}t_l$ and 
\begin{equation*}
\begin{split}
&(\etild'{}^2_{t'+1}\ddd \etild'{}^2_{t}):=
(\etild^2_{t'+1}\ddd \etild^2_{t})\Gtild_m,\\
&(\etild'{}^2_{t+1}\ddd \etild'{}^2_{t+t_{m+1}}):=
(\etild^2_{t+1}\ddd \etild^2_{t+t_{m+1}}){\Gtild'_{m+1}\kern-3ex{}^{-1}}\ .
\end{split}
\end{equation*}
By \eqref{c1083} of \Lem{\ref{proc108}} the
$\{\etild^1_1,\etild^2_{1}\ddd \etild^2_{t'},\etild'{}^2_{t'+1}\ddd \etild'{}^2_{t+t_{m+1}},
\etild^2_{t+t_{m+1}+1}\ddd 
\etild^2_a,\etild^3_1\ddd \etild^3_b\}$
is a \Wei\ basis of $I_K$. The remaining part of the proof of our assertion is almost
the same as in the proofs of \Lems{\ref{proc107} and \ref{proc109}}. 
In the course of that arugument one sees $t''_m=t-t''$.
The details are left to the readers.
\end{pf}

\begin{rem}\lb{proc124}
We can carry out the reasoning in the proof of \Lem{\ref{proc109}},
starting with the \Wei\ basis of $I_K$ obtained in the final stage of 
the proof of \Lem{\ref{proc107}}.
We may therefor assume that the \Wei\ basis of $I$ stated in \Lem{\ref{proc109}}
and its standard relation matrices are obtained from those of $I_K$ with
the same properties by the substitution $(\zg_{ij})=\Gg$
for sufficiently general $\Gg\in\GL(4,k)$. The same holds also for the
\Wei\ basis stated in \Lem{\ref{proc7}}.
\end{rem}

\begin{lem}\lb{proc101}
Let $\ag$ and $\bg$ be positive integers with $\ag\leq\bg<a$, $j'$ a positive integer
with $\bg+j'\leq a$,
$c_i$ the $(i+j',i)$ component of $U_1$ for each $\ee{\ag}{i}{\bg}$,
and $v_i$ the $i$-th column of $\tpose{U_4}$ for each $\ee{1}{i}{a}$.
Suppose that 
\begin{equation*}
U_1=
\begin{bmatrix}
&&&&&&&\vspace{-3ex}\\
&\hdotsfor{8}&\vspace{2ex}\\
0&\dots&0&c_\ag&\ast&\hdotsfor{2}&\ast&\dots&\ast\\
0&\dots&0&0&c_{\ag+1}&\ast&\dots&\ast&\dots&\ast\\
\vdots&&&\ddots&\ddots&\ddots&\ast&\ast&\dots&\ast\\
0&\hdotsfor{3}&0&0&c_\bg&\ast&\dots&\ast\vspace{1ex}\\
&\hdotsfor{8}&\vspace{1ex}
\end{bmatrix}
\end{equation*}
and that $c_i\in k^\ast$ for all $\ee{\ag}{i}{\bg}$.
Then 
\begin{equation*}
v_i\in (x_3,x_4)\Imasup{R}{\tpose{U_3},\tpose{U_5}}+
\Imasup{R}{\seq{v}{\bg+1}{a}}
\end{equation*}
for all $\ee{\ag}{i}{a}$.
\end{lem}

\begin{pf}
For each $\ee{1}{j}{a}$, let $h_j$ (resp. $h'_j$) denote the $j$-th row of $U_4$
(resp.  the $j$-th row of $U_2$). Likewise let $u_{ij}$ denote 
the $(i,j)$ component of $U_1$ for each pair $i,j$ with $\ee{1}{i}{a},\ \ee{1}{j}{a}$.
By the equality $\lg_2\lg_3=0$, we have
\begin{equation*}
(U_1,U_2,U_4)
\begin{bmatrix}
-U_4\\-U_5\\U_3
\end{bmatrix}
=0.
\end{equation*}
Let $i$ be an integer with $\ee{\ag}{i}{\bg}$.
Then, from the $(i+j')$-th row of the above equality, we get
\begin{equation*}
-c_ih_i-\smm{j=i+1}{a}u_{i+j'j}h_j-h'_{i+j'}U_5+h_{i+j'}U_3=0.
\end{equation*}
Transposing both sides of this equality, we obtain
\begin{equation*}
-c_iv_i-\smm{j=i+1}{a}u_{i+j'j}v_j-
\tpose{U_5}\tpose{h'_{i+j'}}+\tpose{U_3}\tpose{h_{i+j'}}=0,
\end{equation*}
as $v_j=\tpose{h_j}$ for all $j$. By \eqref{eq32}
the components of $h_j$ and $h'_j$ ($\ee{1}{j}{a}$) are
contained in $(x_3, x_4)\ktwo$ and $(x_3, x_4)\kthree$ respectively.
Hence 
\begin{equation*}
v_i\in (x_3,x_4)\Imasup{R}{\tpose{U_3},\tpose{U_5}}+
\Imasup{R}{\seq{v}{i+1}{a}}.
\end{equation*}
Our assertion follows from this by descending induction on $i$.
\end{pf}

\begin{lem}\lb{proc116}
Let $M,\ E$ be \fg\ graded modules over $R$ with $M\sset E$ and $\seq{v}{1}{q}$
homogeneos elements of $E$. Let further $\Phg=(f_{ij})$ be a $q\tm q$ matirx with 
components in $R$
such that $f_{ii}$
is a homogeneous element of $k[x_2]\bslash\{0\}+(\seqtwo)R$
for all $\ee{1}{i}{q}$, $f_{ij}\in(\seqtwo)R$ for all $i,j$ with $i<j$, and
the degree of $\prdd{i=l}{q}f_{ii}$ is less than or equal to 
$q-l+1$ for all $\ee{1}{l}{q}$.
Denote by $M_0$ the $\ktwo$-submodule of $E$ generated by
$x_2^{\rg}v_i\ (\ee{1}{i}{q},\ \en{0}{\rg}{q-i+1})$ over $\ktwo$.
Suppose that every component of $(\seq{v}{1}{q})\Phg$ lies in $M$ and that
$x_1v_i\in\Imasup{\kthree}{\seq{v}{1}{q}}+M$ for all $\ee{1}{i}{q}$.
Then $\Imasup{R}{\seq{v}{1}{q}}\sset M_0+M$. 
In particular, $x_iM_0\in M_0+M$ for $i=1,2$.
\end{lem}

\begin{pf}
Let $M_1:=\Imasup{\kthree}{\seq{v}{1}{q}}$. Then $M_0\sset M_1$. Since 
$Rv_j\sset M_1+M$ by the hypothesis that $x_1v_i$ is an element of $M_1+M$ 
for all $\ee{1}{i}{q}$,
it is enough to show that $M_1\sset M_0+M$. For each $\ee{1}{i}{q}$, let 
$\begin{smallbmatrix}
\Phg''_i\\
\Phg'_i
\end{smallbmatrix}$ 
be the last $q-i+1$ columns of $\Phg$, where $\Phg'_i$ is a square matrix.
By the assumption on $\Phg$, the components of $\Phg''_i$ are elemsnts of 
$(x_3,x_4)R$ and $\det(\Phg'_i)\in k^\ast x_2^{l_i}+(x_3,x_4)R$ with $\ee{0}{l_i}{q-i+1}$.
Besides, every component of $(\seq{v}{1}{q})
\begin{smallbmatrix}
\Phg''_i\\
\Phg'_i
\end{smallbmatrix}$ 
lies in $M$. Multiplying by $\det(\Phg'_i)\Phg'_i{}^{-1}$ on the right, we obtain
$\det(\Phg'_i)v_i\in \smm{j=1}{i-1}(x_3,x_4)Rv_j +M$. Hence
$x_2^{q-i+1}v_i\in\smm{j=1}{i}(x_3,x_4)Rv_j +M\sset
(x_3,x_4)M_1+M$. Since 
$x_2^{\rg}v_i\in M_0$ for $\en{0}{\rg}{q-i+1}$, this implies that 
$x_2^\rg v_i\in M_0+(x_3,x_4)M_1+M$ for all $\rg\geq0$.
This being true for all $\ee{1}{i}{q}$, we find that 
$M_1+M\sset (M_0+M)+(x_3,x_4)(M_1+M)$. Since $\gradn{M_1+M}{l}=0$
for all $l\lless 0$, we get $M_1+M=M_0+M$ by Nakayama's lemma. 
This proves our assertion.
\end{pf}

\begin{thm}\lb{proc118}
Assume $\ch{k}=0$. 
Suppose that $a\geq 2$, $b\geq1$ and that $X$ is contained in an irreducible 
surface of degree $a$.
Let $\og$ and $\seq{t}{1}{\og}$ be the integers defined in Sectinon \ref{sec2}
and suppose $\og>1$.
Let $m$ be an integer with $\en{1}{m}{\og}$, $t:=\smm{l=1}{m}t_l$, and $a'$ be an integer with 
$\ee{t+1}{a'}{a}$.
If $n_i=n_t+i-t$ for all $\ee{t}{i}{a'}$  and $n_{a+1}<n_{a'}$,
then
$t(t-1)/2 > p$, where $p:=\max\{\ i\ |\ n_{a+i}<n_{a'},\ \ee{1}{i}{b}\ \}$.
\end{thm}

\begin{pf}
To prove our assertion, we use \Wei\ basis of $I$ as stated in 
\Lems{\ref{proc107} and \ref{proc109}} or in \Lem{\ref{proc7}},
according as $a'\geq t+2$ or $a'=t+1$.
\case{1}
We first consider the case $a'\geq t+2$.
Let $\{e^1_1,e^2_1\ddd e^2_a,e^3_1\ddd e^3_b\}$
and $\lg_2$ be a \Wei{} basis of $I$ with respect to $\seqfour$
and its standard relation matrix respectively satisfying the conditions described in 
\Lems{\ref{proc107} and \ref{proc109}}.
Let further
\begin{equation*}
\lg_3:=
\begin{bmatrix}
-U_4\\-U_5\\U_3
\end{bmatrix}
\end{equation*}
and $v_j$ be the $j$-th column of $\tpose{U_4}$ for each $\ee{1}{j}{a}$.
Recall that $v_j$ is an element of $(x_3,x_4)\dss{i=1}{b}\ktwo(n_{a+i}+2)$ for all
$\ee{1}{j}{a}$. Let $M$ denote the module
$(\seqtwo)\Imasup{R}{\tpose{U_3},\tpose{U_5}}+\Imasup{R}{\seq{v}{a'}{a}}$.
Since $\lg_2\lg_3=0$ and the condition \eqref{eq219} holds, 
it follows from \Lem{\ref{proc101}} that 
\begin{equation}
v_j\in M
\lb{eq226}
\end{equation}
for all $\ee{t}{j}{a}$. This means
$\Imasup{R}{\seq{v}{t}{a}}\sset M$.
Moreover,
$x_1v_j\ (\ee{1}{j}{a})$ are contained in 
$(\seqtwo)\Imasup{R}{\tpose{U_3},\tpose{U_5}}+\Imasup{\kthree}{\seq{v}{1}{a}}$
by \Lem{\ref{proc200}}. Hence
\begin{equation}
x_1v_j\in\Imasup{\kthree}{\seq{v}{1}{t-1}}+M
\lb{eq240}
\end{equation}
for all $\ee{1}{j}{a}$. 
Suppose that $t\geq2$ for a while.
By \Lem{\ref{proc122}} and \Rem{\ref{proc124}}, there are a 
$a\tm (a+1)$ matrix $G'$ with components in $k[x_2]$ and a $G\in\GL(t-1,k[x_2])$ such that
$(u_{ij}):=G'\lt(\bclmsc{U_{01}}{U_1}\clmsc{\phantom{1}}{t\ddd a}\rt)G$
have the properties \eqref{c1221} and \eqref{c1222} stated there with $q=t-1$.
Let
$v'_1\ddd v'_{t-1}$ be the elements of $(x_3,x_4)\dss{i=1}{b}\ktwo(n_{a+i}+2)$
defined by
\begin{equation*}
(v'_1\ddd v'_{t-1}):=(\seq{v}{1}{t-1})\tpose{G^{-1}}.
\end{equation*}
First we have $\Imasup{\kthree}{\seq{v}{1}{t-1}}=\Imasup{\kthree}{v'_1\ddd v'_{t-1}}$.
The $x_1v'_j$ therefore lies in $\Imasup{\kthree}{v'_1\ddd v'_{t-1}}+M$
for all $\ee{1}{j}{t-1}$ by \eqref{eq240}.
We have
\begin{equation*}
\lt(\begin{bmatrix}
G'&0\\
0&1_{b}
\end{bmatrix}
\lg_2
\begin{bmatrix}
G&0\\
0&1_{a+2b-t+1}
\end{bmatrix}\rt)\lt(
\begin{bmatrix}
G^{-1}&0\\
0&1_{a+2b-t+1}
\end{bmatrix}\lg_3\rt)=0.
\end{equation*}
It follows therefore from the first $a$ rows that 
\begin{multline*}
-\lt(G'\lt(\bclm{U_{01}}{U_1}\clmsc{\phantom{1}}{t\ddd a}\rt)G\rt)
\tpose{(v'_1\ddd v'_{t-1})}-
\lt(G'\lt(\bclm{U_{01}}{U_1}\clmsc{\phantom{1}}{1\ddd t-1}\rt)\rt)
\tpose{(v_t\ddd v_{a})}\\
-G'\bclm{U_{02}}{U_2}U_5+G'\bclm{0}{U_4}U_3=0.
\end{multline*}
Since $(u_{ij})=G'\lt(\bclmsc{U_{01}}{U_1}\clmsc{\phantom{1}}{t\ddd a}\rt)G$ and
$\Imasup{R}{\seq{v}{t}{a}}\sset M$, 
we find by transposing the above equality that every component of
$(v'_1\ddd v'_{t_1-1})\Phg$
lies in $M$, where $\Phg$ is the transpose of the square matrix consisting of the
last $t-1$ rows of $(u_{ij})$, namely,
\begin{equation*}
\Phg=
\begin{bmatrix}
u_{a-t+2\,1}&\hdots&u_{a1}\\
\vdots&&\vdots\\
u_{a-t+2\,t-1}&\hdots&u_{a\,t-1}
\end{bmatrix}.
\end{equation*}
Notice that $\Phg$ satisfies the hypotheses stated in \Lem{\ref{proc116}}
by \Lem{\ref{proc122}}.
Let $M_0$ be the module generated by
$x_2^{\rg}v'_i\ (\ee{1}{i}{t-1},\ \en{0}{\rg}{t-i})$ over $\ktwo$.
Then, with the use of \Lem{\ref{proc116}},
we see 
$\Imasup{R}{\seq{v}{1}{t-1}}=\Imasup{R}{v'_1\ddd v'_{t-1}}\sset M_0+M$
and $x_iM_0\sset M_0+M$ for $i=1,2$. 
For each pair $i,\rg$ with $\ee{1}{i}{t-1},\ \en{0}{\rg}{t-i}$, there is a 
$\pg_{i,\rg}\in(\seqtwo)\dss{i=1}{b}\ktwo(n_{a+i}+2)$
such that
\begin{equation}
x_2^{\rg}v'_i\in(\seqtwo)\Imasup{\kthree}{\tpose{U_5}}+\pg_{i,\rg}
\lb{eq244}
\end{equation}
by \Lem{\ref{proc117}} with $V_5=\tpose{U_5}$, since $x_2^{\rg}v'_i\in\Imasup{\kthree}{\seq{v}{1}{t-1}}$.
Let $P_0$ (resp. $P_1)$ denote the $\ktwo$-module 
generated by $\pg_{i,\rg}$ $(\ee{1}{i}{t-1},\ \en{0}{\rg}{t-i})$
(resp. $v_1\ddd v_{t-1}$). 
Notice first that $P_0+M=M_0+M$ by \eqref{eq244}. 
Hence
\begin{equation*}
x_iP_0\sset P_0+M
\end{equation*}
for $i=1,2$ by what we have just seen. 
Since each $v_l$ $(\ee{1}{l}{t-1})$ is a liniear combination of $v'_1\ddd v'_{t-1}$ over 
$k[x_2]$ by its definition, it follows from \eqref{eq244} with $\rg=0$
that there are constants $s_{i\ng}\in k$ such that 
\begin{equation*}
\begin{split}
v_l&\equiv\sm{i\geq1,\ \ng\geq0}s_{i\ng}x_2^\ng\pg_{i,0}\
(\mod\ \Imasup{\kthree}{\tpose{U_5}})\\
&\equiv\sm{i\geq1,\ \ng\geq0}s_{i\ng}(\tpose{\Uovercc\kernsub{5}})^\ng\pg_{i,0}\
(\mod\ \Imasup{\kthree}{\tpose{U_5}})\\
\end{split}
\end{equation*}
again by \Lem{\ref{proc117}}. Consequently,
$v_l=\sm{i\geq1,\ \ng\geq0}s_{i\ng}(\tpose{\Uovercc\kernsub{5}})^\ng\pg_{i,0}
\in\sm{\ng\geq0}(\tpose{\Uovercc\kernsub{5}})^\ng P_0$
for all $\ee{1}{l}{t-1}$. Hence $P_1\sset\sm{\ng\geq0}(\tpose{\Uovercc\kernsub{5}})^\ng P_0$.
On the other hand, transposing the relation 
\begin{equation*}
-U_{21}U_{4}-\Uovercc\kernsub{3}\Uovercc\kernsub{5}+
\Uovercc\kernsub{5}\Uovercc\kernsub{3}=
(U_{21},U_3,U_5)
\begin{bmatrix}
-U_4\\
-U_5\\
U_3
\end{bmatrix}
=0
\end{equation*}
which follows from $\lg_2\lg_3=0$, we see that the columns of 
$[\tpose{\Uovercc\kernsub{3}},\tpose{\Uovercc\kernsub{5}}]=
\tpose{\Uovercc\kernsub{3}}\tpose{\Uovercc\kernsub{5}}
-\tpose{\Uovercc\kernsub{5}}\tpose{\Uovercc\kernsub{3}}$
are contained in the module generated by 
$\seq{v}{1}{t-1},\seq{v}{a'}{a}$ over $\ktwo$, since 
$U_{21}\clmsc{\phantom{1}}{1,2,\ddd t-1,a',a'+1\ddd a}=0$.
Let $\Xg$ denote the module over $\ktwo$ generated by the columns of 
$[\tpose{\Uovercc\kernsub{3}},\tpose{\Uovercc\kernsub{5}}]$ and 
$P$ the $\ktwo$-module generated by $\seq{v}{a'}{a}$. 
Then, as seen above,
\begin{equation*}
\Xg\sset P_1+P\sset \sm{\ng\geq0}(\tpose{\Uovercc\kernsub{5}})^\ng P_0+P.
\end{equation*}
So far we have assumed that $t\geq2$. In the case where $t=1$, let
$P_0$ be the zero module. Then, it can be verified directly that
$\Xg\sset P=\sm{\ng\geq0}(\tpose{\Uovercc\kernsub{5}})^\ng P_0+P$ and
that $x_i P_0\sset P_0+M$ for $i=1,2$.
Now, we can apply the results of Section \ref{sec3} to the present situation
with $\Vovercc\kernsub{3}=\tpose{\Uovercc\kernsub{3}}$,
$\Vovercc\kernsub{5}=\tpose{\Uovercc\kernsub{5}}$,
$V_3=\tpose{U_3}$, $V_5=\tpose{U_5}$, $Q=b$ and $d_i=n_{a+i}\ (\ee{1}{i}{b})$.
Since $v_j$ lies in $(\seqtwo)\dss{i=1}{b}\ktwo(n_{a+i}+2)$ and satisfies \eqref{eq226}
for $\ee{t}{j}{a}$, we find by \Cor{\ref{proc110}} that 
\begin{equation}
(\tpose{\Uovercc\kernsub{5}})^\rg v_j\in
P_0+\sm{m\geq0,\ n\geq0}
L_{m,n}P
\lb{eq227}
\end{equation}
for all $\rg\geq0$ and $\ee{t}{j}{a}$ with
\begin{equation*}
L_{m,n}:=\sm{
	\begin{smallmatrix}
	\row{p}{1}{m+n}\ \text{such that}\\
	\#\{i|p_i=3\}=m\ \andt\ \#\{i|p_i=5\}=n
	\end{smallmatrix}
}
\lt(\prdd{j=1}{m+n}\tpose{\Uovercc}_{p_j}\rt).
\end{equation*}
 This holds also for $\ee{1}{j}{t-1}$ by
\Prop{\ref{proc105}}, since 
$v_j\in P_1\sset\sm{\ng\geq0}(\tpose{\Uovercc\kernsub{5}})^\ng P_0$
for this range.

\par
Suppose $n_{a+1}<n_{a'}$ and put 
$p:=\max\{\ i\ |\  n_{a+i}<n_{a'},\ \ee{1}{i}{b}\ \}$. Then $p\geq1$. 
Recall that
\begin{equation*}
\Ext{3}{R}{R/I}{R}\cong\dss{i=1}{b}R(n_{a+i}+2)/
\Imasup{R}{\tpose{U_3},\tpose{U_5},\seq{v}{1}{a}}
\end{equation*}
by \Lem{\ref{proc200}}.
Since $X$ is locally \CM\ by assumption, the length of this module must be finite.
Let $N_p:=0^p\opdot\dss{i=p+1}{b}\ktwo(n_{a+i}+2)$.
For each $\ee{1}{j}{a}$, we have 
$\deg(v_j)=-1-n_j$ by our standard method of constructing $\lg_2$. 
In other words
\begin{equation*}
\Dg(v_j)=\tpose{(n_{a+1}+1-n_j\ddd n_{a+p}+1-n_j\ddd n_{a+b}+1-n_j)}.
\end{equation*}
Since $n_{a+i}+1\leq n_{a'}$ for $\ee{1}{i}{p}$ and since each component of
$v_j$ lies in $(\seqtwo)\ktwo$, we see $v_j\in N_p$ for all $\ee{a'}{j}{a}$.
By the maximality of $p$, three cases are possible: (i) $p=b$,
(ii) $p<b$ and $n_{a+p}+1=n_{a+p+1}=n_{a'}$,
(iii) $p<b$ and $n_{a+p}+1<n_{a+p+1}$.
In the case (i), $P$ must be zero since $N_p=N_b=0$. Hence 
\begin{equation*}
\Ext{3}{R}{R/I}{R}\cong\lt(\dss{i=1}{b}\ktwo(n_{a+i}+2)\rt)/P_0
\end{equation*}
by \Lem{\ref{proc200}} and \eqref{eq227}. 
Since $P_0$ is generated by $t(t-1)/2$ elements contained in 
$(x_3,x_4)\dss{i=1}{b}\ktwo(n_{a+i}+2)$ over $\ktwo$,
the length of $\Ext{3}{R}{R/I}{R}$ cannot be finite if $t(t-1)/2\leq b=p$.
Hence $t(t-1)/2>p$.
We next consider the case (ii).
Let $q':=\max\{\ i\ |\ n_{a+p-i+1}=\cdots= n_{a+p-1}=n_{a+p}\ \}$.
Since the basic sequence $(a;n_1\ddd n_a;
n_{a+1}\ddd n_{a+b})$ is
an invariant of $X$, which does not vary under any sufficiently general
homogeneous transformation of variables, the presentation of the module $\Ext{3}{R}{R/I}{R}$
described in \Lem{\ref{proc200}}
holds for all sufficiently general $\seqfour$.
The number of the  minimal generators of the
$\ktwo$-module $\Ext{3}{R}{R/I}{R}$
of degree $-1-n_{a+p+1}$ is therefore $q'$, and remains unchanged for any small
homogeneous transformation of variables 
$x_1,x_2,x_3,x_4$. 
For each $j$ with $\ee{a'}{j}{a}$, let $v''_j:=x_3^{n_j-n_{a+p+1}}v_j$.
Then $\deg(v''_j)=-1-n_{a+p+1}$ and 
$v''_j\in \gradn{N_p}{-1-n_{a+p+1}}\cap \Imasup{\ktwo}{\seq{v}{1}{a}}$, 
since $\deg(v_j)=-1-n_j$ and $v_j\in N_p$ for all $\ee{a'}{j}{a}$.
We find thererfore by \Lem{\ref{proc3}} and the case \eqref{cond101} 
of \Lem{\ref{proc44}} that
\begin{equation*}
(\tpose{\Uovercc\kernsub{3}}+s\tpose{\Uovercc\kernsub{5}})^\ng v''_j\in N_p
\end{equation*}
for all $\ng\geq0$ and all $s\in k$.
Taking the coefficient of $s^n$ in the expansion of the above relation
with $\ng=m+n$, we obtain $L_{m,n}v''_j\in N_p$.
Hence $L_{m,n}v_j\in N_p$ for all $\ee{a'}{j}{a}$, $m\geq0$, and $n\geq0$.
Thus
\begin{equation}
\sm{m\geq0,\ n\geq0}
L_{m,n}P\sset N_p.
\lb{eq228}
\end{equation}
Let $\prr:\dss{i=1}{b}\ktwo(n_{a+i}+2)\lra \dss{i=1}{p}\ktwo(n_{a+i}+2)$
be the natural projection to the first $p$ components. Then
\begin{equation*}
\begin{split}
&\len{\ktwo}{\Ext{3}{R}{R/I}{R}}\\
&\hskip2em= l_{\ktwo}\lt(
\lt(\dss{i=1}{b}\ktwo(n_{a+i}+2)\rt)/\lt(
\smm{\rg=0}{b-1}\Imasup{\ktwo}{(\tpose{\Uovercc\kernsub{5}})^\rg\cdot\tpose{U_4}}
\rt)\rt)
\\
&\hskip2em\geq l_{\ktwo}\lt(
\lt(\dss{i=1}{p}\ktwo(n_{a+i}+2)\rt)/\prr\lt(
\smm{\rg=0}{b-1}\Imasup{\ktwo}{(\tpose{\Uovercc\kernsub{5}})^\rg\cdot\tpose{U_4}}
\rt)\rt)
\\
&\hskip2em\geq l_{\ktwo}\lt(
\lt(\dss{i=1}{p}\ktwo(n_{a+i}+2)\rt)/\prr(P_0)\rt)
\end{split}
\end{equation*}
by \eqref{eq227} and \eqref{eq228}. 
This length cannot be finite if $t(t-1)/2\leq p$, since $P_0$ 
is generated by $t(t-1)/2$ elements contained in 
$(x_3,x_4)\dss{i=1}{b}\ktwo(n_{a+i}+2)$ over $\ktwo$.
Hence $t(t-1)/2>p$.
In the remaining case (iii), we can obtain \eqref{eq228}
in the same way as above by \Lem{\ref{proc3}} and the case \eqref{cond100} 
of \Lem{\ref{proc44}}, and reach our assertion. 
\case{2}
In the case $a'=t+1$, we make use of a \Wei\ basis stated in \Lem{\ref{proc7}}
and can argue in exactly the same manner as above to obtain $t''(t''+1)/2>p$ if 
$n_{a+1}<n_{a'}$. Since $t''\leq t-1$, this implies $t(t-1)/2>p$.
\end{pf}

\begin{cor}\lb{proc119}
Assume $\ch{k}=0$. Suppose 
that $X$ is contained in an irreducible surface of degree $a\geq2$
and that $b\geq1$. Let $a'$ be an integer with $\ee{2}{a'}{a}$. Then $n_{a+1}\geq n_{a'}$,
if 
\begin{parenumr}
\item
$n_i=n_1+i-1$ for all $\ee{1}{i}{a'}$, or
\item
$a'\geq3$ and $n_i=n_1+i-2$ for all $\ee{2}{i}{a'}$.
\end{parenumr}
\end{cor}

\begin{pf}
Suppose that $n_{a+1}<n_{a'}$. We use \Thm{\ref{proc118}} with $m=1$. Then $t=t_1$
and $p\geq1$.
If $n_i=n_1+i-1$ for all $\ee{1}{i}{a'}$, then $t_1=1$ and $0=t_1(t_1-1)/2>p\geq1$
by \Thm{\ref{proc118}}. Likewise if 
$a'\geq3$ and $n_i=n_1+i-2$ for all $\ee{2}{i}{a'}$, then $t_1=2$ and 
$1=t_1(t_1-1)/2>p\geq1$. Thus we are led to a contadiction.
\end{pf}

In the remaining part of this section, we describe other miscellaneous results
that can be proved using the existence of an irreducible homogeneous 
polynomial of minimal degree.

\begin{lem}\lb{proc33}
Assume that $a\geq2$ and that $X$ is contained in an irreducible surface of degree $a$.
Suppose $n_i=n_1+i-1$ for all $\ee{1}{i}{a}$. Then, there is a \Wei{} basis 
$\{e^1_1,e^2_1\ddd e^2_a,e^3_1\ddd e^3_b\}$ of $I$ with respect to $\seqfour$
such that $U_{21}=0$. Moreover, it also have the properties
stated in \Lems{\ref{proc107} and \ref{proc109}} or in \Lem{\ref{proc7}} with 
$m=t=1$.
\end{lem}

\begin{pf}
Let $\{e^1_1,e^2_1\ddd e^2_a,e^3_1\ddd e^3_b\}$ be the \Wei{} basis 
of $I$ stated in \Lems{\ref{proc107} and \ref{proc109}}
with $m=t=1$ for $a\geq 3$, or in \Lem{\ref{proc7}} with $m=t=1$ for $a=2$.
Denote by $J$ the ideal $(e^1_1,e^2_1)\sset R$. 
Then, 
\begin{equation*}
U_{21}\clmlt{\ }{a}=0\quad\andt\quad
\rank_{k}\lt(U_{1}\clmlt{\ }{a}\ (\mod\ (\seqfour))\rt)=a-1.
\end{equation*}
These equalities imply that 
$e^2_2\ddd e^2_a\in J$, since $(e^1_1,e^2_1\ddd e^2_a,e^3_1\ddd e^3_b)\lg_2=0$.
Hence $J=(e^1_1,e^2_1\ddd e^2_a)$.
We have $n_{a+1}\geq n_a$ by \Cor{\ref{proc119}} and 
$a\leq n_1<n_2<\cdots<n_a$ by hypothesis, so that $\gradn{J}{a}=\gradn{I}{a}$.
Since there is no unit in the first two rows of $\lg_2$, the ideal $J$ is 
minimally generated by $e^1_1$ and $e^2_1$. Besides,
$\gradn{I}{a}$ contains an irreducible polynomial by hypothesis.
Consequently, $J$ must be a complete intersection of two homogeneous polynomials
of of degree $a$ and $n_1$. Moreover
\begin{equation*}
I':=Re^1_1+\lt(\smm{l=1}{a}\kthree e^2_l\rt)=
Re^1_1\op\lt(\dss{l=1}{a}\kthree e^2_l\rt)\sset J,
\end{equation*}
where $I'$ is considered as a $\kthree$-module here.
Since $n_l=n_1+l-1$ by hypothesis, it follows from the above equality that
\begin{equation*}
\dimk{\gradn{I'}{t}}=
\clm{t-a+3}{3}_{+}+\smm{l=1}{a}\clm{t-n_1-l+3}{2}_{+}
\end{equation*}
for $t\in\Zbf$. Here, given integers $\mg$ and $\ng$, $\clmsc{\mg}{\ng}_{+}:=\mg!/(\mg-\ng)!\ng!$
if $\mg\geq\ng\geq0$ and $\clmsc{\mg}{\ng}_{+}:=0$ otherewise. Similarly, we see 
\begin{equation*}
\dimk{\gradn{J}{t}}=
\clm{t-a+3}{3}_{+}+\smm{l=1}{a}\clm{t-n_1-l+3}{2}_{+}
\end{equation*}
for $t\in\Zbf$, since the basic sequence of $J$ is 
$(a;n_1,n_1+1\ddd n_1+a-1)$ by \Lem{\ref{proc32}}.
Hence the subset $I'$ coincides with the ideal $J$. 
This implies that $x_1e^2_a\in J=I'$, therefore 
\begin{equation*}
x_1e^2_a=g^1_1e^1_1+\smm{l=1}{a}g^2_le^2_l
\end{equation*}
with some $g^1_1\in R,\ g^2_l\in\kthree\ (\ee{1}{l}{a})$.
Thus the last column of $U_{21}$ must also be zero, which proves our assertion.
\end{pf}

\begin{lem}\lb{proc47}
Let 
\begin{equation*}
A=\begin{bmatrix}
0&\phg_1\\
\phg_2&\phg_3
\end{bmatrix}\quad\andt\quad
B=\begin{bmatrix}
0&\psg_1\\
\psg_2&\psg_3
\end{bmatrix}
\end{equation*}
be matrices with entries in $\ktwo$ such that $[A,B]=0$. 
Suppose that $\phg_i\ (\ee{1}{i}{3})$ and
$\psg_i\ (\ee{1}{i}{3})$ are homogeneous with
$\deg{\phg_1}=\deg{\psg_1}\leq\deg{\phg_3}=\deg{\psg_3}=1\leq\deg{\phg_2}=\deg{\psg_2}$,
$\deg{\phg_1}-\deg{\phg_3}=\deg{\phg_3}-\deg{\phg_2}$.
Then, there are $\phg,\ \psg\in\ktwo$ and a $2\tm2$ matrix $C$ with homogeneous
entries in $\ktwo$ such that $A=\phg C$ and $B=\psg C$.
\end{lem}

\begin{pf}
By the hypothesis $[A,B]=0$, we have a system of equations 
\begin{equation*}
\lt\{
\begin{aligned}
\phg_1\psg_2-\phg_2\psg_1&=0\\
\phg_2\psg_3-\phg_3\psg_2&=0\\
\phg_3\psg_1-\phg_1\psg_3&=0
\end{aligned}\rt.
\end{equation*}
which is symmetric with respect to the permutation of 
$\{1,2,3\}$. 
\case{1}
We consider first the case where $\deg \phg_i=\deg \psg_i=1$ for all $\ee{1}{i}{3}$.
Suppose $\gcd(\phg_1,\phg_2)=1$.
Then $(\psg_1,\psg_2)=c(\phg_1,\phg_2)$ with some $c\in k$ and the
equations $\phg_3(c\phg_1)-\phg_1\psg_3=0$ and $\phg_2\psg_3-\phg_3(c\phg_2)=0$
imply $\psg_3=c\phg_3$. Hence $B=cA$. 
Suppose $\gcd(\phg_1,\phg_3)=1$ or $\gcd(\phg_2,\phg_3)=1$. 
Then in the same way as above, we see $B=cA$ with some $c\in k$.
Suppose next that $\phg_1\neq0$, $\phg_2=c_2\phg_1$, $\phg_3=c_3\phg_1$ with 
$c_2,\ c_3\in k$. Then
the equations $\phg_1\psg_2-(c_2\phg_1)\psg_1=0$ and $(c_3\phg_1)\psg_1-\phg_1\psg_3=0$
imply $\psg_2=c_2\psg_1$ and $\psg_3=c_3\psg_1$ and we get
\begin{equation*}
A=\phg_1\begin{bmatrix}
0&1\\c_2&c_3
\end{bmatrix}\quad\andt\quad
B=\psg_1\begin{bmatrix}
0&1\\c_2&c_3
\end{bmatrix}.
\end{equation*}
If $\phg_2\neq0$, $\phg_1=c_1\phg_2$, $\phg_3=c_3\phg_2$ with 
$c_1,\ c_3\in k$, or $\phg_3\neq0$, $\phg_2=c_2\phg_3$, $\phg_1=c_1\phg_3$ with 
$c_2,\ c_1\in k$, then we obtain similar expressions as above.
Finally, if $\phg_1=\phg_2=\phg_3=0$, then $A=0B$.
\case{2}
Next we consider the case 
$\deg{\phg_1}=\deg{\psg_1}<\deg{\phg_3}=\deg{\psg_3}=1<\deg{\phg_2}=\deg{\psg_2}$.
Suppose $\phg_1\neq0$. Then $\phg_1$ must be an element of $k$, so that
$\psg_3=(\psg_1/\phg_1)\phg_3,\ \psg_2=(\psg_1/\phg_1)\phg_2$. Hence $B=(\psg_1/\phg_1)A$.
In the same way, $A=(\phg_1/\psg_1)B$ if $\psg_1\neq0$. Suppose $\phg_1=\psg_1=0$.
Then $(\phg_2,\phg_3)$ and $(\psg_2,\psg_3)$ are proportional.
Hence $A=\phg C$ and $B=\psg C$ with some $2\tm2$ matrix $C$ with homogeneous
entries in $\ktwo$ and $\phg,\ \psg\in\ktwo$.
\end{pf}

\begin{prop}\lb{proc12}
Suppose that $X$ is contained in an irreducible surface of degree $a$ and that
$b>0$. Suppose further that $n_i=n_1+i-1$ for all $\ee{1}{i}{a}$.
Then $b>2$.  Moreover, if there is an integer  $b_0$ with $\nn{0}{b_0}{b}$ such that
$n_{a+b_0}+1<n_{a+b_0+1}$, then $b_0>2$.
\end{prop}

\begin{pf}
We use the \Wei\ basis stated in \Lem{\ref{proc33}}. 
Let $v_j$ be the $j$-th column of $\tpose{U_4}$ for each $\ee{1}{j}{a}$, 
$P$ the $\ktwo$-module generated by $v_a$, and 
$\Xg$ the module over $\ktwo$ generated by the columns of 
$[\tpose{\Uovercc\kernsub{3}},\tpose{\Uovercc\kernsub{5}}]$. 
Since $U_{21}=0$, it follows from the relation $\lg_2\lg_3=0$ that
$[\tpose{U_3},\tpose{U_5}]=[\tpose{\Uovercc\kernsub{3}},\tpose{\Uovercc\kernsub{5}}]=
[\tpose{U_3},\tpose{\Uovercc\kernsub{5}}]=0$ and that $\Xg=0$.
Besides,  $v_j\in(\seqtwo)\Imasup{R}{\tpose{U_3},\tpose{U_5}}+\Imasup{R}{v_a}$
for all $\ee{1}{j}{a}$ as in \eqref{eq226} by \Lem{\ref{proc101}}.
Hence $(\tpose{\Uovercc\kernsub{5}})^\rg v_j\in
(\seqtwo)\Imasup{R}{\tpose{U_3},\tpose{U_5}}+
\Imasup{R}{(\tpose{\Uovercc\kernsub{5}})^\rg v_a}$.
Since $x_1^\mg x_2^\ng(\tpose{\Uovercc\kernsub{5}})^\rg v_a\equiv
(\tpose{\Uovercc\kernsub{3}})^\mg (\tpose{\Uovercc\kernsub{5}})^\ng
(\tpose{\Uovercc\kernsub{5}})^\rg v_a=
(\tpose{\Uovercc\kernsub{3}})^\mg (\tpose{\Uovercc\kernsub{5}})^{\ng+\rg}v_a$ modulo 
$(\seqtwo)\Imasup{R}{\tpose{U_3},\tpose{U_5}}$ for all $\mg\geq0,\ \ng\geq0$
by \Lem{\ref{proc117}}, there is a $\phg\in\sm{m\geq0,\ n\geq0}
(\tpose{\Uovercc\kernsub{3}})^m(\tpose{\Uovercc\kernsub{5}})^n P$ such that
$(\tpose{\Uovercc\kernsub{5}})^\rg v_j\equiv\phg$ modulo 
$(\seqtwo)\Imasup{R}{\tpose{U_3},\tpose{U_5}}$.
We find therefore  by \Cor{\ref{proc111}} that 
$(\tpose{\Uovercc\kernsub{5}})^\rg v_j-\phg\in\sm{m\geq0,\ n\geq0}
(\tpose{\Uovercc\kernsub{3}})^m(\tpose{\Uovercc\kernsub{5}})^n P$.
Hence
\begin{equation}
(\tpose{\Uovercc\kernsub{5}})^\rg v_j\in
\sm{m\geq0,\ n\geq0}\ktwo 
(\tpose{\Uovercc\kernsub{3}})^m(\tpose{\Uovercc\kernsub{5}})^n v_a
\lb{eq230}
\end{equation}
for all $\rg\geq0$ and $\ee{1}{j}{a}$.
Let $b'$ denote $b_0$ if there is a $b_0$ with $\nn{0}{b_0}{b}$ such that
$n_{a+b_0}+1<n_{a+b_0+1}$, otherwise let $b':=b$.
Then, there are $b'\tm b'$ matrices 
$\Uovercc{}'\kernnsub{3}$ and $\Uovercc{}'\kernnsub{5}$ with
homogeneous components in $\ktwo$
such that
\begin{equation*}
\tpose{\Uovercc\kernsub{3}}=\begin{bmatrix}
\tpose{\Uovercc{}'\kernnsub{3}}&0\\
\ast&\ast
\end{bmatrix},\quad
\tpose{\Uovercc\kernsub{5}}=\begin{bmatrix}
\tpose{\Uovercc{}'\kernnsub{5}}&0\\
\ast&\ast
\end{bmatrix}
\end{equation*}
for reasons of the degrees of the components.
Notice that $[\tpose{\Uovercc{}'\kernnsub{3}},\tpose{\Uovercc{}'\kernnsub{5}}]=0$.
Let $\prr:\dss{i=1}{b}\ktwo(n_{a+i}+2)\lra \dss{i=1}{b'}\ktwo(n_{a+i}+2)$
be the natural projection to the first $b'$ components,
$v'_j:=\prr(v_j)$ ($\ee{1}{j}{a}$) and $\tpose{U'_4}:=(v'_1\ddd v'_a)$.
Then
\begin{equation}
\len{R}{\Ext{3}{R}{R/I}{R}}
\geq
l_\ktwo\lt(\dss{i=1}{b'}\ktwo(n_{a+i}+2)/
\sm{\rg\geq0}\Imasup{\ktwo}{\tpose{\Uovercc{}'\kernnsub{5}})^\rg\cdot\tpose{U'_4}}\rt)
\lb{eq231}
\end{equation}
by \Lem{\ref{proc200}}.
Suppose $b'=1$ or $b'=2$.
When $b'=1$, the matrices $\tpose{\Uovercc{}'\kernnsub{3}}$
and $\tpose{\Uovercc{}'\kernnsub{5}}$ are just elements of $\ktwo$, so that
\begin{equation*}
\tpose{\Uovercc{}'\kernnsub{5}^\rg} v'_j\in\sm{m\geq0,\ n\geq0}\ktwo 
(\tpose{\Uovercc{}'\kernnsub{5}})^m(\tpose{\Uovercc{}'\kernnsub{3}})^n v'_a
\sset\ktwo v'_a
\end{equation*}
by \eqref{eq230}.
Since $v'_a$ is also merely an element of $(\seqtwo)\ktwo$,
the length of $\Ext{3}{R}{R/I}{R}$ must be infinite by \eqref{eq231}, which is a contradiction.
When $b'=2$, let
\begin{equation*}
\tpose{\Uoverccnatural{3}}:=g_1 1_2-\tpose{\Uovercc{}'\kernnsub{3}}=
\begin{bmatrix}
0&\phg_1\\\phg_2&\phg_3
\end{bmatrix},\quad
\tpose{\Uoverccnatural{5}}:=g_2 1_2-\tpose{\Uovercc{}'\kernnsub{5}}=
\begin{bmatrix}
0&\psg_1\\\psg_2&\psg_3
\end{bmatrix},
\end{equation*}
where $g_1$ (resp. $g_2$) denotes the $(1,1)$ component of
$\Uovercc{}'\kernnsub{3}$ (resp. $\Uovercc{}'\kernnsub{5}$). 
Since the degrees of the entries of $\Uoverccnatural{3}$ and
$\Uoverccnatural{5}$ are determined by $n_{a+1}$ and $n_{a+2}$,
one sees that $\phg_i\ (\ee{1}{i}{3})$ and $\psg_ i\ (\ee{1}{i}{3})$ are homogeneous 
elements of $\ktwo$ such that 
$\deg{\phg_i}=\deg{\psg_j}=1$ for all $i,j$ or
$\deg{\phg_1}=\deg{\psg_1}<\deg{\phg_3}=\deg{\psg_3}=1<\deg{\phg_2}=\deg{\psg_2}$,
$\deg{\phg_1}-\deg{\phg_3}=\deg{\phg_3}-\deg{\phg_2}$.
Moreover, $[\tpose{\Uoverccnatural{3}},\tpose{\Uoverccnatural{5}}]=0$
by $[\tpose{\Uovercc{}'\kernnsub{3}},\tpose{\Uovercc{}'\kernnsub{5}}]=0$.
By \Lem{\ref{proc47}},
there are $\phg,\ \psg\in\ktwo$ and a $2\tm2$ matrix $C$ with homogeneous
entries in $\ktwo$ such that $\tpose{\Uoverccnatural{3}}=\phg C$ and 
$\tpose{\Uoverccnatural{5}}=\psg C$. We find therefore
\begin{equation*}
\begin{split}
(\tpose{\Uovercc{}'\kernnsub{5}})^\rg v'_j
&\in\sm{m\geq0,\ n\geq0}\ktwo 
(\tpose{\Uovercc{}'\kernnsub{3}})^m(\tpose{\Uovercc{}'\kernnsub{5}})^n v'_a\\
&=\sm{m\geq0,\ n\geq0}\ktwo 
(g_1 1_2-\tpose{\Uoverccnatural{3}})^m(g_2 1_2-\tpose{\Uoverccnatural{5}})^n v'_a\\
&\sset\sm{m\geq0,\ n\geq0}\ktwo 
(\tpose{\Uoverccnatural{3}})^m(\tpose{\Uoverccnatural{5}})^n v'_a
\sset\sm{m\geq0,\ n\geq0}\ktwo 
C^{m+n} v'_a\\
&=\ktwo v'_a+\ktwo C v'_a
\end{split}
\end{equation*}
for all $\rg\geq0$ and $\ee{1}{j}{a}$ by \eqref{eq230} and Hamilton-Cayley's formula.
Since the components of $v'_a$ lie in $(x_3,x_4)\ktwo$, 
the length of $\Ext{3}{R}{R/I}{R}$ must be infinite by \eqref{eq231}.
In any case we are let to a contradiction. Thus $b'>2$.
\end{pf}

\begin{prop}\lb{proc34}
Suppose that $X$ is contained in an irreducible surface of degree $a$ and that
$b>1$.  Suppose further that
$n_i=n_1+i-1$ for all $\ee{1}{i}{a}$. Then, $n_{a+j}\neq n_{a+1}+j-1$
for some $\ee{1}{j}{b}$. Moreover, if there is an integer $b_0$ with
$\nn{1}{b_0}{b}$ such that $n_{a+b_0}+1<n_{a+b_0+1}$, then $n_{a+j}\neq n_{a+1}+j-1$
for some $\ee{1}{j}{b_0}$.
\end{prop}

\begin{pf}
We use the \Wei\ basis and the standard relation matrix $\lg_2$ stated in \Lem{\ref{proc33}}. 
Let $b''$ denote $b_0$ if there is a $b_0$ with $\nn{0}{b_0}{b}$ such that
$n_{a+b_0}+1<n_{a+b_0+1}$, otherwise let $b'':=b$.
Then, there are $b''\tm b''$ matrices 
$U''_3$ and $U''_5$ with homogeneous components in $R$ such that
\begin{equation*}
\tpose{U_3}=\begin{bmatrix}
\tpose{U''_3}&0\\
\ast&\ast
\end{bmatrix},\quad
\tpose{U_5}=\begin{bmatrix}
\tpose{U''_5}&0\\
\ast&\ast
\end{bmatrix}
\end{equation*}
for reasons of the degrees of the components.
Suppose that $n_{a+j}=n_{a+1}+j-1$
for all $\ee{1}{j}{b''}$. 
Let $b'$ be the minimum of $j$ $(\en{1}{j}{b''})$ such that the $(j,j+1)$ components of
$\tpose{U''_3}$ and $\tpose{U''_5}$ are both zero.  In case there is no such $j$, let
$b':=b''$.
Since a unit can appear only as a $(j,j+1)$ component 
$(\en{1}{j}{b''})$ in $\tpose{U''_3}$ and $\tpose{U''_5}$,
one has
\begin{align*}
&\tpose{U''_3}\clmlt{b'+1\ddd b''}{1,2,\ldots,b'}=
\tpose{U''_5}\clmlt{b'+1\ddd b''}{1,2,\ldots,b'}=0\quad\andt\\
&\rankr_k\lt(\tpose{U''_3}\clmlt{b'+1\ddd b''}{b'+1\ddd b''}
+s\tpose{U''_5}\clmlt{b'+1\ddd b''}{b'+1\ddd b''}\rt)\ 
(\mod\ (\seqfour))\ =b'-1
\end{align*}
for some $s\in k$.
Put $\tpose{U'_3}:=\tpose{U''_3}\clmlsc{b'+1\ddd b''}{b'+1\ddd b''}$ and
$\tpose{U'_5}:=\tpose{U''_5}\clmlsc{b'+1\ddd b''}{b'+1\ddd b''}$.
Then, they are $b'\tm b'$ matrices with homogeneous components in $R$ such that
\begin{align}
&\tpose{U_3}=\begin{bmatrix}
\tpose{U'_3}&0\\
\ast&\ast
\end{bmatrix},\quad
\tpose{U_5}=\begin{bmatrix}
\tpose{U'_5}&0\\
\ast&\ast
\end{bmatrix},\notag\\
&\rk{k}{\tpose{U'_3}+s\tpose{U'_5}}\ (\mod\ (\seqfour))\ =b'-1.
\lb{eq11}
\end{align}
Since $U_{21}=0$, one has $[\tpose{U_3},\tpose{U_5}]=0$, so that 
$[\tpose{U'_3},\tpose{U'_5}]=[\tpose{U'_3}+s\tpose{U'_5},\tpose{U'_5}]=0$.
Taking the degrees of the components into account, one sees that 
\begin{align*}
\Imasup{R}{\tpose{U'_3},\tpose{U'_5}}&=
\Imasup{R}{\tpose{U'_3}+s\tpose{U'_5},\tpose{U'_5}}\\
&=\Imr^R\lt(\tpose{U'_3}+s\tpose{U'_5},\tpose{U'_5}
\clmlt{\phantom{1} }{1,2,\ldots,b'-1}\rt).
\end{align*}
On the other hand, since the $j$-th column $v_j$ of $\tpose{U_4}$ is contained in
$(\seqtwo)\Imasup{R}{\tpose{U_3},\tpose{U_5}}+\Imasup{R}{v_a}$
for all $\ee{1}{j}{a}$ by \Lem{\ref{proc101}}, one has
$\Imasup{R}{\tpose{U_3},\tpose{U_5},\tpose{U_4}}=
\Imr^R\lt(\tpose{U_3},\tpose{U_5},v_a\rt)$.
Let $\prr:\dss{i=1}{b}\ktwo(n_{a+i}+2)\lra \dss{i=1}{b'}\ktwo(n_{a+i}+2)$
be the natural projection to the first $b'$ components,
$v'_j:=\prr(v_j)$ ($\ee{1}{j}{a}$) and $\tpose{U'_4}:=(v'_1\ddd v'_a)$,
where $\tpose{U_4}=(\seq{v}{1}{a})$. Then, by what we have seen,
\begin{equation*}
\prr(\Imasup{R}{\tpose{U_3},\tpose{U_5},\tpose{U_4}})=
\Imr^R\lt(\tpose{U'_3}+s\tpose{U'_5},\tpose{U'_5}
\clmlt{\phantom{1}}{1,2,\ldots,b'-1},v'_a\rt).
\end{equation*}
By \eqref{eq11}, therefore,
\begin{equation*}
\prr(\Coksup{R}{\tpose{U_3},\tpose{U_5},\tpose{U_4}})\cong
R/(g,g',g'')
\end{equation*}
with suitable homogeneous polynomials $g,g',g''$. Since the length of
this module is infinite and since
\begin{equation*}
\begin{split}
\len{R}{\Ext{3}{R}{R/I}{R}}&=\len{R}{\Coksup{R}{\tpose{U_3},\tpose{U_5},\tpose{U_4}}}\\
&\geq\len{R}{\prr(\Coksup{R}{\tpose{U_3},\tpose{U_5},\tpose{U_4}})},
\end{split}
\end{equation*}
we are led to a contradiction. Hence $n_{a+j}\neq n_{a+1}+j-1$
for some $\ee{1}{j}{b''}$. This proves our assertion.
\end{pf}

\section{\mathversion{bold}The case $b=1$}
\lb{sec5}

Let the notation be the same as in Section \ref{sec4}.
We can give a fairly precise characterization of the basic sequences 
of integral curves with $b=1$.

\begin{lem}\lb{proc36}
Assume that $b=1$ and that $X$ does not contain any line as an
irreducible component.
Let $g_1:=U_3$, $g_2:=U_5$. Then 
\begin{equation*}
\rankr_{R/(g_1,g_2)}\lt(
\begin{bmatrix}
U_{01}\\
U_1\\
U_{21}
\end{bmatrix}
\quad(\mod\ (g_1,g_2))\quad
\rt)\geq a-1.
\lb{eq77}
\end{equation*}
\end{lem}

\begin{pf}
Let $\{e^1_1,e^2_1\ddd e^2_a,e^3_1\}$ be a \Wei\ basis of $I$ with respect to 
$\seqfour$ and let $\lg_1$ and $\lg_2$ be as in Section \ref{sec1}.
Let $\pt$ be a point of $\Pthree$.
Since $X$ is locally \CM, the ideal sheaf $\Ical_X$ of $X$ has a free resolution
\begin{equation*}
0\lra\Ocal_X^{a+1}|_\Wcal\xra{\ \lgbar\ }
\Ocal_X^{a+2}|_\Wcal\xra{\ \lg_1\ }\Ical_X|_\Wcal\lra0
\end{equation*}
on a neighborhood $\Wcal$ of $\pt$, where $\lgbar$ is a matrix obtained from
$\lg_2$ by deleting one column. By Hilbert-Burch
theorem, each $e^i_l$ is an $(a+1)\tm(a+1)$ minor of the 
$(a+2)\tm(a+1)$ matrix $\lgbar$ on $\Wcal$. If our assertion did not hold,
then 
\begin{equation*}
\rankr_{R/(g_1,g_2)}\lt(
\ \lg_2\ (\mod\ (g_1,g_2))\ \rt) < a+1.
\end{equation*}
This would mean that $e^i_l\equiv0\ (\mod\ (g_1,g_2))$ for all $i,\ l$, i.e. 
$\Ical_X\sset(g_1,g_2)\Ocal_X$, contradicting the assumption that 
$X$ does not contain any line as an irreducible component.
\end{pf}

\begin{lem}\lb{proc37}
Assume that $b=1$ and that $X$ is contained in an irreducible 
surface of degree $a$. Let $(\seq{v}{1}{a}):=\tpose{U_4}$,
$p_i:=\deg(v_i)=n_{a+1}+1-n_i \ (\ee{1}{i}{a})$, and 
\begin{equation}
0\lra\dss{j=1}{a-1}\ktwo(-q_j)\xra{\ F\ }
\dss{i=1}{a}\ktwo(-p_i)\xra{\ \tpose{U_4}\ }
\afrak\lra0
\lb{eq78}
\end{equation}
a free resolution of $\afrak:=\Imasup{\ktwo}{\tpose{U_4}}$,
where $p_1\geq p_2\geq\cdots\geq p_a$ and 
$\seq{q}{1}{a-1}$ are suitable integers with 
$q_1\geq q_2\geq\cdots\geq q_{a-1}$.
Let further 
$n'_1\ddd n'_\og$ $(\og\geq1)$ be a strictly increasing sequence of integers
such that $\{n'_1\ddd n'_\og\}=\{\ n_i\ |\ \ee{1}{i}{a}\ \}$ and 
let $i_\gg:=\min\{\ i\ |\ n_i=n'_\gg,\ \ee{1}{i}{a}\ \}$ for each 
$\ee{1}{\gg}{\og}$.
Denote by $\Delta_{F}$ (resp. $\Delta_{U_1}$) 
the $a\tm(a-1)$ matrix (resp. $a\tm a$ matrix) whose
$(i,j)$ component is $q_j-p_i$ (resp. $n_j+1-n_i$).
Then, for each $\gg$ $(\ee{2}{\gg}{\og})$, there is a $j_\gg\ (\ee{1}{j_\gg}{a-1})$
such that the transpose of the $j_\gg$-th column of $\Delta_F$
coincides with the $i_\gg$-th row of $\Delta_{U_1}$.
\end{lem}

\begin{pf}
Let $g_1:=U_3$, $g_2:=U_5$.
Since $g_1$ and $g_2$ are linear forms with $g_1-x_1, g_2-x_2\in\ktwo$,
we can write
\begin{align}
&U_{01}=g_2U'_{01}+U''_{01},\quad
U_1=g_11_a+g_2U'_1+U''_1,
\lb{eq79}\\
&\witht\quad
U'_{01},\ U'_1\in\MATr(\kthree),\quad 
U''_{01},\ U''_1\in\MATr(\ktwo).
\notag
\end{align}
Let $\gg$ be an integer with $\ee{2}{\gg}{\og}$.
By \Lem{\ref{proc26}}, 
there must exist an integer $i'_\gg$ $(\ee{1}{i'_\gg}{a})$
with $n_{i'_\gg}=n_{i_\gg}=n'_{\gg}$
such that the $i'_\gg$-th row of $U_1$ contains a unit.
Since the unit of $U_1$ lies in $U''_1$, this means that
the $i'_\gg$-th row of $U''_1$ contains a unit.
On the other hand, it follows from the relation $\lg_2\lg_3=0$ 
that $U''_1U_4=0$, i.e. $\tpose{U_4}\tpose{U''_1}=0$.
The $i'_\gg$-th column of $\tpose{U''_1}$ is therefore a
linear combination over $\ktwo$ of the columns of $F$
by \eqref{eq78}. Since the $i'_\gg$-th row contains a unit, 
this implies that there is a $j_\gg\ (\ee{1}{j_\gg}{a-1})$
such that the $j_\gg$-th column of $F$ contains a unit.
Consequently, the $j_\gg$-th column of $\Delta_F$
coincides with the $i'_\gg$-th column of 
$\Dg(\tpose{U''_1})=\Dg(\tpose{U_1})=\tpose{\Delta_{U_1}}$.
Hence the 
transpose of the $j_\gg$-th column of $\Delta_F$
coincides with the $i_\gg'$-th row of $\Delta_{U_1}$.
Besides, the $i_\gg$-th row and the $i'_\gg$-th row of 
$\Delta_{U_1}$ are the same. Thus the 
transpose of the $j_\gg$-th column of $\Delta_F$
coincides with the $i_\gg$-th row of $\Delta_{U_1}$.
\end{pf}

\begin{lem}\lb{proc38}
Assume that $b=1$, $X$ is contained in an irreducible 
surface of degree $a$, and that $X$ does not contain any line as an
irreducible component. Denote by $\Delta_{U_{01}}$ and
$\Delta_{U_{21}}$ the row vectors $(n_1+1-a,n_2+1-a,\ldots,n_a+1-a)$
and $(n_1+1-n_{a+1},n_2+1-n_{a+1},\ldots,n_a+1-n_{a+1})$ respectively.
With the notation of \Lem{\ref{proc37}},
let $\Delta^{(1)}_{U_1}:=\Delta_{U_1}\clmlsc{i_2,i_3,\ldots,i_\og}{\phantom{1}}$
and
$\Delta^{(1)}_F=(\phi_{ij}):=\Delta_F\clmlsc{\phantom{1}}{j_2,j_3,\ldots,j_\og}$.
Rearrange the rows of the matrix
\begin{equation*}
\begin{bmatrix}
\Delta_{U_{01}}\\
\Delta^{(1)}_{U_1}\\
\Delta_{U_{21}}
\end{bmatrix}
\end{equation*}
in such a way that the last column becomes a decreasing sequence
of integers, and denote the resulting matrix by $\Delta=(\dg_{ij})$, i.e.
$\dg_{1a}\geq \dg_{2a}\geq\cdots\geq \dg_{a+3-\og\,a}$.
Then $\dg_{ia}\geq\phi_{ai}$ for all $\ee{1}{i}{a-\og}$.
\end{lem}

\begin{pf}
Let the notation be as in the proof of \Lem{\ref{proc37}}.
The $i'_\gg$-th column of $\tpose{U''_1}$ is a linear combination
of the columns of $F$ over $\ktwo$ in which the coefficient of the
$j_\gg$-th column of $F$ is a nonzero constant. 
Hence we may assume that the $j_\gg$-th column of $F$
coincides with the $i'_\gg$-th column of $\tpose{U''_1}$ for all $\ee{2}{\gg}{\og}$,
after a suitable column operations on $F$ if necessary.
Let $U''{}^{(2)}_1$ be the $(\og-1)\tm a$ matrix whose $(\gg-1)$-th row
is the $i'_\gg$-th row of $U_1$ for all $\ee{2}{\gg}{\og}$
and $F^{(2)}$ be the $a\tm(\og-1)$ matrix whose $(\gg-1)$-th column
is the $j_\gg$-th column of $F$ for all $\ee{2}{\gg}{\og}$.
Further, let 
$U''{}^{(1)}_1:=U''_1\clmlsc{i'_2,i'_3,\ldots,i'_\og}{\phantom{1}}$ and
$F^{(1)}:=F\clmlsc{\phantom{1}}{j_2,j_3,\ldots,j_\og}$.
Note that $\Dg(F^{(1)})=\Delta^{(1)}_F$, $\tpose{U''{}^{(2)}_1}=F^{(2)}$, and that
$F=(F^{(2)},F^{(1)})$ 
(resp. $U_1''=\lt[\begin{smallmatrix}
U''{}^{(2)}_1\\
U''{}^{(1)}_1
\end{smallmatrix}\rt]$)
up to permutation of the columns (resp. rows).
Let $\Thg$ be the invertible matrix representing a
permutation of $\{1,2,\ldots,a+3-\og\}$ such that
\begin{equation*}
\Thg\begin{bmatrix}
\Delta_{U_{01}}\\
\Delta^{(1)}_{U_1}\\
\Delta_{U_{21}}
\end{bmatrix}=\Delta,
\quad\text{and let}\quad 
V'':=\Thg\begin{bmatrix}
U''_{01}\\
U''{}^{(1)}_1\\
U_{21}
\end{bmatrix},\quad
U'':=\begin{bmatrix}
U''{}^{(2)}_1\\
V''
\end{bmatrix}.
\end{equation*}
We have $V'',\ U''\in\MATr(\ktwo)$. Moreover, $U''$ is obtained from
\begin{equation*}
\begin{bmatrix}
U''_{01}\\
U''_1\\
U_{21}
\end{bmatrix}
\end{equation*}
by rearranging its rows.
Since the $i_\gg$-th row and the 
$i'_\gg$-th row of $\Delta_{U_1}$ are the same,
$\Dg(U''{}^{(1)}_1)=\Delta^{(1)}_{U_1}$, so that $\Dg(V'')=\Delta$.
Besides, $V''U_4=0$, i.e. $\tpose{U_4}\tpose{V''}=0$, by  the relation $\lg_2\lg_3=0$.
Each column of $\tpose{V''}$ is therefore a linear combination of
the columns of $F$ over $\ktwo$ by \eqref{eq78}, so that
each row of $V''$ is a linear combination of
the rows of $\tpose{F}$ over $\ktwo$.
Suppose there is a $t$ with $\ee{1}{t}{a-\og}$ such that
$\dg_{ta}<\phi_{at}$. Then, comparing the degrees,
we find that, for every $i\geq t$, the $i$-th row of $V''$ is 
a linear combination over $\ktwo$ of the rows of
$\tpose{\lt(F^{(1)}\clmlsc{\phantom{1}}{1,2,\ldots,t}\rt)}$ 
and $\tpose{F^{(2)}}$. Moreover, the rows of $U''{}^{(2)}_1$
are the rows of $\tpose{F^{(2)}}$.
This implies that the rows of $U''$ are linear combinations over 
$\ktwo$ of the first $t-1$ rows of $V''$,
the $a-\og-t$ rows of $\tpose{\lt(F^{(1)}\clmlsc{\phantom{1}}{1,2,\ldots,t}\rt)}$ 
and of the $\og-1$ rows of $\tpose{F^{(2)}}$. Hence
\begin{multline*}
\rankr_{R/(g_1,g_2)}
\lt(\begin{bmatrix}
U_{01}\\
U_1\\
U_{21}
\end{bmatrix}
\quad(\mod\ (g_1,g_2))\quad
\rt)\\
=\rk{\ktwo}{U''}
\leq (t-1)+(a-\og-t)+(\og-1)=a-2,
\end{multline*}
which contradicts \Lem{\ref{proc36}}.
Thus $\dg_{ia}\geq\phi_{ai}$ for all $\ee{1}{i}{a-\og}$.
\end{pf}

\begin{thm}\lb{proc39}
Assume that $b=1$, $X$ is contained in an irreducible 
surface of degree $a$, and that $X$ does not contain any line as an
irreducible component.  
Let $n'_1\ddd n'_\og$ and $i_\gg$ $(\ee{1}{\gg}{\og})$ be as in \Lem{\ref{proc37}}.
Delete the terms $n'_2\ddd n'_\og$ from the sequence $(a,\seq{n}{1}{a},n_{a+1})$
and then rearrange the remaining terms so that they make a nondecreasing sequence.
Let $(n''_1\ddd n''_{a+3-\og})$ denote the resulting sequence.
Then 
\begin{equation*}
n_{a+1}\leq
a-2+\smm{j=1}{a}n_j -\smm{\gg=2}{\og}n'_\gg
-\smm{i=1}{a-\og}n''_i.
\end{equation*}
\end{thm}

\begin{pf}
Let the notation 
be as in \Lems{\ref{proc37} and \ref{proc38}}.
Denote by $(\xg_{\gg 1}\ddd \xg_{\gg a})$ the $i_\gg$-th row of 
$\Delta_{U_1}$ for each $\ee{2}{\gg}{\og}$.
Then, we have an exact sequence of the form
\begin{equation*}
0\lra\dss{j=1}{a-1}\ktwo(-q'_j)\xra{\ (F^{(2)}, F^{(1)})\ }
\dss{i=1}{a}\ktwo(-p_i)\xra{\ \tpose{U_4}\ }
\afrak\lra0
\end{equation*}
by \eqref{eq78}, where $(q'_1\ddd q'_{a-1})=(\seq{q}{1}{a-1})$
up to permutation.
Let $\Delta^{(2)}_{U_1}$ (resp. $\Delta^{(2)}_F$) be the 
$(\og-1)\tm a$ (resp. $a\tm(\og-1)$) matrix whose 
$(\gg-1)$-th row (resp. column)
is the $i_\gg$-th row (resp. $j_\gg$-th column) of $\Delta_{U_1}$
(resp. $\Delta_F$) for all $\ee{2}{\gg}{\og}$. 
Then, $\Delta^{(2)}_F=\tpose{\Delta^{(2)}_{U_1}}$ by \Lem{\ref{proc37}},
so that
\begin{equation*}
\Dg(F^{(2)}, F^{(1)})=(\Delta^{(2)}_F,\Delta^{(1)}_F)=
\begin{bmatrix}
\xg_{21}&\cdots&\xg_{\og1}&\phi_{11}&\cdots&\phi_{1\,a-\og}\\
\vdots&\cdots&\vdots&\vdots&\cdots&\vdots\\
\xg_{2a}&\cdots&\xg_{\og a}&\phi_{a1}&\cdots&\phi_{a\,a-\og}
\end{bmatrix}.
\end{equation*}
By Hilbert-Burch theorem, the first component
of $\tpose{U_4}$ is the determinant of the $(a-1)\tm(a-1)$ matrix
$ (F^{(2)}, F^{(1)})\clmlsc{1}{\phantom{1}}$ up to multiplication by a constant.
Comparing the degrees, we see 
\begin{equation}
n_{a+1}+1-n_1=p_1=\smm{\gg=2}{\og}\xg_{\gg\,\gg}
+\smm{i=1}{a-\og}\phi_{\og+i\,i}.
\lb{eq81}
\end{equation}
Observe that
$\dg_{ia}+p_a=\dg_{ij}+p_j$, $\phi_{ai}+p_a=\phi_{ji}+p_j$
for all $\ee{1}{i}{a-\og}$ and $\ee{1}{j}{a}$, so that $\dg_{ij}\geq\phi_{ji}$
by \Lem{\ref{proc38}}. With the use of this inequality, we find 
by \eqref{eq81} that 
\begin{equation*}
n_{a+1}+1-n_1\leq\smm{\gg=2}{\og}\xg_{\gg\,\gg}
+\smm{i=1}{a-\og}\dg_{i\,\og+i}.
\end{equation*}
Since the components of $\Delta_{U_{01}}$, $\Delta_{U_{21}}$ and $\Delta_{U_1}$ are
as stated in \Lems{\ref{proc37} and \ref{proc38}}, we see that
$\xg_{\gg j}=n_j+1-n_{i_\gg}=n_j+1-n'_\gg$, $\dg_{ij}=n_j+1-n''_i$ for all $\gg,i,j$.
Putting them into the above inequality, we get our assertion.
\end{pf}

\begin{rem}\lb{proc48}
\no{i}
One needs rearrangements 
as stated in the above theorem only in case $n_{a+1}<n_i$ for some $\ee{1}{i}{a}$.
\indno{ii}
In section 4 of the paper \cite{A2}, the case $U_{21}=0$ is treated  
in detail. Since we do not assume it in this section, the results
loco citato are somewhat different from those given here.
\end{rem}

\section{Comparison of our results with Cook's assertions by numerical
computations with the help of a computer}
\lb{sec7}

Let $X$ be a curve in $\Pthree$ and $I$ the saturated homogeneous ideal of $X$
in $R$. Assume that $\ch{k}=0$. 
Let $\rg$ be an arbitrary nonnegative integer.
With the notation of \Rem{\ref{proc40}},
we consider the strongly Borel fixed monomial ideal in $k[x_1,x_2]$ generated by 
\begin{eqnarray}
\{\ x_1^a,x_1^{a-1}x_2^{\beta_{a-1}}\ddd x_1x_2^{\beta_1},x_2^{\beta_0}\ \}
\cup
\{\ x_1^{t_l}x_2^{u_l}\ |\ 
\ee{1}{l}{b},\ 
\beta_{t_lu_l}\leq \rg\ \},
\label{eq84}
\end{eqnarray}
and let 
$x_1^{s_\rg},x_1^{s_\rg-1}x_2^{\mu_{s_\rg-1}(\rg)}
\ddd x_1x_2^{\mu_1(\rg)},x_2^{\mu_0(\rg)}$
$(1\leq\mu_{s_\rg-1}(\rg)<\mu_{s_\rg-2}(\rg)<\cdots<\mu_0(\rg))$ be
its minimal generators, where the linear forms $\seqfour$ are chosen
sufficiently generally. Put 
$\nu_l(\rg):=\deg(x_1^{s_\rg-l}x_2^{\mu_{s_\rg-l}(\rg)})=s_\rg-l+\mu_{s_\rg-l}(\rg)$
for $\ee{1}{l}{s_\rg}$. 
We have $\nu_l(\rg)\leq\nu_{l+1}(\rg)$ for all $\en{1}{l}{s_\rg}$ and 
$s_\rg\leq\nu_1(\rg)$.
Now assume further that $X$ is integral.
With our notation, the assertions
in the main theorem of \cite{C} can be formulated in the 
following manner:
\begin{enumr}
\item\label{cond19}
$\nu_l(\rg)\leq\nu_{l+1}(\rg)\leq\nu_l(\rg)+1$ for all $\en{1}{l}{s_\rg}$,
\item\label{cond20}
if further $s_\rg<a$, then $s_\rg\leq\nu_1(\rg)\leq s_\rg+1$.
\end{enumr}
Unfortunately, there is an error in the proof of the above assertions
(see \cite[Section 4]{DS}).
In \cite{A10}, however, a proof of \eqref{cond19} is given for a special case.

\begin{thm}[{\cite[\Thm{5}]{A10}}]\label{proc112}
With the notation above, the assertion \eqref{cond19} 
is true for all $\rg$ such that $s_\rg=s_0$.
\end{thm}

\par
Let $d(X)$ and $g(X)$ denote the degree and the arithmetic genus of $X$ respectively,
and let $(a;n_1\ddd n_a;\seq{n}{a+1}{a+b})$ be the basic sequence of $I$.
Since
\begin{multline*}
1-g(X)+\ng d(X)=\dimk{\gradn{R/I}{\ng}}=\dimk{\gradn{R}{\ng}}-\dimk{\gradn{I}{\ng}}\\
=\binom{\ng+3}{3}
-\binom{\ng-a+3}{3}-\smm{l=1}{a}\binom{\ng-n_l+2}{2}
-\smm{l=1}{b}\binom{\ng-n_{a+l}+1}{1}
\end{multline*}
for all $\ng\ggreater0$, direct computation shows that 
\begin{align*}
&d(X)=\smm{l=1}{a}n_l-\frac{1}{2}a(a-1)-b,\\
&g(X)=1+\smm{l=1}{a}\frac{1}{2}n_l(n_l-3)-\smm{l=1}{b}n_{a+l}+b
-\frac{1}{6}a(a-1)(a-5)
\end{align*}
 (see \cite[\Rem{1.9}]{A3}). 
We want to consider sequences satisfying these equalities.

\par
Let $(d,g)$ be a pair of positive integer $d$ and a nonnegative integer $g$,
and $\Bs=(a;n_1\ddd n_a;\seq{n}{a+1}{a+b})$  be a sequence of 
integers such that 
\begin{equation}
\begin{gathered}
0<a\leq n_1\leq\cdots\leq n_a,\ 
n_1\leq n_{a+1}\leq\cdots\leq n_{a+b},\ b\geq0,\\ 
a\leq d-1,\ n_l\leq d-1\ (\ee{1}{l}{a+b}),\ d=D(\Bs),\ g=G(\Bs).
\end{gathered}
\lb{eq101}
\end{equation}
See Introduction for the definitions of $D(\Bs)$ and $G(\Bs)$. 
When $\Bs=\BR{I}$, we have $d(X)=D(\BR{I})$ and $g(X)=G(\BR{I})$ by the above formulas.
Moreover, $a\leq d(X)-1$ and $n_l\leq d(X)-1$
for all $\ee{1}{l}{a+b}$ by Castelnuovo's regularity theorem combined with
the results of \cite{BS}.
Hence $\BR{I}$ satisfies \eqref{eq101}.
Moreover the generic initial ideal $\inr^x(I)$ is strongly Borel fixed by \cite{U}.
With this in mind, we shall say that a sequence $\Bs$ 
satisfies A-conditions for $(d,g)$, if there is a strongly Borel fixed monomial
ideal in $k[\seqfour]$ giving the sequence $\Bs$ as in \Rem{\ref{proc40}}, and if
$\Bs$ satisfies \eqref{eq101} and all the numerical conditions described in
\Lem{\ref{proc26}}, \Thm{\ref{proc118}}, \Cor{\ref{proc119}}, 
\Props{\ref{proc12}, \ref{proc34}},
and \Thm{\ref{proc39}}.
Likewise, we shall say that $\Bs$ satisfies C-conditions for $(d,g)$,
if $\Bs$ satisfies \eqref{eq101} and is given by a strongly Borel fixed monomial
ideal which fulfills Cook's assertions \eqref{cond19} and \eqref{cond20}.
We can find all the sequences $\Bs$
satisfying A-conditions and C-conditions respectively
with the help of a computer. 
We can find the sequences satisfying both conditions, too.

\par
Since there are only  a finite number of possibilities of $a,\ b$, $n_l\ (\ee{1}{l}{a+b})$ 
with the property \eqref{eq101} for a fixed $(d,g)$, one can extract from them all $\Bs$'s 
satisfying the numerical conditions described in the lemma, the propositions, the theorems 
and the corollary
mentioned just above, carrying out simple numerical testings.

\par
In order to get the generators of strongly Borel fixed monomial ideals with a given $\Bs$, 
we need some more knowledge of the
structure of such ideals. Assume that $\Bs=(a;n_1\ddd n_a;\seq{n}{a+1}{a+b})$ 
satisfies \eqref{eq101}. Let $\bg_l\ (\ee{0}{l}{a-1})$, $t_l,u_l,\bg_{t_lu_l}\ (\ee{1}{l}{b})$
be integers such that 
$n_l=a-l+\bg_{a-l}\ (\ee{1}{l}{a})$, $n_{a+l}=t_l+u_l+\bg_{t_lu_l}\ (\ee{1}{l}{b})$,
$1\leq\bg_{a-1}<\bg_{a-2}<\cdots<\bg_0$,
$t_l<a,\ u_l<\bg_{t_l},\ \bg_{t_lu_l}>0$ for $\ee{1}{l}{b}$,
and $(t_l,u_l)\neq(t_{l'},u_{l'})$ for $l,l'$ with $l\neq l'$.
As explained in \Rem{\ref{proc40}}, we can treat our problem in this setting.
Let
\begin{equation*}
J:=Rx_1^a+\smm{l=1}{a}\kthree x_1^{a-l}x_2^{\bg_{a-l}}+
\smm{l=1}{b}\ktwo x_1^{t_l}x_2^{u_l}x_3^{\bg_{t_lu_l}}.
\end{equation*}
The right hand side is a $\ktwo$-submodule of $R$ and in fact the 
sum is direct over $\ktwo$. Notice that $J$ is a strongly Borel fixed monomial 
ideal in $R$ if and only if so is $J\cap k[x_1,x_2,x_3]$ in $k[x_1,x_2,x_3]$.
For nonnegative integers $n$ and $\gg$ with $n\geq\gg$, let
\begin{equation*}
\eg(n,\gg):=\max\{\ u\ |\ \ee{0}{u}{n-\gg}\ \andt\ 
x_1^{n-\gg-u}x_2^ux_3^\gg\in \gradn{J}{n}\ \},
\end{equation*}
where $\gradn{J}{n}$ denotes as usual the subsets of $J$ consisting of 
homogeneous elements of degree $n$. 
In case the set on the right hand side is empty, we put $\eg(n,\gg):=-\infty$
for convenience sake.
It is easy to verify that $J$ is a strongly Borel fixed monomial ideal in $R$ if and only if
\begin{align}
{}&x_1^{n-\gg-u}x_2^ux_3^\gg\in \gradn{J}{n}\ \fallt\ u,n,\gg\ 
\witht\ 0\leq u\leq\eg(n,\gg),\ \ee{0}{\gg}{n},\ \andt
\lb{eq234}\\
{}&\eg(n-1,\gg-1)\leq\eg(n,\gg)\leq\eg(n,\gg-1)-1\ \fallt\ n,\gg\ 
\witht\ \ee{1}{\gg}{n}.\
\lb{eq235}
\end{align}
Rewriting the first one of the above conditions in terms of the integers 
$t_l,u_l,\bg_{t_lu_l}\ (\ee{1}{l}{b})$, we obtain the following necessary and
sufficient condition for $J$ to be a strongly Borel fixed ideal.

\begin{lem}\lb{proc113}
The set $J$ is a strongly Borel fixed monomial ideal in $R$ if and only if 
the condition \eqref{eq235} holds and
\begin{equation}
\begin{split}
&\{\ (t_l,u_l,\bg_{t_lu_l})\ |\ n_{a+l}=n\ \andt\ 
\bg_{t_lu_l}=\gg\ (\ee{1}{l}{b}) \}\\
&\hskip2em=\{\ (n-\gg-u,u,\gg)\ |\ u\geq0\ \andt\ 
\ee{\eg(n-1,\gg-1)+1}{u}{\eg(n,\gg)} \}
\end{split}
\lb{eq236}
\end{equation}
for all $n,\gg$ with $n\geq1,\ \ee{1}{\gg}{n}$.
\end{lem}

\begin{pf}
Let $J':=Rx_1^a+\smm{l=1}{a}\kthree x_1^{a-l}x_2^{\bg_{a-l}}$ and 
$J'':=\smm{l=1}{b}\ktwo x_1^{t_l}x_2^{u_l}x_3^{\bg_{t_lu_l}}$.
Then $J'$ is a strongly Borel fixed monomial ideal in $R$.
We have
\begin{equation*}
\begin{split}
&\gradn{J}{n}=\gradn{J'}{n}\op(kx_3+kx_4)\gradn{J''}{n-1}
\op\ds{\ee{1}{l}{b},\ n_{a+l}=n}k x_1^{t_l}x_2^{u_l}x_3^{\bg_{t_lu_l}},
\\
&\gradn{J}{n-1}=\gradn{J'}{n-1}\op\gradn{J''}{n-1}
\end{split}
\end{equation*}
over $k$.
Assuming that $J$ is a strongly Borel fixed monomial ideal, we show \eqref{eq236}. 
When $n<a$, both sides of \eqref{eq236} are empty. Suppose that $n\geq a$
and that $\gg\geq1$.
Let $x_1^{n-\gg-u}x_2^ux_3^\gg$ be an element of 
$\gradn{J'}{n}\op(kx_3+kx_4)\gradn{J''}{n-1}$.
Then it lies in $x_3\gradn{J}{n-1}$, 
since $\gradn{J''}{n-1}\sset\gradn{J}{n-1}$ and 
$\gradn{J'}{n}\cap k[x_1,x_2,x_3]x_3\sset x_3\gradn{J'}{n-1}$.
Namely, $x_1^{n-\gg-u}x_2^ux_3^{\gg-1}$ lies in $\gradn{J}{n-1}$, so that
$u\leq\eg(n-1,\gg-1)$.
Besides,
$x_1^{n-\gg-\eg(n-1,\gg-1)}x_2^{\eg(n-1,\gg-1)}x_3^{\gg}
=x_3(x_1^{n-\gg-\eg(n-1,\gg-1)}x_2^{\eg(n-1,\gg-1)}x_3^{\gg-1})
\in x_3\gradn{J}{n-1}$ in this case.
Hence 
$\max\{\ u\ |\ \ee{0}{u}{n-\gg}\ \andt\ 
x_1^{n-\gg-u}x_2^ux_3^\gg\in \gradn{J'}{n}\op(kx_3+kx_4)\gradn{J''}{n-1}\ \}
=\eg(n-1,\gg-1)$ unless $\eg(n-1,\gg-1)\neq-\infty$.
The right hand side of \eqref{eq236}
therefore coincides with the left hand side by \eqref{eq234}.
Conversely, assuming that both \eqref{eq235} and \eqref{eq236} hold,
we show \eqref{eq234} by induction on $n$. The case $n<a$ is trivial.
Suppose that $n\geq a$ and that \eqref{eq234} is true for smaller values of 
$n$. When $\gg=0$, we see $x_1^{n-u}x_2^u\in\gradn{J'}{n}$
for all $\ee{0}{u}{\eg(n,0)}$ since $J'$ is a strongly Borel fixed monomial ideal.
Let $\gg$ be an integer with $\ee{1}{\gg}{n}$.
It is enough to consider the case where $\eg(n,\gg)\geq0$.
If $\eg(n-1,\gg-1)=-\infty$, then 
$x_1^{n-\gg-u}x_2^ux_3^\gg\in\gradn{J}{n}$ for all $\ee{0}{u}{\eg(n,\gg)}$
by \eqref{eq236}.
Otherwise, $x_1^{n-\gg-u}x_2^ux_3^{\gg-1}\in \gradn{J}{n-1}$ 
for all $\ee{0}{u}{\eg(n-1,\gg-1)}$ by the 
induction hypothesis, which implies that 
$x_1^{n-\gg-u}x_2^ux_3^{\gg}\in \gradn{J}{n}$ for all such $u$.
Hence we see \eqref{eq234} by \eqref{eq236}.
\end{pf}

Assume that \eqref{eq234} -- \eqref{eq236} are valid.
Let $\dg(n,\gg):=\max\{-1,\eg(n,\gg)\}-\max\{-1,\eg(n-1,\gg-1)\}$ for $\ee{1}{\gg}{n}$
and $\chg(n):=\#\{\ l\ |\ n_{a+l}=n,\ \ee{1}{l}{b}\ \}$ for $n\geq0$. 
Then $\dg(n,\gg)=\{\ (t_l,u_l,\bg_{t_lu_l})\ |\ n_{a+l}=n\ \andt\ 
\bg_{t_lu_l}=\gg\ (\ee{1}{l}{b})\ \}$ by \eqref{eq236}.
One sees
\begin{align}
&\smm{\gg=1}{n}\dg(n,\gg)=\chg(n),
\lb{eq237}\\
&\eg(n,\gg)=\eg(n-1,\gg-1)\ \ift\ \dg(n,\gg)=0, \andt
\lb{eq238}\\
&\eg(n,\gg)=\dg(n,\gg)+\max\{-1,\eg(n-1,\gg-1)\}\ \ift\ \dg(n,\gg)>0.
\lb{eq239}
\end{align}
In particular, $\eg(n,\gg)=\eg(n-1,\gg-1)$ for all $\ee{1}{\gg}{n}$ 
such that $\chg(n)=0$. 
Moreover, the $\eg(n,0)\ (n\geq0)$ are completely determined by $J'$ and are nonnegative for 
$n\geq a$. 

\par
Given $\Bs$ satisfying \eqref{eq101}, the $\bg_{a-l}$ can be obtained 
by the simple formula $n_l=a-l+\bg_{a-l}$ for all $\ee{1}{l}{a}$.
The above observation indicates that we can find $t_l,u_l,\bg_{t_lu_l}\ (\ee{1}{l}{b})$
with the help of a computer in the following manner.
\begin{enumr}
\item\lb{cond104}
Find all possible nonnegative integers 
$\dg(n,\gg)$ with $\ee{1}{\gg}{n}$ satisfying \eqref{eq237} for all
$n\geq1$. In fact, we have only to consider the finite cases where $\chg(n)>0$. We put $\dg(n,\gg)=0$
for all $\ee{1}{\gg}{n}$ with $\chg(n)=0$.
\item\lb{cond102}
Define $\eg(n,\gg)$ by \eqref{eq238} and \eqref{eq239} for all $\ee{0}{\gg}{n}$
inductively on $\gg$, starting with $\eg(n,0)\ (n\geq0)$.
\item\lb{cond103}
Check if $\eg(n,\gg)$ satisfies \eqref{eq235} for all $\ee{1}{\gg}{n}\leq n_{a+b}$.
\item
Define $t_l,u_l,\bg_{t_lu_l}\ (\ee{1}{l}{b})$ by \eqref{eq236}.
\end{enumr}
Since $\eg(n,\gg)=\eg(n-1,\gg-1)$ for all $n,\gg$ with $\ee{1}{\gg}{n},\ n>n_{a+b}$
by our constructions \eqref{cond104} and \eqref{cond102}, 
it is enogh to check if $\eg(n,\gg)$ sagisfies \eqref{eq235} only for $\ee{1}{\gg}{n}\leq n_{a+b}$.
Once the integers 
$\bg_l\ (\ee{0}{l}{a-1})$, $t_l,u_l,\bg_{t_lu_l}\ (\ee{1}{l}{b})$ are obtained,
it is an easy matter to verify if $J$ stisfies Cook's assertions \eqref{cond19} and \eqref{cond20}.

\par
The programs we have used to find all $\Bs$'s satisfying
A- or C-conditions are written in C in which only numerical computations
of integers are carried out. The sequences we get by this program
do not necessarily correspond to integral curves in $\Pthree$.

\par
At present, the numbers of the outputs seem too big.
That is, we are far from the answer to the problem stated in Introduction.
The tables below help us to see how A-conditions and C-conditions 
differ from each  other.

\par
Tables 1 and 2 show the numbers of the basic sequences
which satisfy A-conditions or C-conditions or both.
If $g$ is close to the upper bound $1+d(d-3)/6$ of the 
genus of a nonsingular irreducible curve in $\Pthree$ not
contained in any quadric surfaces (see \cite[Introduction]{GP2}),
A-conditions seem stronger than C-conditions.
In general, however, both of the implications 
\lqq A-conditions $\Rightarrow$ C-conditions\rqq
and \lqq C-conditions $\Rightarrow$ A-conditions\rqq\ are false.

\par
Table 3 shows the basic sequences with $(d,g) = (10,9)$
that satisfy C-conditions and the corresponding Borel fixed
monomial ideals, where a triplet $(s,t,u)$ indicates the 
monomial $x_1^sx_2^tx_3^u$.
Among them, only the four sequences
\begin{equation*}
(3;4,5,5;7),\ (3;5,5,5;5,6),\ (4;4,4,4,5;6),\ (4;4,4,5,5;5,5)
\end{equation*}
satisfy A-conditions. Other ones are ruled out by
\Prop{\ref{proc12}} and \Cor{\ref{proc119}}.

\newpage

{\footnotesize

\newcommand{\tableone}{
\begin{tabular}[t]{cccc}
\multicolumn{4}{l}{Table 1: Numbers of outputs\vspace{3pt}}\\
\hline\hline
$(d,g)$&C&A&both\\ 
\hline
(3,0)&1&1&1\\ 
(4,0)&1&1&1\\ 
(5,0)&3&2&2\\ 
(6,0)&4&3&3\\ 
(7,0)&9&6&6\\ 
(8,0)&20&13&12\\ 
(9,0)&48&28&27\\ 
(10,0)&111&69&68\\ 
(11,0)&250&142&137\\ 
(12,0)&570&348&311\\ 
(13,0)&1380&804&731\\ 
\hline\hline
\end{tabular}
}

\newcommand{\tabletwo}{
\smash{\begin{tabular}[t]{cccc}
\multicolumn{4}{l}{\footnotesize
Table 2: Numbers of outputs}\vspace{3pt}\\
\hline\hline
$(d,g)$&C&A&both\\ 
\hline
(26,100)&1&1&1\\ 
(26,99)&2&1&1\\ 
(26,98)&3&1&1\\ 
(26,97)&5&2&2\\ 
(26,96)&6&3&3\\ 
\hline
(27,109)&1&1&1\\ 
(27,108)&2&1&1\\ 
(27,107)&2&1&1\\ 
(27,106)&3&2&2\\ 
(27,105)&4&3&3\\ 
\hline
(28,117)&1&1&1\\ 
(28,116)&1&1&1\\ 
(28,115)&3&2&2\\ 
(28,114)&4&2&2\\ 
(28,113)&5&3&3\\ 
\hline
(29,126)&2&2&2\\ 
(29,125)&2&1&1\\ 
(29,124)&3&1&1\\ 
(29,123)&5&2&2\\ 
(29,122)&6&3&3\\ 
\hline
(30,136)&1&1&1\\ 
(30,135)&2&1&1\\ 
(30,134)&2&1&1\\ 
(30,133)&3&2&2\\ 
(30,132)&5&4&4\\
\hline\hline
\end{tabular}}
}

\newcommand{\tablethree}{
\smash{\begin{tabular}[t]{p{140pt}p{80pt}}
\multicolumn{2}{p{200pt}}
{Table 3: C-conditions for $(d,g)=(10,9)$}
\vspace{3pt}\\
\hline\hline
Borel fixed monomial ideals& 
Basic sequences\\
\hline
(3,0,0)&(3;4,5,5;7)\\ 
(2,2,0)(1,4,0)(0,5,0)&\\ 
(1,3,3)&\\ 
\hline
(3,0,0)&(3;4,5,6;6,6)\\ 
(2,2,0)(1,4,0)(0,6,0)&\\ 
(0,5,1)(1,3,2)&\\ 
\hline
(3,0,0)&(3;5,5,5;5,6)\\ 
(2,3,0)(1,4,0)(0,5,0)&\\ 
(2,2,1)(1,3,2)&\\ 
\hline
(3,0,0)&(3;5,5,6;5,5,6)\\ 
(2,3,0)(1,4,0)(0,6,0)&\\ 
(2,2,1)(1,3,1)(0,5,1)&\\ 
\hline
(4,0,0)&(4;4,4,4,5;6)\\ 
(3,1,0)(2,2,0)(1,3,0)(0,5,0)&\\ 
(0,4,2)&\\ 
&\\
(4,0,0)&\\ 
(3,1,0)(2,2,0)(1,3,0)(0,5,0)&\\ 
(3,0,3)&\\
\hline
(4,0,0)&(4;4,4,5,5;4,6)\\ 
(3,1,0)(2,2,0)(1,4,0)(0,5,0)&\\ 
(3,0,1)(1,3,2)&\\
\hline
(4,0,0)&(4;4,4,5,5;5,5)\\ 
(3,1,0)(2,2,0)(1,4,0)(0,5,0)&\\ 
(1,3,1)(0,4,1)&\\ 
&\\
(4,0,0)&\\ 
(3,1,0)(2,2,0)(1,4,0)(0,5,0)&\\ 
(1,3,1)(3,0,2)&\\
\hline
(4,0,0)&(4;4,4,5,6;4,5,6)\\ 
(3,1,0)(2,2,0)(1,4,0)(0,6,0)&\\ 
(3,0,1)(1,3,1)(0,5,1)&\\
\hline
(4,0,0)&(4;4,5,5,5;4,5,5)\\ 
(3,1,0)(2,3,0)(1,4,0)(0,5,0)&\\ 
(3,0,1)(2,2,1)(1,3,1)&\\ 
\hline\hline
\end{tabular}}
}

\begin{array}[t]{ll}
\tableone&\hspace{30pt}\tablethree
\vspace{30pt}\\
\tabletwo&
\end{array}
}

\newpage

\def\same{{$\rule{15mm}{0.4pt}$\,}}

\def\AUntitled{{\it Untitled}, in preparation.}

\def\Aone{{\it Preparatory structure theorem for ideals defining 
space curves},  Publ. RIMS, Kyoto Univ. {\bf 19} (1983), 493 -- 518.}

\def\Atwo{{\it On the structure of arithmetically Buchsbaum curves in \
$\Pbf^3_k$}, Publ. RIMS, Kyoto Univ. {\bf 20} (1984), 793 -- 837.}

\def\Athree{{\it Examples of nonsingular irreducible curves which give
reducible singular points of \ $\hbox{\rm red(H}_{d,g})$}, Publ. RIMS, Kyoto
Univ. {\bf 21} (1985), 761 -- 786.}

\def\Afour{{\it Curves in  $\Pbf^3$  whose ideals are simple in a
certain numerical sense}, Publ. RIMS, Kyoto Univ. 
{\bf 23} (1987), 1017 -- 1052.}

\def\Afive{{\it Integral arithmetically Buchsbaum curves in \
$\Pbf^3$},  J. Math. Soc. Japan {\bf 41}, No. 1 (1989), 1 -- 8.}

\def\Asix{{\it Application of the generalized Weierstrass preparation 
theorem to the study of homogeneous ideals}, 
Trans. AMS {\bf 317} (1990), 1 -- 43.}

\def\Aseven{{\it Maximal quasi-Buchsbaum graded modules over 
polynomial rings with  $\#\{i|\Hm[i,M] \neq 0, \ i < \dim[R]\} \leq 2$}, 
in Japanese, Proc. 12th Sympos. Commutative Algebra, Kyoto, October 30
-- November 2, 1990,  pp. 110 -- 125.}

\def\Aeight{{\it Free complexes defining maximal \qBbm\ graded modules 
over polynomial rings}, J. Math. Kyoto Univ. {\bf 33}, 
No. 1 (1993), 143 -- 170.}

\def\Anine{{\it On the classification of homogeneous ideals of height two 
in polynomial rings}, Proc.  36th Sympos. Algebra, 
Okayama, July 29 -- August 1, 1991, pp. 129 -- 151.}

\def\Aten{{\it On the basic sequences of homogeneous ideals of height 
three}, Proc. 15th Sympos. Commutative Algebra,
Kashikojima, October 19 -- 22, 1993, pp. 196 -- 207.}

\def\Aeleven{{\it Generators of graded modules associated with
linear filter-regular sequences},
J. Pure Appl. Algebra {\bf 114} (1996), 1 -- 23.}

\def\Atwelve{{\it Basic sequences of homogeneous ideals in 
polynomial rings},  J. Algebra {\bf 190} (1997), 329 -- 360.}

\def\Athirteen{{\it Cohen-Macaulay normal local domains whose
associated graded rings have no depth},  
J. Math. Soc. Japan {\bf 39}, No. 1 (1987), 27 -- 31.}

\def\Afourteen{{\it Basic sequence and Nollet's $\thg_X$ of a
homogeneous ideal of height two}, preprint (August, 1996).}

\def\Afifteen{{\it Existence of homogeneous ideals fitting into long \Bou\ 
sequences},  Proc. Amer. Math. Soc. {\bf 127} (1999), 3461 -- 3466.}

\def\Asixteen{{\it Maximal \Bbm\ modules over \Gor\ local rings},
preprint (October, 1997).}

\def\Aseventeen{{\it Generic \Grob\ bases and \Wei\ bases of homogeneous
submodules of graded free modules}, 
J. Pure Appl. Algebra {\bf 152} (2000), 3 -- 16.}

\def\Aeighteen{{\it Homogeneous prime ideals and graded modules
fitting into long \Bou\ sequences}, J. Pure Appl. Algebra {\bf 162} (2001),
1 -- 21.}

\def\Anineteen{{\it On the basic sequence of a curve in $\Pthree$}, 
in Japanese, Kohkyuhroku 621, Research Institute for Mathematical
Sciences, Kyoto University,  April 1987,  pp. 73 -- 90.}

\def\Atwenty{{\it Inequalities satisfied by the basic sequence
of an integral curve in $\Pthree$}, preprint (November, 2002)}

\def\Atwentyone{{\it Verification of the connectedness
of space curve invariants for a special case}, Comm. Alg. {\bf 32}, No. 10
(2004), 3739 -- 3744.}

\def\BStill{{\it A criterion for detecting
$m$-regularity}, Invent. math. {\bf 87} (1987), 1 -- 11.}

\def\Cook{\textit{The connectedness of space curves invariants}, 
Compositio Math. {\bf 111} (1998), 221 -- 244.}

\def\GPesk{{\it Genre des courbes de l'\'espace
projectif}, in \lq\lq Algebraic Geometry", Lecture Notes in Math. {\bf
687}, Springer-Verlag, Berlin $\cdot$ Heidelberg $\cdot$ New York, 1978,
pp. 31 -- 59.}

\def\GPesktwo{{\it Genre des courbes de l'\'espace
projectif (II)}, Ann. scient. Ec. Norm. Sup., $4^{\textrm{e}}$ s\`erie,
{\bf 15} (1982), 401 -- 418.}

\def\DSch{{\it Non-general type surfaces in $\Pbf^4$: Some remarks on
bounds and constructions},  J. Symb. Comp. {\bf 29} (2000), 545 -- 583.}

\def\Urab{{\it On Hironaka's Monoideal}, Publ. RIMS, Kyoto Univ. {\bf
15} (1979), 279 -- 287.}

\end{document}